\documentclass[11pt,final]{siamltex}
\usepackage{amsmath}
\usepackage{algorithm}
\usepackage{algpseudocode}
\usepackage{bm}
\usepackage{amssymb,version}
\usepackage{cases}
\usepackage{color}
\usepackage{verbatim}
\usepackage{graphicx}
\newtheorem{rem}{Remark}[section]
\allowdisplaybreaks

\usepackage{hyperref}

\begin{document}
	\title{A Unified Framework for First- and Second-Order Energy-Stable Parametric Finite Element Methods for Planar Willmore Flow
		\thanks{This work is supported in part by the National Natural Science Foundation of China (Grant Nos. 12671476, 12271302, 12131014) and Shandong Provincial Natural Science Foundation for Outstanding Youth Scholar (Grant No. ZR2024JQ030). }
    }
	
	\author{ Xiaoli Li
		\thanks{School of Mathematics, Shandong University, Jinan, Shandong, 250100, P.R. China. Email: xiaolimath@sdu.edu.cn  }.
		\and Yujing Yan
		\thanks{Corresponding author. School of Mathematics, Shandong University, Jinan, Shandong, 250100, P.R. China. Email: yyjsust@163.com}.    
    }
	\maketitle
	\begin{abstract}
	The construction of structure-preserving schemes for planar Willmore flow that are second-order accurate in time and avoid large coupled nonlinear solves remains challenging. In this work, we develop a unified framework for first- and second-order fully discrete energy-stable parametric finite element schemes for planar Willmore flow with a length-penalization term. The framework is based on a new equivalent reformulation involving a global auxiliary variable. 
		The main novelty lies in the second-order construction, which combines
	midpoint curvature quantities weighted by elementwise arc-length Jacobian
	ratios with a compatible discrete tangent to obtain exact discrete chain
	rules for the bending energy and curve length. 	A further contribution is to show that both schemes satisfy discrete	modified energy dissipation laws that yield the same type of a priori	bounds on the discrete solution as the original energy law.	They also retain Barrett--Garcke--N{\"u}rnberg tangential motion to promote favorable mesh distribution. Each stage reduces to two linear solves with exactly the same coefficient matrix and one scalar nonlinear solve. We rigorously prove unique solvability of the linear subproblems under suitable assumptions and derive computable sufficient conditions for local existence and uniqueness of the scalar solution near unity. Numerical experiments demonstrate the expected temporal accuracy, energy dissipation, good mesh quality and computational efficiency.
	\end{abstract}
	
	\begin{keywords}
		planar Willmore flow, energy stability, second-order time discretization, unique solvability
	\end{keywords}
	\begin{AMS}
		65M60, 65M12, 35K55, 53C44
	\end{AMS}
	\section{Introduction}
	The elastic bending energy of a planar curve, also known as the planar Willmore energy, measures the elastic cost of bending a thin rod. Its origins lie in Galileo’s work on beams and the classical elastica studied by the Bernoullis and Euler \cite{levien2008elastica}. Willmore investigated the corresponding energy for surfaces \cite{willmore1993riemannian}. Related energies arise in membrane mechanics \cite{helfrich1973elastic,bonito2010biomembranes}, surface design and restoration \cite{clarenz2004restoration,xu2006surface}, and phase-field models of strongly anisotropic interfaces \cite{torabi2009phase,chen2018willmore}. These applications motivate accurate and efficient numerical methods that capture the evolving geometry while preserving energy dissipation. 
	
	In this work, we consider the \(L^2\)-gradient flow of the length-penalized bending energy
	$$
	E(\Gamma(t))
	=\frac12\int_{\Gamma(t)}\varkappa^2\,\mathrm ds
	+\lambda|\Gamma(t)|,
	\qquad \lambda\geq0,
	$$
	where \(\Gamma(t)\subset\mathbb R^2\) is a smooth closed curve, \(\mathrm ds\) is the arc-length element, \(|\Gamma(t)|\) denotes the curve length, and \(\lambda\) is a fixed length-penalization parameter. Let \(\boldsymbol X\) parametrize \(\Gamma(t)\), with unit tangent \(\boldsymbol\tau=\partial_s\boldsymbol X\) and unit normal \(\boldsymbol n\). We define the  curvature by \(\partial_s\boldsymbol\tau=\varkappa\boldsymbol n\) and the normal velocity by \(V=\partial_t\boldsymbol X\cdot\boldsymbol n\), where \(\partial_s\) denotes differentiation with respect to arc length. The gradient flow is then given by
	$$V=-\partial_{ss}\varkappa-\frac12\varkappa^3+\lambda\varkappa,$$
	and satisfies the energy dissipation law
	$$\frac{\mathrm d}{\mathrm dt}E(\Gamma(t))
	=-\int_{\Gamma(t)}V^2\,\mathrm ds\leq0.$$
The main numerical difficulties lie in resolving the nonlinear
fourth-order evolution, preserving the gradient-flow structure
and energy dissipation under full discretization, and maintaining
a well-distributed mesh throughout the evolution.
	
	Numerical approaches to Willmore flow include level set methods \cite{droske2004level}, phase-field methods \cite{du2004phase,franken2013phase,bretin2015phasefield,rumpf2025hybrid}, and front-tracking methods \cite{bgn2016spontaneous,bgn2017boundary,dziuk2008computational}. Among the latter, parametric finite element methods (PFEMs) represent the evolving interface \(\Gamma(t)\) explicitly through a moving mesh. Dziuk \cite{dziuk2008computational} developed a PFEM for Willmore flow with semidiscrete energy dissipation, while Deckelnick and Dziuk \cite{deckelnick2009elastic} proposed and analyzed an energy-stable semidiscrete approximation of length-penalized elastic flow. Insufficient control of tangential mesh motion can lead to mesh deterioration. To improve mesh quality, Barrett, Garcke, and N{\"u}rnberg (BGN) developed parametric formulations that implicitly induce tangential	redistribution~\cite{BARRETT2007441,bgn2008willmore}. Subsequent work	retained this tangential motion and established energy stability for
	elastic flow at the spatially semidiscrete level~\cite{bgn2012elastic}. To preserve energy stability under full discretization, a recent breakthrough was achieved by Bao and Li \cite{li2025sinum}, who developed an unconditionally energy-stable, fully discrete scheme for planar Willmore flow using a curvature evolution equation. As reported by Jiang et al.~\cite{jiang2026structure}, insufficient tangential mesh motion can lead to mesh distortion in computations involving complex curves. Several energy-stable methods use tangential motion to maintain mesh quality during the evolution. For mean curvature flow and surface diffusion, Gao et al.~\cite{gao2026mdr} proposed a BGN--MDR scheme that improves mesh quality while preserving energy stability. For planar Willmore flow, Jiang et al.~\cite{jiang2026structure}	developed energy-stable PFEMs that improve mesh quality through	explicit tangential velocity~\cite{kemmochi2025tangential}
	but require coupled nonlinear solves. Garcke et al.~\cite{garcke2025stable} incorporated BGN tangential motion	into two fully discrete energy-stable schemes with favorable mesh properties, one linear and one nonlinear. In the linear scheme, however, the discrete bending energy is evaluated using the arc-length measure from the previous time step. In subsequent work~\cite{gnz2026splitting}, the same authors developed a framework inspired by the arbitrary Lagrangian--Eulerian (ALE) formulation that accommodates spontaneous curvature and open hypersurfaces with boundary.

	The fully discrete planar Willmore flow schemes described above are first-order
	accurate in time. Extending these schemes to higher order in time while
	retaining discrete energy stability remains challenging.  BGN-based PFEMs have been extended through second-order schemes \cite{JIANG2024113220}, BDF discretizations of orders two to four \cite{jiang2024bdf}, and a predictor–corrector method requiring no additional mesh regularization \cite{jiang2025predictor}. Jiang et al. \cite{jiang2026imex} developed second- and third-order IMEX Runge–Kutta schemes whose stages supply the prediction geometry. These methods address curve shortening, mean curvature flow, surface diffusion, and related evolutions. Arbitrary-order structure-preserving methods use continuous Petrov–Galerkin time discretization for mean curvature flow and surface diffusion \cite{zhang2026arbitrary}, and local pullback maps with Runge–Kutta collocation for surface diffusion and volume-preserving mean curvature flow \cite{bao2026pullback}. For planar Willmore flow, Pan and Jia \cite{pan2026lagrange} introduced a Lagrange multiplier and a supplementary energy evolution equation to establish bending-energy dissipation for BDF1 and BDF2. Their implementation requires coupled nonlinear solves.
	
	The objective of this paper is to develop efficient first- and second-order schemes for planar Willmore flow that combine unconditional energy dissipation and good mesh quality while avoiding large coupled nonlinear solves. Inspired by the zero-energy-contribution (ZEC) technique \cite{yang2021zec}, we introduce a global auxiliary variable into the formulation of Garcke et al. \cite{garcke2025stable}. This yields a unified framework for constructing and analyzing both schemes. The first-order method combines backward Euler time discretization with an implicit–explicit treatment of the reformulated equations. The second-order construction poses two main challenges:
	(i) three-level schemes, such as BDF2 and Crank--Nicolson leapfrog,
	do not automatically yield an energy dissipation law between
	consecutive time levels; and	(ii) Crank--Nicolson discretization involves a unit normal and	an arc-length derivative that depend on the unknown midpoint
	curve, introducing nonlinear geometric coupling.  A predicted midpoint geometry reduces this coupling, but the energies at consecutive time levels remain defined with different arc-length measures. Recovering their difference therefore requires compatible discrete chain rules connecting these geometries while retaining second-order accuracy.	To address these challenges, we adopt a predictor–corrector Crank–Nicolson strategy and introduce weighted midpoint curvature quantities that account for the arc-length measures at both time levels when evaluating curvature terms on the predicted midpoint geometry. Together with a compatible discrete tangent, these quantities yield exact discrete chain rules for the bending and length contributions while retaining second-order midpoint consistency. 	The resulting schemes satisfy unconditional energy dissipation laws and retain BGN tangential mesh redistribution. The main contributions are summarized as follows:
	
	\begin{itemize}
		\item  We construct first- and second-order PFEMs within a unified global auxiliary framework. The second-order scheme combines midpoint curvature quantities	weighted by elementwise arc-length Jacobian ratios with a compatible	discrete tangent to obtain exact discrete chain rules for the bending	energy and curve length while preserving second-order consistency in time. Both schemes retain BGN tangential motion to promote mesh redistribution.
		
		\item We prove an unconditional modified-energy dissipation inequality for the
		first-order scheme and an exact modified-energy dissipation identity for
		the second-order corrector, both of which yield the same type of a priori
		bounds on the discrete solution as the original energy dissipation law. Each stage requires two linear solves with the same coefficient matrix
		and one scalar nonlinear solve. 
		
		\item The second-order midpoint consistency of the weighted geometric
quantities is established. We also prove unique solvability of the linear subproblems under
		suitable assumptions and derive computable sufficient conditions for the
		scalar equation to have a unique solution near unity. Numerical experiments
		confirm that these conditions hold in the cases tested.
	\end{itemize}
	
	To the best of our knowledge, this is the first fully discrete parametric finite element framework for planar Willmore flow that unifies first- and second-order schemes within a global auxiliary formulation, with unconditional modified-energy dissipation and  local existence and uniqueness of the scalar solution. These properties are achieved within a unified discretization and solution framework.
	
	The rest of the paper is organized as follows. In Section \ref{part2}, we present the reformulation. In Section \ref{part3}, we introduce the numerical schemes. In Section \ref{part4}, we establish energy stability, describe the solution procedures, and analyze solvability. Numerical results and conclusions are presented in Sections \ref{part5} and \ref{part6}, respectively. Explicit representations of the first- and second-order
	auxiliary equations are provided in	Appendices~\ref{app:quartic_coefficients}
	and~\ref{app:second_auxiliary}, respectively.
	
	\section{Mathematical formulations}\label{part2}
	
	We first recall the curvature evolution formulation of \cite{garcke2025stable}, which provides the basis for our global auxiliary reformulation.	Let $\Gamma(t)$ be a smooth closed curve with a regular parametrization $\boldsymbol X:\mathbb I\times[0,T]\to\mathbb R^2$, where $\mathbb I=\mathbb R/\mathbb Z$ and $|\boldsymbol X_\rho|>0$. We write
	\[
	\mathrm ds=|\boldsymbol X_\rho|\,\mathrm d\rho,\quad
	\partial_s=|\boldsymbol X_\rho|^{-1}\partial_\rho,\quad
	\boldsymbol\tau=\partial_s\boldsymbol X,\quad
	\boldsymbol n=-\boldsymbol\tau^\perp,
	\]
where the subscript $\rho$ denotes differentiation with respect
to the reference parameter $\rho$, and $(\cdot)^\perp$ denotes
clockwise rotation by $\pi/2$.
	
	We retain two curvature variables: $\varkappa$ is governed by a curvature evolution equation, while $\kappa$ is determined by the geometric identity $\kappa\boldsymbol n=\partial_{ss}\boldsymbol X$. The two curvature variables coincide at the continuous level	but are approximated differently at the discrete level:
	$\varkappa$ enters the bending energy, while $\kappa$ appears
	in the BGN curvature relation. We consider the following system:
	\begin{align}
		V&=-\partial_{ss}\varkappa-\frac12\varkappa^3+\lambda\kappa,\label{le1}\\
		\partial_t\varkappa&=\partial_{ss}V+V\varkappa^2+(\partial_t\boldsymbol X\cdot\boldsymbol\tau)\partial_s\varkappa,\label{le2}\\
		\partial_t\boldsymbol X\cdot\boldsymbol n&=V,\label{le3}\\
		\kappa\boldsymbol n&=\partial_{ss}\boldsymbol X.\label{le4}
	\end{align}
	The system is completed by prescribing the initial parametrization \(\boldsymbol X(\cdot,0)=\boldsymbol X_0\) and the initial curvature \(\varkappa(\cdot,0)=\varkappa_0\), where \(\varkappa_0\) is the geometric curvature of the initial curve.

	\subsection{Weak formulation and energy dissipation}
	We introduce the spaces
	\[
	L^2(\Gamma(t))
	:=\left\{u:\Gamma(t)\to\mathbb R:
	\int_{\Gamma(t)}|u|^2\,\mathrm ds<\infty\right\}
	\]
	and
	\[
	H^1(\Gamma(t))
	:=\left\{u\in L^2(\Gamma(t)):
	\partial_su\in L^2(\Gamma(t))\right\}.
	\]
	The $L^2$ inner product is
	\[
	(u,v)_{\Gamma(t)}
	:=\int_{\Gamma(t)}uv\,\mathrm ds,
	\qquad u,v\in L^2(\Gamma(t)).
	\]
	The vector-valued spaces $[L^2(\Gamma(t))]^2$ and
	$[H^1(\Gamma(t))]^2$ are defined componentwise.
	
	Testing \eqref{le1}--\eqref{le4} with
	$\varphi,\chi\in H^1(\Gamma(t))$,
	$\xi\in L^2(\Gamma(t))$, and
	$\boldsymbol\eta\in[H^1(\Gamma(t))]^2$, respectively,
	and using integration by parts together with the antisymmetric
	reformulation in \cite[Eq.~(2.17)]{garcke2025stable},
	we obtain the following weak formulation. Given the initial curve $\Gamma_0$ parameterized by $\boldsymbol X_0$, with
	$\boldsymbol X(\cdot,0)=\boldsymbol X_0$,
	$\varkappa(\cdot,0)=\varkappa_0$, for any $t\geq0$, find
	$\bigl(\boldsymbol X(\cdot,t),V(\cdot,t),\varkappa(\cdot,t),\kappa(\cdot,t)\bigr)\in
	[H^1(\Gamma(t))]^2\times H^1(\Gamma(t))\times H^1(\Gamma(t))\times L^2(\Gamma(t))$
	such that
	\begin{align}
		&(V, \varphi )_{\Gamma(t)} = (\partial_s \varkappa, \partial_s \varphi )_{\Gamma(t)} - \frac{1}{2} (\varkappa^3, \varphi )_{\Gamma(t)} + \lambda (\kappa, \varphi)_{\Gamma(t)}, \quad \forall \varphi \in H^1({\Gamma(t)}),\label{sc1_ela_sch1}\\
		&(\partial_t\varkappa, \chi )_{\Gamma(t)} + \frac{1}{2} ( \partial_s(\partial_t\boldsymbol{X}) \cdot \boldsymbol{\tau}, \varkappa \chi)_{\Gamma(t)} = -(\partial_sV, \partial_s\chi )_{\Gamma(t)}+ \frac{1}{2} (V \varkappa^2, \chi)_{\Gamma(t)} \nonumber\\
		&\qquad\qquad\qquad\qquad\qquad\qquad\qquad\qquad + \frac{1}{2} (\partial_t\boldsymbol X \cdot \boldsymbol{\tau}, (\partial_s\varkappa \chi - \varkappa \partial_s\chi))_{\Gamma(t)} \quad \forall \chi \in H^1({\Gamma(t)}),\\
		&( \partial_t\boldsymbol{X} \cdot \boldsymbol{n} , \xi )_{\Gamma(t)} = (V, \xi)_{\Gamma(t)} \quad \forall \xi \in L^2({\Gamma(t)}),\\
		&(\kappa \boldsymbol{n}, \boldsymbol{\eta} )_{\Gamma(t)} = -(\partial_s \boldsymbol{X}, \partial_s \boldsymbol{\eta} )_{\Gamma(t)} \quad \forall \boldsymbol{\eta} \in [H^1({\Gamma(t)})]^2.\label{sc1_ela_sch4}
	\end{align}
	
	The weak formulation satisfies the following  energy dissipation law.
	
	\begin{theorem}
		Let $(\boldsymbol X,V,\varkappa,\kappa)$ be a	weak solution of \eqref{sc1_ela_sch1}--\eqref{sc1_ela_sch4}. Then
		\[
		\frac{\mathrm d}{\mathrm dt}
		\left[
		\frac12(\varkappa,\varkappa)_{\Gamma(t)}
		+\lambda|\Gamma(t)|
		\right]
		=-(V,V)_{\Gamma(t)}\leq0.
		\]
	\end{theorem}
	
	\begin{proof}
		Testing the four equations with $\varphi=V$, $\chi=\varkappa$,
		$\xi=-\lambda\kappa$, and
		$\boldsymbol\eta=\lambda\partial_t\boldsymbol X$, respectively,
		and cancelling the antisymmetric term gives $0={}(V,V)_{\Gamma(t)}
		+(\partial_t\varkappa,\varkappa)_{\Gamma(t)}+\frac12
		(\partial_s\partial_t\boldsymbol X\cdot\boldsymbol\tau,
		\varkappa^2)_{\Gamma(t)}
		+\lambda
		(\boldsymbol\tau,\partial_s\partial_t\boldsymbol X)_{\Gamma(t)}.$
		This yields the desired energy identity. See \cite{garcke2025stable} for details.
	\end{proof}

	\subsection{A new reformulation} 
	To construct energy-stable fully discrete schemes, we introduce a time-dependent global auxiliary variable $Q(t)$ governed by the following ODE:
	\begin{align}
		&\frac{1}{2}\frac{dQ^2}{dt}= \frac{1}{2} ( \partial_s(\partial_t\boldsymbol{X}) \cdot \boldsymbol{\tau}, \varkappa^2)_{\Gamma(t)} (Q-1),\label{ode1}\\
		&Q|_{t=0}=1.\label{ode2}
	\end{align}
	Introducing $Q$ into \eqref{sc1_ela_sch1}--\eqref{sc1_ela_sch4}
	and coupling the resulting system with \eqref{ode1}--\eqref{ode2}, we obtain
	the following reformulation:
	\begin{align}
		&(V, \varphi )_{\Gamma(t)} = (\partial_s \varkappa, \partial_s \varphi )_{\Gamma(t)} - \frac{1}{2} (\varkappa^3, \varphi )_{\Gamma(t)} + \lambda (\kappa, \varphi)_{\Gamma(t)}, \quad \forall \varphi \in H^1({\Gamma(t)}),\label{zec1}\\
		&(\partial_t\varkappa, \chi )_{\Gamma(t)} + Q\frac{1}{2} ( \partial_s(\partial_t\boldsymbol{X}) \cdot \boldsymbol{\tau}, \varkappa \chi)_{\Gamma(t)} = -(\partial_sV, \partial_s\chi )_{\Gamma(t)}+ \frac{1}{2} (V \varkappa^2, \chi)_{\Gamma(t)} \nonumber\\
		&\qquad\qquad\qquad\qquad\qquad\qquad\qquad\qquad + \frac{1}{2} (\partial_t\boldsymbol X \cdot \boldsymbol{\tau}, (\partial_s\varkappa \chi - \varkappa \partial_s\chi))_{\Gamma(t)} \quad \forall \chi \in H^1({\Gamma(t)}),\label{zec2}\\
		&( \partial_t\boldsymbol{X} \cdot \boldsymbol{n} , \xi )_{\Gamma(t)} = (V, \xi)_{\Gamma(t)} \quad \forall \xi \in L^2({\Gamma(t)}),\label{zec3}\\
		&(\kappa \boldsymbol{n}, \boldsymbol{\eta} )_{\Gamma(t)} = -(\partial_s \boldsymbol{X}, \partial_s \boldsymbol{\eta} )_{\Gamma(t)} \quad \forall \boldsymbol{\eta} \in [H^1({\Gamma(t)})]^2,\label{zec4}\\
		&\frac{1}{2}\frac{dQ^2}{dt}= \frac{1}{2} ( \partial_s(\partial_t\boldsymbol{X}) \cdot \boldsymbol{\tau}, \varkappa^2)_{\Gamma(t)} (Q-1),\label{zec5}
	\end{align}
	with the initial conditions
	\begin{equation}\label{zec6}
		\boldsymbol X(\cdot,0)=\boldsymbol X_0,\qquad
		\varkappa(\cdot,0)=\varkappa_0,\qquad
		Q(0)=1.
	\end{equation}
	
	\subsubsection{Equivalence and energy dissipation}
	We next show that the reformulated system is equivalent to the original weak formulation and satisfies the same energy dissipation law.
	
	To establish equivalence between the reformulated system \eqref{zec1}--\eqref{zec6} and the original  weak formulation \eqref{sc1_ela_sch1}--\eqref{sc1_ela_sch4}, we first prove that the auxiliary initial-value problem admits the unique solution $Q\equiv1$.
	
	\begin{lemma}\label{thm:Q_unique}
		Let $T>0$ and assume that
		$C(t)=\frac{1}{2}\left(\partial_s(\partial_t\boldsymbol{X})\cdot\boldsymbol{\tau},\varkappa^2\right)_{\Gamma(t)}\in C([0,T])$.
		Then the initial value problem
		\begin{align}	\label{eq:Q_ivp}
			\frac12\frac{\mathrm d}{\mathrm dt}Q^2(t)&=C(t)(Q(t)-1),\quad t\in[0,T],\\
			Q(0)&=1,\label{qw1}
		\end{align}
		admits a unique solution $Q(t)\equiv1, t\in[0,T]$. 
	\end{lemma}
	\begin{proof}
		We first observe that $Q(t)\equiv1$ is a solution of
		\eqref{eq:Q_ivp}--\eqref{qw1}. Hence, existence follows.
		
		To prove uniqueness, let $Q\in C^1([0,T])$ be any solution
		of \eqref{eq:Q_ivp}--\eqref{qw1}. Then
		\begin{equation}\label{eq:Q_product}
			Q(t)Q'(t)=C(t)(Q(t)-1).
		\end{equation}
		Define
		\[
		Y(t)=(Q(t)-1)e^{Q(t)}.
		\]
		Since $Q\in C^1([0,T])$ and $Q(0)=1$, we have
		$Y\in C^1([0,T])$ and $Y(0)=0$.
		Differentiating $Y$ and using \eqref{eq:Q_product}, we obtain
		\[
		\begin{aligned}
			Y'(t)=Q'(t)e^{Q(t)}	+(Q(t)-1)e^{Q(t)}Q'(t)=e^{Q(t)}Q(t)Q'(t)=C(t)(Q(t)-1)e^{Q(t)}	=C(t)Y(t).
		\end{aligned}
		\]
		Multiplying by the integrating factor
		$\exp(-\int_0^t C(s)\,\mathrm ds)$ gives
		\[
		\frac{\mathrm d}{\mathrm dt}
		\left[
		Y(t)\exp\left(-\int_0^t C(s)\,\mathrm ds\right)
		\right]=0.
		\]
		Therefore,
		\[
		Y(t)=Y(0)\exp\left(\int_0^t C(s)\,\mathrm ds\right)=0,
		\qquad t\in[0,T].
		\]
		Since $e^{Q(t)}>0$, the identity
		$(Q(t)-1)e^{Q(t)}=0$ implies $Q(t)=1$ throughout $[0,T]$.
		This proves uniqueness.
	\end{proof}
	
	Hence, the initial-value problem \eqref{eq:Q_ivp}--\eqref{qw1} admits the unique solution $Q\equiv1$ on $[0,T]$. Consequently, \eqref{zec2} reduces to its original form, and the auxiliary equation \eqref{zec5} holds identically. Since the remaining equations are unchanged, \eqref{zec1}--\eqref{zec6} reduces to the original weak formulation \eqref{sc1_ela_sch1}--\eqref{sc1_ela_sch4}. Conversely, any solution of the original weak formulation extends to the reformulated system by setting $Q\equiv1$. Thus, the two formulations are equivalent.
	
	\medskip    
	The reformulated system \eqref{zec1}--\eqref{zec6} also	satisfies an energy dissipation law. Define the modified energy
	\[
	\mathcal E(t)
	:=\frac12(\varkappa,\varkappa)_{\Gamma(t)}
	+\lambda|\Gamma(t)|+\frac12(Q^2-1).
	\]
	
	\begin{theorem}
		The solutions of the scheme \eqref{zec1}-\eqref{zec6} satisfy
		\[
		\frac{\mathrm d}{\mathrm dt}\mathcal E(t)
		=-(V,V)_{\Gamma(t)}\leq0.
		\]
	\end{theorem}
	
	\begin{proof}
		Testing \eqref{zec1}--\eqref{zec4} with
		$\varphi=V$, $\chi=\varkappa$, $\xi=-\lambda\kappa$,
		and $\boldsymbol\eta=\lambda\partial_t\boldsymbol X$,
		respectively, yields
		\begin{align}
			&(V, V)_{\Gamma(t)} = (\partial_s \varkappa, \partial_s V )_{\Gamma(t)} - \frac{1}{2} (\varkappa^3, V )_{\Gamma(t)} + \lambda (\kappa, V)_{\Gamma(t)}, \label{1zec1}\\
			&(\partial_t\varkappa, \varkappa )_{\Gamma(t)} + \frac{Q}{2} ( \partial_s(\partial_t\boldsymbol{X}) \cdot \boldsymbol{\tau}, \varkappa^2)_{\Gamma(t)} = -(\partial_sV, \partial_s\varkappa )_{\Gamma(t)}+ \frac{1}{2} (V \varkappa^2, \varkappa)_{\Gamma(t)}, \label{1zec2}\\
			&-\lambda( \partial_t\boldsymbol{X} \cdot \boldsymbol{n} , \kappa )_{\Gamma(t)} = -\lambda(V, \kappa)_{\Gamma(t)},\label{1zec3}\\
			&\lambda(\kappa \boldsymbol{n}, \partial_t \boldsymbol{X} )_{\Gamma(t)} = -\lambda(\partial_s \boldsymbol{X}, \partial_s (\partial_t \boldsymbol{X}) )_{\Gamma(t)}.\label{1zecs4}
		\end{align}
		Here, the antisymmetric tangential term vanishes for
		$\chi=\varkappa$. Moreover, \eqref{zec5} gives
		\begin{align}
			\frac{d}{dt}(\frac{1}{2}Q^2)= \frac{Q}{2} ( \partial_s(\partial_t\boldsymbol{X}) \cdot \boldsymbol{\tau}, \varkappa^2)_{\Gamma(t)} -  \frac{1}{2} ( \partial_s(\partial_t\boldsymbol{X}) \cdot \boldsymbol{\tau}, \varkappa ^2)_{\Gamma(t)} .\label{1zec5}
		\end{align}
		Combining \eqref{1zec1}--\eqref{1zec5}, we obtain
		\[
		\begin{aligned}
			\frac{\mathrm d}{\mathrm dt}\mathcal E(t)
			&=(\partial_t\varkappa,\varkappa)_{\Gamma(t)}
			+\frac12
			\left(\partial_s(\partial_t\boldsymbol X)\cdot\boldsymbol\tau,
			\varkappa^2\right)_{\Gamma(t)}
			+\lambda(\boldsymbol\tau,
			\partial_s(\partial_t\boldsymbol X))_{\Gamma(t)}
			+\frac12\frac{\mathrm d}{\mathrm dt}(Q^2-1)\\
			&=-(V,V)_{\Gamma(t)}\leq0.
		\end{aligned}
		\]
		This completes the proof.
	\end{proof}
	
	Since $Q\equiv1$, the modified energy coincides with the original energy at the continuous level.
	
	
	\section{First- and second-order fully discrete PFEM schemes}\label{part3}
	In this section, we develop first- and second-order fully discrete PFEM schemes for the reformulated system \eqref{zec1}--\eqref{zec6}. Both schemes satisfy unconditional dissipation laws for a modified discrete	energy, yielding the same type of a priori bounds on the discrete solution as the original energy dissipation law. At each stage, the discrete equations reduce to linear subproblems and a scalar nonlinear equation.
	
	\subsection{Notation and discretization}
	We now introduce the spatial and temporal discretizations used in the fully discrete schemes. Let
	$\mathbb I=[0,1]=\bigcup_{j=1}^{N}I_j$, $N\geq 3$,
	be a partition of $\mathbb I$ with $I_j=[\rho_{j-1},\rho_j]$. We define the piecewise linear finite element space
	$$V^h:=\left\{u\in C(\mathbb I):\left.u\right|_{I_j}\text{ is linear},\ j=1,\ldots,N;\ u(\rho_0)=u(\rho_N)\right\}\subseteq H^1(\mathbb I).$$
	
For a polygonal curve \(\Gamma^h=\Gamma[\boldsymbol X^h]\) parameterized by \(\boldsymbol X^h\in[V^h]^2\), we define the mass-lumped inner product by the composite trapezoidal rule:
	\begin{align}
		(u,v)^h_{\Gamma^h}
		:=
		\frac12
		\sum_{j=1}^{N}
		\left|\boldsymbol X^h(\rho_j,t)-\boldsymbol X^h(\rho_{j-1},t)\right|
		\left[(u\cdot v)(\rho_j^-)+(u\cdot v)(\rho_{j-1}^+)\right].
	\end{align}
	Here, $u$ and $v$ are scalar- or vector-valued piecewise continuous functions with possible jumps at the nodes $\{\rho_j\}_{j=1}^{N}$, and
	$u(\rho_j^\pm):=\lim_{\rho\to\rho_j^\pm}u(\rho)$.
	
	Let $\Delta t>0$ be a uniform time step and define
	$t_m:=m\Delta t$, $m=0,\ldots,M$, $M\Delta t=T$.
	At each time level $t_m$, let $\Gamma^m:=\Gamma[\boldsymbol X^m]$, where
	$\boldsymbol X^m\in[V^h]^2$,
	$V^m,\varkappa^m,\kappa^m\in V^h$,
	$Q^m\in\mathbb R$,
	denote the fully discrete approximations of
	$\boldsymbol X(\cdot,t_m)$, $V(\cdot,t_m)$, $\varkappa(\cdot,t_m)$, $\kappa(\cdot,t_m)$, and $Q(t_m)$, respectively. The discrete curve $\Gamma^m$ is a polygonal curve consisting of $N$ ordered line segments
	\begin{equation}
		\boldsymbol h_j^m:=\boldsymbol X^m(\rho_j)-\boldsymbol X^m(\rho_{j-1}),\qquad j=1,\ldots,N,
		\label{eq:discrete_edges}
	\end{equation}
	where periodic indexing is used. We assume throughout that
	$|\boldsymbol h_j^m|>0$, $j=1,\ldots,N$,
	so that $\Gamma^m$ is nondegenerate.

	The discrete arc-length derivative $\partial_s$, the unit tangent $\boldsymbol\tau^m$, and the unit normal $\boldsymbol n^m$ on $\Gamma^m$ are defined elementwise by
	\begin{equation}
		\left.\partial_s u^h\right|_{I_j}
		:=
		\frac{u^h(\rho_j)-u^h(\rho_{j-1})}{|\boldsymbol h_j^m|},
		\qquad
		\boldsymbol\tau_j^m:=\frac{\boldsymbol h_j^m}{|\boldsymbol h_j^m|},
		\qquad
		\boldsymbol n_j^m:=-(\boldsymbol\tau_j^m)^\perp.
	\end{equation}
	Consequently, for piecewise continuous scalar- or vector-valued functions $u$ and $v$, the mass-lumped inner product on $\Gamma^m$ is defined by
	\begin{equation}
		(u^h,v^h)^h_{\Gamma^m}
		:=
		\frac12\sum_{j=1}^{N}|\boldsymbol h_j^m|
		\left[(u^h\cdot v^h)(\rho_{j-1}^{+})+(u^h\cdot v^h)(\rho_j^{-})\right],
		\qquad
		\forall\,u^h,v^h\in V^h\ \text{or}\ [V^h]^2.
		\label{eq:mass_lumped_inner_product}
	\end{equation}
	
	With the above notation, we introduce two fully discrete PFEM schemes for the reformulated system: a first-order scheme and a second-order predictor--corrector scheme. Both are designed to be unconditionally energy stable.
	
	\subsection{First-order energy-stable scheme}
	We first introduce a first-order, fully discrete PFEM scheme for the
	reformulated system based on the backward Euler time discretization.
	The global auxiliary variable $Q$ is discretized simultaneously  to preserve the energy-dissipation law at the fully discrete level. 
	
	Given initial data $\boldsymbol X^0\in[V^h]^2$,
	$\varkappa^0\in V^h$, and $Q^0\in\mathbb R$,
	for each $m=0,\ldots,M-1$, find $ \boldsymbol X^{m+1}\in[V^h]^2,
	V^{m+1},\varkappa^{m+1},\kappa^{m+1}\in V^h,
	Q^{m+1}\in\mathbb R$ 	such that the following equations hold for all
	$\phi^h,\chi^h,\xi^h\in V^h$ and
	$\boldsymbol\eta^h\in[V^h]^2$:
	\begin{align}
		& \left({V}^{m+1} , \phi^{h}\right)_{\Gamma^{m}}^{h} = \left(\partial_s \varkappa^{m+1}, \partial_s\phi^h \right)_{\Gamma^{m}}^{h}  - \frac{1}{2}\left((\varkappa^m)^2 \varkappa^{m+1},\; \phi^h \right)_{\Gamma^{m}}^{h} +  \lambda\left(\kappa^{m+1},\; \phi^h \right)_{\Gamma^{m}}^{h}, \label{ful1} \\
		& \left( \frac{{\varkappa}^{m+1}-{\varkappa}^{m} }{\Delta t} ,\, \chi^{h} \right)_{\Gamma^{m}}^{h} = -\left( \partial_sV^{m+1}, \partial_s\chi^h\right)_{\Gamma^{m}}^{h} + \frac{1}{2} \left( (\varkappa^m)^2 V^{m+1},\; \chi^h \right)_{\Gamma^{m}}^{h} \nonumber\\
		&\quad\qquad\quad \qquad\qquad\qquad+ \frac{1}{2}\left(\boldsymbol{\tau}^m \cdot \frac{\boldsymbol{X}^m - \boldsymbol{X}^{m-1}}{\Delta t},\; \left[\partial_s\varkappa^{m+1} \chi^h - \varkappa^{m+1} \partial_s\chi^h\right]\right)_{\Gamma^{m}}^{h} \nonumber\\ 
		&\quad\qquad\quad \qquad\qquad\qquad -\frac{1}{2}Q^{m+1}\left( \partial_s (\frac{\boldsymbol{X}^m - \boldsymbol{X}^{m-1}}{\Delta t}) \cdot \partial_s \boldsymbol{X}^m \varkappa^{m},\;  \chi^h \right)_{\Gamma^{m}}^{h}, \label{ful2}\\ 
		& (\boldsymbol{n}^m \cdot \frac{\boldsymbol{X}^{m+1}-\boldsymbol{X}^m}{\Delta t} ,\; \xi^h )_{\Gamma^{m}}^{h}=( V^{m+1},\; \xi^h )_{\Gamma^{m}}^{h},\label{ful3}\\ 	
		&\left(\kappa^{m+1} \boldsymbol{n}^m,\; \boldsymbol{\eta}^h \right)_{\Gamma^{m}}^{h} = -\left( \partial_s \boldsymbol{X}^{m+1},\; \partial_s \boldsymbol{\eta}^h \right)_{\Gamma^{m}}^{h},\label{ful4}\\
		&\frac{{Q}^{m+1}-{Q}^m  }{\Delta t}Q^{m+1} =\frac{1}{2}Q^{m+1}\left( \partial_s (\frac{\boldsymbol{X}^m - \boldsymbol{X}^{m-1}}{\Delta t}) \cdot \partial_s \boldsymbol{X}^m \varkappa^{m},\;  \varkappa^{m+1} \right)_{\Gamma^{m}}^{h} \nonumber\\
		&\quad \quad
		\quad \qquad \qquad \qquad \qquad - \frac{1}{2}\left( \partial_s (\frac{\boldsymbol{X}^{m+1} - \boldsymbol{X}^{m}}{\Delta t}) \cdot \partial_s \boldsymbol{X}^{m+1} \varkappa^{m+1},\;  \varkappa^{m+1} \right)_{\Gamma^{m}}^{h} .\label{ful5}
	\end{align}
	For a given $Q^{m+1}$, \eqref{ful1}--\eqref{ful4}
	form a linear system for
	$(\boldsymbol X^{m+1},V^{m+1},
	\varkappa^{m+1},\kappa^{m+1})$, coupled with the scalar
	nonlinear equation \eqref{ful5}.
	The linear decomposition, scalar nonlinear solve, and
	complete implementation procedure are described in Section~\ref{part4}.
	
	\subsection{Second-order energy-stable  scheme}
	A standard Crank--Nicolson discretization involves a midpoint normal and arc-length derivative that depend on the unknown new curve, introducing nonlinear geometric coupling.
	Moreover, the bending energies at consecutive time levels are evaluated with different arc-length measures, so standard midpoint approximations do not directly recover their difference. To address these difficulties, we combine a predicted midpoint geometry with suitably weighted geometric quantities to retain second-order temporal consistency and establish a dissipation law for the modified discrete energy. In the predictor step, we apply the first-order scheme from the previous subsection with
	a half time step $\Delta t/2$ to obtain
	$\widetilde{\boldsymbol X}^{m+\frac12}$ and
	$\widetilde\varkappa^{m+\frac12}$.
	The corresponding predicted midpoint curve is denoted by
	$\widetilde\Gamma^{m+\frac12}:=\Gamma[\widetilde{\boldsymbol X}^{m+\frac12}]$. The unit normal vector and the arc-length derivative associated with the predicted curve $\widetilde\Gamma^{m+\frac12}$ are defined by $\widetilde{\boldsymbol n}^{m+\frac12}:=
	-\left(\frac{\widetilde{\boldsymbol X}^{m+\frac12}_\rho}{|\widetilde{\boldsymbol X}^{m+\frac12}_\rho|}\right)^\perp, \widetilde\partial_s:=
	\frac{1}{|\widetilde{\boldsymbol X}^{m+\frac12}_\rho|}\partial_\rho.$ In the corrector step, a Crank--Nicolson-type discretization is formulated on $\widetilde\Gamma^{m+\frac12}$, together with suitably designed weighted geometric quantities, to achieve second-order temporal consistency and preserve a discrete energy dissipation law.
	
	To formulate the corrector step, we introduce the elementwise ratios of
	the arc-length Jacobians at time levels $m+1$ and $m$ to the corresponding
	Jacobian on the predicted midpoint curve:
	\begin{equation}
		J_+^{m+\frac12}
		:=\frac{|\boldsymbol X_\rho^{m+1}|}
		{|\widetilde{\boldsymbol X}_\rho^{m+\frac12}|},
		\qquad
		J_-^{m+\frac12}
		:=\frac{|\boldsymbol X_\rho^m|}
		{|\widetilde{\boldsymbol X}_\rho^{m+\frac12}|}.
	\end{equation}
	Using these ratios, we define the weighted midpoint curvature
	and the associated weighted time difference by
	\begin{equation}
		\begin{aligned}
			\varkappa_J^{m+\frac12}
			:=\frac12\left(
			\sqrt{J_+^{m+\frac12}}\,\varkappa^{m+1}
			+\sqrt{J_-^{m+\frac12}}\,\varkappa^m
			\right),~
			D_t^J\varkappa^{m+\frac12}
			:=\frac{
				\sqrt{J_+^{m+\frac12}}\,\varkappa^{m+1}
				-\sqrt{J_-^{m+\frac12}}\,\varkappa^m
			}{\Delta t}.
		\end{aligned}
	\end{equation}
	We also define the discrete tangent
	\begin{equation}
		\boldsymbol\tau_J^{m+\frac12}
		:=\frac{
			\boldsymbol X_\rho^{m+1}+\boldsymbol X_\rho^m
		}{
			|\boldsymbol X_\rho^{m+1}|+|\boldsymbol X_\rho^m|
		}.
	\end{equation}
	
	These weighted quantities are designed to provide second-order
	approximations of the corresponding geometric quantities at
	$t_{m+\frac12}$ while retaining the discrete identities required for
	the energy-dissipation law. The predictor--corrector procedure is specified as follows.
	
	\medskip
	\paragraph{Step 1: Predictor}
	Given $\boldsymbol X^m,\boldsymbol X^{m-1}\in[V^h]^2$,
	$\varkappa^m\in V^h$, and $Q^m\in\mathbb R$,
	the first-order scheme is applied with the half time step $\Delta t/2$ to determine
	$\widetilde{\boldsymbol X}^{m+\frac12}\in[V^h]^2$,
	$\widetilde V^{m+\frac12},\widetilde\varkappa^{m+\frac12},\widetilde\kappa^{m+\frac12}\in V^h$, $\widetilde Q^{m+\frac12}\in\mathbb R$,
	such that, for all $(\phi^h,\chi^h,\xi^h,\boldsymbol\eta^h)\in V^h\times V^h\times V^h\times[V^h]^2$,
	the following equations hold:
	\begin{align}
		&(\widetilde V^{m+\frac12},\phi^h)^h_{\Gamma^m}={}(\partial_s\widetilde\varkappa^{m+\frac12},
		\partial_s\phi^h)^h_{\Gamma^m}-\frac12 ((\varkappa^m)^2\widetilde\varkappa^{m+\frac12},
		\phi^h)^h_{\Gamma^m}+\lambda (\widetilde\kappa^{m+\frac12},\phi^h)^h_{\Gamma^m},\label{ss1}\\
		&\left( \frac{\widetilde\varkappa^{m+\frac12}-\varkappa^m}{\Delta t/2},\chi^h\right)^h_{\Gamma^m}={}-(\partial_s\widetilde V^{m+\frac12},\partial_s\chi^h)^h_{\Gamma^m}+\frac12((\varkappa^m)^2\widetilde V^{m+\frac12}, \chi^h)^h_{\Gamma^m}\nonumber\\
		&\qquad\qquad\qquad\qquad\qquad+
		\frac12 \left( \boldsymbol\tau^m\cdot\frac{\boldsymbol X^m-\boldsymbol X^{m-1}}{\Delta t},
		\partial_s\widetilde\varkappa^{m+\frac12}\chi^h-\widetilde\varkappa^{m+\frac12}\partial_s\chi^h\right)^h_{\Gamma^m}\nonumber\\
		&\qquad\qquad\qquad\qquad\qquad-\frac12\widetilde Q^{m+\frac12}\left(\partial_s\frac{\boldsymbol X^m-\boldsymbol X^{m-1}}{\Delta t}\cdot\partial_s\boldsymbol X^m\,\varkappa^m, \chi^h\right)^h_{\Gamma^m},\label{ss2}\\
		&\left(\boldsymbol n^m\cdot \frac{\widetilde{\boldsymbol X}^{m+\frac12}-\boldsymbol X^m}
		{\Delta t/2},\xi^h \right)^h_{\Gamma^m}={}(\widetilde V^{m+\frac12},\xi^h)^h_{\Gamma^m},\label{ss3}\\
		&(\widetilde\kappa^{m+\frac12}\boldsymbol n^m,
		\boldsymbol\eta^h)^h_{\Gamma^m}={}-(\partial_s\widetilde{\boldsymbol X}^{m+\frac12}, \partial_s\boldsymbol\eta^h)^h_{\Gamma^m},\label{ss4}\\
		&\frac{\widetilde Q^{m+\frac12}-Q^m}{\Delta t/2}
		\widetilde Q^{m+\frac12}={}\frac12\widetilde Q^{m+\frac12}
		\left(\partial_s\frac{\boldsymbol X^m-\boldsymbol X^{m-1}}{\Delta t}\cdot\partial_s\boldsymbol X^m\,\varkappa^m, \widetilde\varkappa^{m+\frac12}\right)^h_{\Gamma^m}\nonumber\\
		&\quad\quad\qquad\qquad\qquad-\frac12\left(\partial_s\frac{\widetilde{\boldsymbol X}^{m+\frac12}-\boldsymbol X^m}{\Delta t/2}\cdot\partial_s\widetilde{\boldsymbol X}^{m+\frac12}
		\,\widetilde\varkappa^{m+\frac12},\widetilde\varkappa^{m+\frac12}\right)^h_{\Gamma^m}.\label{ss5}
	\end{align}

	\paragraph{Step 2: Corrector}
	Based on the predicted curve $\widetilde\Gamma^{m+\frac12}$, find
	$\boldsymbol X^{m+1}\in[V^h]^2$, $V^{m+1},\varkappa^{m+1}$, $ \kappa^{m+1}\in V^h$, $Q^{m+1}\in\mathbb R$,
	such that, for all
	$(\phi^h,\chi^h,\xi^h,\boldsymbol\eta^h)\in V^h\times V^h\times V^h\times[V^h]^2$,
	the following Crank--Nicolson-type equations hold:
	\begin{align}
		&(V^{m+\frac12},\phi^h)^h_{\widetilde\Gamma^{m+\frac12}}={}
		(\widetilde\partial_s \varkappa^{m+\frac12},	\widetilde\partial_s\phi^h)^h_{\widetilde\Gamma^{m+\frac12}} -\frac12	((\widetilde\varkappa^{m+\frac12})^2
		\varkappa^{m+\frac12},\phi^h)^h_	{\widetilde\Gamma^{m+\frac12}}	+\lambda	(\kappa^{m+\frac12},\phi^h)^h_	{\widetilde\Gamma^{m+\frac12}},	\label{sec1}\\
		&	(\frac{\varkappa^{m+1}-\varkappa^m}{\Delta t},\chi^h)^h_	{\widetilde\Gamma^{m+\frac12}}
		={}	-(\widetilde\partial_sV^{m+\frac12},	\widetilde\partial_s\chi^h)^h_	{\widetilde\Gamma^{m+\frac12}} +\frac12	((\widetilde\varkappa^{m+\frac12})^2	V^{m+\frac12},\chi^h)^h_	{\widetilde\Gamma^{m+\frac12}} \nonumber\\
		&\quad\qquad\qquad\qquad\qquad\qquad+\frac12
		\left(	\widetilde{\boldsymbol\tau}^{m+\frac12}
		\cdot\widetilde{\boldsymbol U}^{m+\frac12},
		\widetilde\partial_s\varkappa^{m+\frac12}\chi^h
		-\varkappa^{m+\frac12}	\widetilde\partial_s\chi^h
		\right)^h_{\widetilde\Gamma^{m+\frac12}}	\nonumber\\
		&\quad\qquad\qquad\qquad\qquad\qquad-\frac12 Q^{m+\frac12}
		\left(	\widetilde\partial_s\widetilde{\boldsymbol U}^{m+\frac12}	\cdot
		\widetilde\partial_s\widetilde{\boldsymbol X}^{m+\frac12}\,\widetilde\varkappa^{m+\frac12},	\chi^h	\right)^h_{\widetilde\Gamma^{m+\frac12}},\label{sec2}\\
		&(\widetilde{\boldsymbol n}^{m+\frac12}
		\cdot  \frac{\boldsymbol{X}^{m+1}-\boldsymbol{X}^m}{\Delta t}, \xi^h)^h_
		{\widetilde\Gamma^{m+\frac12}}	={}	(V^{m+\frac12},\xi^h)^h_
		{\widetilde\Gamma^{m+\frac12}},	\label{sec3}\\
		&(\kappa^{m+\frac12}\widetilde{\boldsymbol n}^{m+\frac12},
		\boldsymbol\eta^h)^h_{\widetilde\Gamma^{m+\frac12}}	={}	-(\widetilde\partial_s {\boldsymbol X}^{m+\frac12},	\widetilde\partial_s\boldsymbol\eta^h)^h_	{\widetilde\Gamma^{m+\frac12}},\label{sec4}\\
		&	\frac{Q^{m+1}-Q^m}{\Delta t}Q^{m+\frac12}	={}	\frac12
		Q^{m+\frac12}	\left( 	\widetilde\partial_s\widetilde{\boldsymbol U}^{m+\frac12}	\cdot	\widetilde\partial_s\widetilde{\boldsymbol X}^{m+\frac12}
		\,\widetilde\varkappa^{m+\frac12},
		\varkappa^{m+\frac12}	\right)^h_{\widetilde\Gamma^{m+\frac12}}\nonumber\\
		&\qquad\qquad\qquad\qquad +(\frac{\varkappa^{m+1}-\varkappa^m}{\Delta t},\varkappa^{m+\frac{1}{2}})^h_	{\widetilde\Gamma^{m+\frac12}}-	\left(	D_t^J\varkappa^{m+\frac{1}{2}},	\varkappa_J^{m+\frac12}	\right)^h_{\widetilde\Gamma^{m+\frac12}} \label{sec6}\\
		&\qquad\qquad\qquad\qquad -\lambda
		\int_{\mathbb I}	\boldsymbol\tau_{J}^{m+\frac12}
		\cdot	\partial_\rho  \frac{\boldsymbol{X}^{m+1}-\boldsymbol{X}^m}{\Delta t}
		\,d\rho  + \lambda (\widetilde\partial_s {\boldsymbol X}^{m+\frac12},	\widetilde\partial_s  \frac{\boldsymbol{X}^{m+1}-\boldsymbol{X}^m}{\Delta t} )^h_	{\widetilde\Gamma^{m+\frac12}},\nonumber
	\end{align}
	where  $V^{m+\frac12}=\frac{ V^{m+1}+V^{m}}{2}$,
	$\varkappa^{m+\frac12}=\frac{\varkappa^{m+1}+\varkappa^m}{2}$,
	$\kappa^{m+\frac12}=\frac{\kappa^{m+1}+\kappa^m}{2},$
	${\boldsymbol X}^{m+\frac12}=\frac{\boldsymbol X^{m+1}+\boldsymbol X^m}{2}$, $	Q^{m+\frac12}:=\frac12(Q^{m+1}+Q^m)$. For $m\geq2$, $\widetilde{\boldsymbol U}^{m+\frac12}:=  \frac{  	2\boldsymbol X^m
		-3\boldsymbol X^{m-1}	+\boldsymbol X^{m-2} }{\Delta t}.$
	
	The weighted quantities
	$\varkappa_J^{m+\frac12}$,
	$D_t^J\varkappa^{m+\frac{1}{2}}$, and
	$\boldsymbol\tau_J^{m+\frac12}$
	are introduced to reconcile the midpoint discretization with the
	change of the discrete arc-length measure between consecutive time
	levels.  As shown below, these quantities satisfy the required second-order
	midpoint consistency relations and play a key role in establishing
	the unconditional energy stability of the second-order scheme.

	\subsection{Temporal consistency of the weighted quantities}
	We first establish the midpoint consistency of the weighted quantities.

	\begin{lemma}	Let
		$t_*:=t_{m+\frac12}:=\frac{t_m+t_{m+1}}{2}$,
		$g(\rho,t):=|\boldsymbol X_\rho(\rho,t)|$,
		$\boldsymbol\tau(\rho,t):=\frac{\boldsymbol X_\rho(\rho,t)}{|\boldsymbol X_\rho(\rho,t)|}$.
		Assume that $\boldsymbol X_\rho$ and $\varkappa$ are three
		times continuously differentiable in time, with derivatives
		uniformly bounded on $\mathbb I\times[0,T]$, and that
		the curve is regular, i.e., $g(\rho,t)\geq c_0>0$ for all $(\rho,t)\in\mathbb  I\times[0,T]$. Write $f_*:=f(\rho,t_*)$, and define	$\partial_{s_*}:=g_*^{-1}\partial_\rho$ 	and the following continuous-reference quantities: $J_+^*:=\frac{g(t_{m+1})}{g(t_{m+\frac12})}$,
		$J_-^*:=\frac{g(t_m)}{g(t_{m+\frac12})}$, 	$\varkappa_J^*:=\frac12\bigl(\sqrt{J_+^*}\,\varkappa(t_{m+1})+\sqrt{J_-^*}\,\varkappa(t_m)\bigr)$,
		$D_t^{J,*}\varkappa:=\frac{\sqrt{J_+^*}\,\varkappa(t_{m+1})-\sqrt{J_-^*}\,\varkappa(t_m)}{\Delta t}$,	$\boldsymbol\tau_J^*:=\frac{\boldsymbol X_\rho(t_{m+1})+\boldsymbol X_\rho(t_m)}{|\boldsymbol X_\rho(t_{m+1})|+|\boldsymbol X_\rho(t_m)|}$.
		Then, as $\Delta t\to0$,
		\begin{align}
			\varkappa_J^* &= \varkappa_*+O(\Delta t^2),\label{eq:kappaJ_consistency}\\[-2pt]
			D_t^{J,*}\varkappa &= \varkappa_{t,*}+\frac12\bigl(\boldsymbol\tau_*\cdot\partial_{s_*}\boldsymbol X_{t,*}\bigr)\varkappa_*+O(\Delta t^2),\label{eq:DtJ_consistency}\\[-2pt]
			\boldsymbol\tau_J^* &= \boldsymbol\tau_*+O(\Delta t^2).\label{eq:tauJ_consistency}
		\end{align}
	\end{lemma}
	
	\begin{proof}
		Let $g_*:=g(\rho,t_*)$. Taylor expansion in time about $t_*$ gives $	g(t_{m+1}) = g_*+\tfrac{\Delta t}{2}g_{t,*}+\tfrac{\Delta t^2}{8}g_{tt,*}+O(\Delta t^3)$, $g(t_m)  = g_*-\tfrac{\Delta t}{2}g_{t,*}+\tfrac{\Delta t^2}{8}g_{tt,*}+O(\Delta t^3)$. 
		It follows that $J_\pm^*=1\pm\frac{\Delta t}{2}\frac{g_{t,*}}{g_*}+O(\Delta t^2)$.
		Since $g_t=g\,\boldsymbol\tau\cdot\partial_s\boldsymbol X_t$,
		we have $\frac{g_{t,*}}{g_*}=\boldsymbol\tau_*\cdot\partial_{s_*}\boldsymbol X_{t,*}$, which yields $J_\pm^* = 1\pm\frac{\Delta t}{2}\bigl(\boldsymbol\tau_*\cdot\partial_{s_*}\boldsymbol X_{t,*}\bigr)+O(\Delta t^2)$. Using $(1+z)^{1/2}=1+\frac12z+O(z^2)$ yields $\sqrt{J_\pm^*} = 1\pm\frac{\Delta t}{4}\bigl(\boldsymbol\tau_*\cdot\partial_{s_*}\boldsymbol X_{t,*}\bigr)+O(\Delta t^2)$.
		
		Next, define $F(t):=\sqrt{\frac{g(t)}{g_*}}\,\varkappa(t)$. By the definitions of $J_\pm^*$,
		$F(t_{m+1})=\sqrt{J_+^*}\,\varkappa(t_{m+1})$,\quad
		$F(t_m)=\sqrt{J_-^*}\,\varkappa(t_m)$.
		Hence,
		$\varkappa_J^*=\frac12\left[F(t_{m+1})+F(t_m)\right]$,\quad
		$D_t^{J,*}\varkappa=\frac{F(t_{m+1})-F(t_m)}{\Delta t}$.
		Since $t_{m+1}=t_*+\frac{\Delta t}{2}$,\ $t_m=t_*-\frac{\Delta t}{2}$, the standard centered-average and centered-difference expansions yield
		$\varkappa_J^*=F(t_*)+O(\Delta t^2)=\varkappa_*+O(\Delta t^2)$,
		which proves \eqref{eq:kappaJ_consistency}, and
		$D_t^{J,*}\varkappa=F_t(t_*)+O(\Delta t^2)$.
		A direct differentiation gives
		$F_t(t_*)=\varkappa_{t,*}+\frac12\frac{g_{t,*}}{g_*}\varkappa_*$.
		Using
		$\frac{g_{t,*}}{g_*}=\boldsymbol\tau_*\cdot\partial_{s_*}\boldsymbol X_{t,*}$,
		we obtain
		$D_t^{J,*}\varkappa=\varkappa_{t,*}+\frac12\bigl(\boldsymbol\tau_*\cdot\partial_{s_*}\boldsymbol X_{t,*}\bigr)\varkappa_*+O(\Delta t^2)$,
		which proves \eqref{eq:DtJ_consistency}.
		
		Finally, Taylor expansion of $\boldsymbol X_\rho$ about $t_*$ yields
		$\boldsymbol X_\rho(t_{m+1})+\boldsymbol X_\rho(t_m)=2\boldsymbol X_{\rho,*}+O(\Delta t^2)$,
		while
		$|\boldsymbol X_\rho(t_{m+1})|+|\boldsymbol X_\rho(t_m)|=2g_*+O(\Delta t^2)$.
		Therefore,
		$\boldsymbol\tau_J^*=\frac{2\boldsymbol X_{\rho,*}+O(\Delta t^2)}{2g_*+O(\Delta t^2)}=\frac{\boldsymbol X_{\rho,*}}{g_*}+O(\Delta t^2)=\boldsymbol\tau_*+O(\Delta t^2)$,	proving \eqref{eq:tauJ_consistency}.
	\end{proof}
	
	\medskip
	\begin{rem}
		The auxiliary equation \eqref{sec6} is a second-order consistent temporal approximation of the continuous auxiliary equation at $t_{m+\frac12}$. Indeed, evaluate the discrete temporal quantities using a sufficiently smooth continuous solution at $t_m$ and $t_{m+1}$. Using the midpoint consistency relations \eqref{eq:kappaJ_consistency}--\eqref{eq:tauJ_consistency}
		together with the standard centered expansions, the auxiliary equation yields
		\begin{equation}
			Q_*Q_{t,*}
			=
			\tfrac12(Q_*-1)\left((\boldsymbol\tau_*\cdot\partial_{s_*}\boldsymbol X_{t,*})\varkappa_*,\varkappa_*\right)_{\Gamma_*}+O(\Delta t^2).
			\label{eq:Q_midpoint_consistency}
		\end{equation}
		Here the two terms involving $\lambda$ cancel at the continuous level owing to the identity
		$(\partial_{s_*}\boldsymbol X_*,\partial_{s_*}\boldsymbol X_{t,*})_{\Gamma_*}\\=\int_{\mathbb I}\boldsymbol\tau_*\cdot \boldsymbol X_{t\rho,*}\,d\rho$.
		Hence the continuous counterpart of the discrete auxiliary equation is
		$Q Q_t=\tfrac12(Q-1)\left((\boldsymbol\tau\cdot\partial_s\boldsymbol X_t)\varkappa,\varkappa\right)_{\Gamma(t)}$,
		which coincides with the continuous auxiliary equation \eqref{eq:Q_ivp} introduced in the reformulated system. 
	\end{rem}
	
	\section{Analysis and implementation of the fully discrete schemes}\label{part4}
	This section presents the stability analysis, implementation, and solvability of the proposed fully discrete schemes. We first establish unconditional stability estimates for a modified discrete energy. We then derive algebraic decompositions that reduce each stage to two linear subproblems with the same coefficient matrix and one scalar nonlinear equation. Finally, we establish unique solvability of the linear subproblems under suitable  assumptions and derive sufficient conditions for local existence and uniqueness of the scalar solution near unity.
	
	\subsection{Energy stability}
	We first establish unconditional stability of the proposed first- and second-order schemes with respect to the modified discrete energy
	\[
	\mathcal E^m:=W^m+\frac12\bigl((Q^m)^2-1\bigr),
	\qquad
	W^m:=\frac12(\varkappa^m,\varkappa^m)^h_{\Gamma^m}+\lambda|\Gamma^m|.
	\]
	\begin{theorem}[Energy stability of the first-order scheme]\label{thm:first_order_energy}
		Let
		$\bigl(V^{m+1},\varkappa^{m+1},\boldsymbol X^{m+1},\kappa^{m+1},\\Q^{m+1}\bigr)\in V^h\times V^h\times[V^h]^2\times V^h\times\mathbb R$
		be a solution of the first-order scheme \eqref{ful1}--\eqref{ful5}. Assume that $\lambda\geq0$ and that
		$\Gamma^m$ and $\Gamma^{m+1}$ are nondegenerate. Then
		\begin{equation}
			\frac{\mathcal E^{m+1}-\mathcal E^m}{\Delta t}\leq-(V^{m+1},V^{m+1})_{\Gamma^m}^h\leq0.
			\label{eq:first_order_energy_stability}
		\end{equation}
		Hence, $\mathcal E^{m+1}\leq\mathcal E^m$, $m=0,\ldots,M-1$, for any $\Delta t>0$.
	\end{theorem}

	\begin{proof}
		Testing \eqref{ful1}--\eqref{ful4} with	$\phi^h=V^{m+1}$, $\chi^h=\varkappa^{m+1}$,
		$\xi^h=-\lambda\kappa^{m+1}$, and	$\boldsymbol\eta^h	=\lambda(\boldsymbol X^{m+1}-\boldsymbol X^m)/\Delta t$,	respectively, and using \eqref{ful5}, we obtain
		\begin{align}\label{energy_eq}
			& \left({V}^{m+1},\,{V}^{m+1}\right)_{\Gamma^{m}}^{h} +\left( \frac{{\varkappa}^{m+1}-{\varkappa}^{m} }{\Delta t} ,\, {\varkappa}^{m+1} \right)_{\Gamma^{m}}^{h} + \lambda \left( \partial_s \boldsymbol{X}^{m+1},\; \partial_s \frac{(\boldsymbol{X}^{m+1}-\boldsymbol{X}^m)}{\Delta t}  \right)_{\Gamma^{m}}^{h}\nonumber\\
			&\qquad\qquad+\frac{{Q}^{m+1}-{Q}^m  }{\Delta t}Q^{m+1}+ \frac{1}{2}\left( \partial_s (\frac{\boldsymbol{X}^{m+1} - \boldsymbol{X}^{m}}{\Delta t}) \cdot \partial_s \boldsymbol{X}^{m+1} \varkappa^{m+1},\;  \varkappa^{m+1} \right)_{\Gamma^{m}}^{h}=0.
		\end{align}	
		Using $a(a-b)\geq\frac12(a^2-b^2)$ and the geometric estimates in \cite{garcke2025stable,li2025sinum}, we obtain
		\[
		\begin{aligned}
			&(\varkappa^{m+1}-\varkappa^m,
			\varkappa^{m+1})_{\Gamma^m}^h+\frac12\left(
			\left[
			\partial_s(\boldsymbol X^{m+1}-\boldsymbol X^m)
			\cdot\partial_s\boldsymbol X^{m+1}
			\right]\varkappa^{m+1},
			\varkappa^{m+1}
			\right)_{\Gamma^m}^h\\
			&\qquad\qquad\qquad\qquad\qquad\qquad\qquad\geq\frac12
			(\varkappa^{m+1},\varkappa^{m+1})_{\Gamma^{m+1}}^h
			-\frac12(\varkappa^m,\varkappa^m)_{\Gamma^m}^h,\\
			&\left(
			\partial_s\boldsymbol X^{m+1},
			\partial_s(\boldsymbol X^{m+1}-\boldsymbol X^m)
			\right)_{\Gamma^m}^h
			\geq|\Gamma^{m+1}|-|\Gamma^m|,\\
			&(Q^{m+1}-Q^m)Q^{m+1}
			\geq\frac12\bigl((Q^{m+1})^2-(Q^m)^2\bigr).
		\end{aligned}
		\]
		Since $\lambda\geq0$, combining these estimates with
		\eqref{energy_eq} yields
		\[
		\frac{W^{m+1}-W^m}{\Delta t}
		+\frac{(Q^{m+1})^2-(Q^m)^2}{2\Delta t}
		\leq-(V^{m+1},V^{m+1})_{\Gamma^m}^h
		\leq0.
		\]
		By the definition of $\mathcal E^m$, this proves
		\eqref{eq:first_order_energy_stability}.
	\end{proof}

	\medskip
	We next establish the energy stability of the second-order	 predictor--corrector scheme.
	\begin{theorem}[Energy stability of the second-order scheme]
		Let	$\bigl(V^{m+1},\varkappa^{m+1},\boldsymbol X^{m+1},\kappa^{m+1},\\ Q^{m+1}\bigr)\in V^h\times V^h\times[V^h]^2\times V^h\times\mathbb R$
		be a solution of the second-order scheme \eqref{ss1}--\eqref{sec6}. Assume that the discrete curves $\Gamma^m$, $\Gamma^{m+1}$, and $\widetilde\Gamma^{m+\frac12}$ are nondegenerate. Then
		\begin{align}
			\frac{\mathcal E^{m+1}-\mathcal E^m}{\Delta t}
			=
			-\left(V^{m+\frac12},V^{m+\frac12}\right)^h_{\widetilde\Gamma^{m+\frac12}}
			\leq 0,
		\end{align}
		Hence, $\mathcal E^{m+1}\leq\mathcal E^m$ for any $\Delta t>0$.
	\end{theorem}
	\begin{proof}
		By the definitions of $\varkappa_J^{m+\frac12}$ and $D_t^J\varkappa^{m+\frac{1}{2}}$,
		\begin{equation}
			\left(D_t^J\varkappa^{m+\frac{1}{2}},\varkappa_J^{m+\frac12}\right)^h_{\widetilde\Gamma^{m+\frac12}}
			=
			\frac{1}{2\Delta t}\left(J_+^{m+\frac12}(\varkappa^{m+1})^2-J_-^{m+\frac12}(\varkappa^m)^2,1\right)^h_{\widetilde\Gamma^{m+\frac12}}.
			\label{eq:proof-chain-kappa-1}
		\end{equation}
		By the definitions of $J_\pm^{m+\frac12}$, the mass-lumped
		inner products satisfy the exact transformation identities
		\[
		\left(J_+^{m+\frac12}u,v\right)^h_
		{\widetilde\Gamma^{m+\frac12}}
		=
		(u,v)^h_{\Gamma^{m+1}},
		\qquad
		\left(J_-^{m+\frac12}u,v\right)^h_
		{\widetilde\Gamma^{m+\frac12}}
		=
		(u,v)^h_{\Gamma^m}.
		\] Hence
		\begin{equation}
			\left(D_t^J\varkappa^{m+\frac12},\varkappa_J^{m+\frac12}\right)^h_{\widetilde\Gamma^{m+\frac12}}
			=
			\frac{1}{2\Delta t}\left[\left(\varkappa^{m+1},\varkappa^{m+1}\right)^h_{\Gamma^{m+1}}-\left(\varkappa^m,\varkappa^m\right)^h_{\Gamma^m}\right].
			\label{eq:J-curvature-chain}
		\end{equation}
		
		Next, by the definition
		$\boldsymbol\tau_J^{m+\frac12}=\frac{\boldsymbol X_\rho^{m+1}+\boldsymbol X_\rho^m}{|\boldsymbol X_\rho^{m+1}|+|\boldsymbol X_\rho^m|}$,
		we obtain the  identity
		\begin{equation}
			\boldsymbol\tau_J^{m+\frac12}\cdot\left(\boldsymbol X_\rho^{m+1}-\boldsymbol X_\rho^m\right)
			=
			\frac{|\boldsymbol X_\rho^{m+1}|^2-|\boldsymbol X_\rho^m|^2}{|\boldsymbol X_\rho^{m+1}|+|\boldsymbol X_\rho^m|}
			=
			|\boldsymbol X_\rho^{m+1}|-|\boldsymbol X_\rho^m|.
			\label{eq:DG-chain-proof}
		\end{equation}
		Therefore,
		\begin{equation}
			\int_{\mathbb I}\boldsymbol\tau_J^{m+\frac12}\cdot\partial_\rho\left(\frac{\boldsymbol X^{m+1}-\boldsymbol X^m}{\Delta t}\right)d\rho
			=
			\frac{1}{\Delta t}\int_{\mathbb I}\left(|\boldsymbol X_\rho^{m+1}|-|\boldsymbol X_\rho^m|\right)d\rho
			=
			\frac{|\Gamma^{m+1}|-|\Gamma^m|}{\Delta t}.
			\label{eq:length-chain-proof}
		\end{equation}
		
		Choosing $\phi^h=V^{m+\frac12}$, $\chi^h=\varkappa^{m+\frac12}$, $\xi^h=-\lambda\kappa^{m+\frac12}$, and $\boldsymbol\eta^h =\lambda(\boldsymbol X^{m+1}-\boldsymbol X^m)/\Delta t$ as test functions in \eqref{sec1}--\eqref{sec4},
		respectively, and summing the resulting identities, we obtain
		\begin{equation}
			\begin{aligned}
				&
				\left(V^{m+\frac12},V^{m+\frac12}\right)^h_{\widetilde\Gamma^{m+\frac12}}
				+\left(\frac{\varkappa^{m+1}-\varkappa^m}{\Delta t},\varkappa^{m+\frac12}\right)^h_{\widetilde\Gamma^{m+\frac12}}	+\lambda\left(\widetilde\partial_s\boldsymbol X^{m+\frac12},\widetilde\partial_s\frac{\boldsymbol X^{m+1}-\boldsymbol X^m}{\Delta t}\right)^h_{\widetilde\Gamma^{m+\frac12}}.\\
				&+\frac12 Q^{m+\frac12}\left(\widetilde\partial_s\widetilde{\boldsymbol U}^{m+\frac12}\cdot\widetilde\partial_s\widetilde{\boldsymbol X}^{m+\frac12}\,\widetilde\varkappa^{m+\frac12},\varkappa^{m+\frac12}\right)^h_{\widetilde\Gamma^{m+\frac12}}=0
			\end{aligned}
			\label{eq:pre-Q-balance}
		\end{equation}
		Using the auxiliary equation \eqref{sec6} together with	the discrete chain rules \eqref{eq:J-curvature-chain}	and \eqref{eq:length-chain-proof}, we obtain
		\begin{align}
			\left(V^{m+\frac12},V^{m+\frac12}\right)^h_{\widetilde\Gamma^{m+\frac12}}
			+\frac{Q^{m+1}-Q^m}{\Delta t}Q^{m+\frac12}+		\frac{1}{2\Delta t}\left[\left(\varkappa^{m+1},\varkappa^{m+1}\right)^h_{\Gamma^{m+1}}-\left(\varkappa^m,\varkappa^m\right)^h_{\Gamma^m}\right]\nonumber\\ +\lambda\frac{|\Gamma^{m+1}|-|\Gamma^m|}{\Delta t}=0.\label{m1}
		\end{align}
		Using the identity	$	\frac{Q^{m+1}-Q^m}{\Delta t}Q^{m+\frac12}	=\frac{(Q^{m+1})^2-(Q^m)^2}{2\Delta t}$, we rewrite \eqref{m1} as
		\begin{equation}
			\begin{aligned}
				&\frac{1}{\Delta t}\Bigl[\frac12(\varkappa^{m+1},\varkappa^{m+1})^h_{\Gamma^{m+1}}+\lambda|\Gamma^{m+1}|+\frac12(Q^{m+1})^2\\
				&\qquad-\frac12(\varkappa^m,\varkappa^m)^h_{\Gamma^m}-\lambda|\Gamma^m|-\frac12(Q^m)^2\Bigr]
				+\left(V^{m+\frac12},V^{m+\frac12}\right)^h_{\widetilde\Gamma^{m+\frac12}}=0.
			\end{aligned}
			\label{eq:second_energy_difference}
		\end{equation}
		By the definition of $\mathcal E^m$, this identity becomes
		\[
		\frac{\mathcal E^{m+1}-\mathcal E^m}{\Delta t}
		=-\left(V^{m+\frac12},V^{m+\frac12}\right)^h_
		{\widetilde\Gamma^{m+\frac12}}
		\leq0.
		\]
		Thus, $\mathcal E^{m+1}\leq\mathcal E^m$, completing the proof.
	\end{proof}
	
	\begin{rem}
		The ratios $J_\pm^{m+\frac12}$ are evaluated elementwise
		in the mass-lumped products. On each element $I_j$, we have
		$
		|\widetilde{\boldsymbol h}_j^{m+\frac12}|
		J_{+,j}^{m+\frac12}=|\boldsymbol h_j^{m+1}|,~
		|\widetilde{\boldsymbol h}_j^{m+\frac12}|
		J_{-,j}^{m+\frac12}=|\boldsymbol h_j^m|.
		$
		These identities exactly transform the edge-length weights
		on $\widetilde\Gamma^{m+\frac12}$ into those on
		$\Gamma^{m+1}$ and $\Gamma^m$, yielding the discrete
		curvature chain rule \eqref{eq:J-curvature-chain}.
	\end{rem}

	\subsection{Efficient implementation}
	This subsection presents efficient solution procedures for the proposed first- and second-order schemes. The main idea is to decompose the coupled discrete systems into linear subproblems and scalar nonlinear equations for the auxiliary variable, thereby avoiding the direct solution of high-dimensional nonlinear systems. We begin with the first-order scheme \eqref{ful1}--\eqref{ful5}.
	
	\subsubsection{Implementation of the  first-order scheme}	
	To decouple the global auxiliary variable $Q^{m+1}$ from the finite element unknowns, the solution is represented in the following form
	\begin{equation}\label{eq:first_order_decomposition}
		\begin{alignedat}{2}
			\boldsymbol X^{m+1}
			&=\boldsymbol X_1^{m+1}
			+Q^{m+1}\boldsymbol X_2^{m+1},
			\qquad&
			\varkappa^{m+1}
			&=\varkappa_1^{m+1}
			+Q^{m+1}\varkappa_2^{m+1},\\
			V^{m+1}
			&=V_1^{m+1}+Q^{m+1}V_2^{m+1},
			&
			\kappa^{m+1}
			&=\kappa_1^{m+1}+Q^{m+1}\kappa_2^{m+1}.
		\end{alignedat}
	\end{equation}
	Substituting \eqref{eq:first_order_decomposition} into \eqref{ful1}--\eqref{ful4} and collecting the terms independent of $Q^{m+1}$ and those dependent of $Q^{m+1}$ yields two linear subproblems.
	
	The $Q^{m+1}$-independent part yields the first linear subproblem: find
	$\bigl(V_1^{m+1},\varkappa_1^{m+1},\boldsymbol X_1^{m+1},\\ \kappa_1^{m+1}\bigr)\in V^h\times V^h\times[V^h]^2\times V^h$
	such that, for all
	$(\phi^h,\chi^h,\xi^h,\boldsymbol\eta^h)\in V^h\times V^h\times V^h\times[V^h]^2$,
	\begin{align}
		&(V_1^{m+1},\phi^h)^h_{\Gamma^m}=
		(\partial_s\varkappa_1^{m+1},\partial_s\phi^h)^h_{\Gamma^m}
		-\frac12((\varkappa^m)^2\varkappa_1^{m+1},\phi^h)^h_{\Gamma^m}
		+\lambda(\kappa_1^{m+1},\phi^h)^h_{\Gamma^m},	\label{ze1}	\\
		&	\left(\frac{\varkappa_1^{m+1}-\varkappa^m}{\Delta t},\chi^h\right)^h_{\Gamma^m}=	-(\partial_sV_1^{m+1},\partial_s\chi^h)^h_{\Gamma^m}	+\frac12((\varkappa^m)^2V_1^{m+1},\chi^h)^h_{\Gamma^m}	\nonumber\\
		&\quad\qquad\qquad\qquad\qquad\quad
		+\frac12\left(\boldsymbol\tau^m\cdot\frac{\boldsymbol X^m-\boldsymbol X^{m-1}}{\Delta t},\partial_s\varkappa_1^{m+1}\chi^h-\varkappa_1^{m+1}\partial_s\chi^h\right)^h_{\Gamma^m},
		\label{ze2}\\
		&	\left(\boldsymbol n^m\cdot\frac{\boldsymbol X_1^{m+1}-\boldsymbol X^m}{\Delta t},\xi^h\right)^h_{\Gamma^m}	=	(V_1^{m+1},\xi^h)^h_{\Gamma^m},
		\label{ze3}	\\
		&(\kappa_1^{m+1}\boldsymbol n^m,\boldsymbol\eta^h)^h_{\Gamma^m}
		=
		-(\partial_s\boldsymbol X_1^{m+1},\partial_s\boldsymbol\eta^h)^h_{\Gamma^m}.
		\label{ze4}
	\end{align}
	The coefficient of $Q^{m+1}$ yields the second linear subproblem: find
	$\bigl(V_2^{m+1},\varkappa_2^{m+1},\boldsymbol X_2^{m+1},\kappa_2^{m+1}\bigr)\in V^h\times V^h\times[V^h]^2\times V^h$
	such that, for all
	$(\phi^h,\chi^h,\xi^h,\boldsymbol\eta^h)\in V^h\times V^h\times V^h\times[V^h]^2$,
	\begin{align}
		&	(V_2^{m+1},\phi^h)^h_{\Gamma^m}	=
		(\partial_s\varkappa_2^{m+1},\partial_s\phi^h)^h_{\Gamma^m}
		-\frac12((\varkappa^m)^2\varkappa_2^{m+1},\phi^h)^h_{\Gamma^m}
		+\lambda(\kappa_2^{m+1},\phi^h)^h_{\Gamma^m},	\label{zze1}\\
		&	\left(\frac{\varkappa_2^{m+1}}{\Delta t},\chi^h\right)^h_{\Gamma^m}	=
		-(\partial_sV_2^{m+1},\partial_s\chi^h)^h_{\Gamma^m}
		+\frac12((\varkappa^m)^2V_2^{m+1},\chi^h)^h_{\Gamma^m}
		\nonumber\\
		&\quad\qquad\qquad\qquad\qquad
		+\frac12\left(\boldsymbol\tau^m\cdot\frac{\boldsymbol X^m-\boldsymbol X^{m-1}}{\Delta t},\partial_s\varkappa_2^{m+1}\chi^h-\varkappa_2^{m+1}\partial_s\chi^h\right)^h_{\Gamma^m}
		\nonumber\\
		&\quad\qquad\qquad\qquad\qquad
		-\frac12\left(\partial_s\frac{\boldsymbol X^m-\boldsymbol X^{m-1}}{\Delta t}\cdot\partial_s\boldsymbol X^m\,\varkappa^m,\chi^h\right)^h_{\Gamma^m},
		\label{zze2}\\
		&\left(\boldsymbol n^m\cdot\frac{\boldsymbol X_2^{m+1}}{\Delta t},\xi^h\right)^h_{\Gamma^m}	=	(V_2^{m+1},\xi^h)^h_{\Gamma^m},	\label{zze3}
		\\
		&	(\kappa_2^{m+1}\boldsymbol n^m,\boldsymbol\eta^h)^h_{\Gamma^m}=
		-(\partial_s\boldsymbol X_2^{m+1},\partial_s\boldsymbol\eta^h)^h_{\Gamma^m}.
		\label{zze4}
	\end{align}
	
	After solving the two linear subproblems, we substitute
	\eqref{eq:first_order_decomposition} into \eqref{ful5}
	to obtain a scalar nonlinear equation for $Q^{m+1}$:
	\begin{equation}\label{eq:first_order_Q}
		F_m(Q^{m+1})=0.
	\end{equation}
	Once $Q^{m+1}$ is determined, the remaining unknowns are
	recovered from \eqref{eq:first_order_decomposition}.
	
	The computational procedure for the first-order scheme is summarized in Algorithm~\ref{alg:first_order}.
	\begin{algorithm}[htbp]
		\caption{First-order energy-stable algorithm}
		\label{alg:first_order}
		\begin{algorithmic}[1]
			\Statex \textbf{Input:} Initial polygonal curve $\Gamma^0$ represented by $\boldsymbol X^0\in[V^h]^2$, time step $\Delta t$, final time $T$.
			\Statex \textbf{Output:} Discrete solution $(\boldsymbol X^M,\varkappa^M,V^M,\kappa^M,Q^M)$, where $M=T/\Delta t$.
			\State \textbf{Initialization:} Compute the initial curvature $\kappa^0$ on $\Gamma^0$ from $\kappa\boldsymbol n=\partial_{ss}\boldsymbol X$ using the least-squares method~\cite{BARRETT2007441,JIANG2024113220}, and initialize $\varkappa^0$ consistently with $\kappa^0$. Set $Q^0=1$.
			\State  Compute the initial normal velocity $V^0$ from
			$(V^0,\phi^h)^h_{\Gamma^0}=(\partial_s\varkappa^0,\partial_s\phi^h)^h_{\Gamma^0}-\tfrac12((\varkappa^0)^3,\phi^h)^h_{\Gamma^0}+\lambda(\kappa^0,\phi^h)^h_{\Gamma^0}$,\quad $\forall\,\phi^h\in V^h$.
			\State Compute the initial normal $\boldsymbol n^0$, define the
			initial vector velocity with zero tangential component by
			$\boldsymbol U^0:=V^0\boldsymbol n^0$, and set
			$\boldsymbol X^{-1}:=\boldsymbol X^0-\Delta t\,\boldsymbol U^0$ so that
			$\frac{\boldsymbol X^0-\boldsymbol X^{-1}}{\Delta t}=\boldsymbol U^0$.
			\For{$m=0,\ldots,M-1$}
			\State \textbf{Linear solves:} Solve \eqref{ze1}--\eqref{ze4} for $(\boldsymbol X_1^{m+1},\varkappa_1^{m+1},V_1^{m+1},\kappa_1^{m+1})$, and solve \eqref{zze1}--\eqref{zze4} for $(\boldsymbol X_2^{m+1},\varkappa_2^{m+1},V_2^{m+1},\kappa_2^{m+1})$.
			\State \textbf{Scalar nonlinear solve:} Substitute \eqref{eq:first_order_decomposition} into \eqref{ful5} and solve $F_m(Q^{m+1})=0$.
			\State \textbf{Update:} Recover $(\boldsymbol X^{m+1},\varkappa^{m+1},V^{m+1},\kappa^{m+1})$ from \eqref{eq:first_order_decomposition}.
			\State Update the discrete curve $\Gamma^{m+1}$ and the associated geometric quantities.
			\EndFor
		\end{algorithmic}
	\end{algorithm}

	\subsubsection{Implementation of the  second-order scheme}
	The predictor uses the first-order solution procedure with	time increment $\Delta t/2$.	For the corrector, the  decomposition	\eqref{eq:first_order_decomposition} yields two linear	subproblems with the same coefficient matrix.	Substituting their solutions into \eqref{sec6} gives	a scalar nonlinear equation for $Q^{m+1}$.
	Once $Q^{m+1}$ is determined, the remaining unknowns  are recovered from the decomposition.

	The main second-order time-stepping scheme is applied for $m\geq2$, since the explicit midpoint velocity
	$\widetilde{\boldsymbol U}^{m+\frac12}=\frac{2\boldsymbol X^m-3\boldsymbol X^{m-1}+\boldsymbol X^{m-2}}{\Delta t}$ requires the solutions at the preceding time levels. Given the prescribed initial solution at $t_0=0$, we compute the starting values at $t_1$ and $t_2$ using an energy-dissipative scheme, such as the first-order scheme proposed above or the energy-stable nonlinear scheme of Garcke, N{\"u}rnberg, and Zhao~\cite{garcke2025stable}. The start-up procedure is required to provide $O(\Delta t^2)$ temporal errors. When the first-order scheme is used, the  variables $V$ and $\kappa$ are reconstructed consistently on the current curves to obtain compatible starting values without altering the modified energy. 	Let $\boldsymbol Y^m:=(\boldsymbol X^m,\varkappa^m,V^m,\kappa^m)$. With the second-order start-up values available, the main
	predictor--corrector scheme is applied for $m\geq2$. The complete
	procedure is summarized in Algorithm~\ref{alg:pc_cn}.

	\begin{algorithm}[htbp]
		\caption{Second-order energy-stable algorithm}
		\label{alg:pc_cn}
		\begin{algorithmic}[1]
			\For{$m=2,\ldots,M-1$}
			\State \textbf{Predictor:} Compute $\bigl(\widetilde{\boldsymbol X}^{m+\frac12},\widetilde\varkappa^{m+\frac12},\widetilde V^{m+\frac12},\widetilde\kappa^{m+\frac12},\widetilde Q^{m+\frac12}\bigr)$ using the first-order scheme with the half time step $\Delta t/2$.
			\State Construct the predicted curve $\widetilde\Gamma^{m+\frac12}$ and  $\widetilde{\boldsymbol\tau}^{m+\frac12}$, $\widetilde{\boldsymbol n}^{m+\frac12}$, and $\widetilde\partial_s$.
			\State Compute the extrapolated midpoint velocity $\widetilde{\boldsymbol U}^{m+\frac12}=\frac{2\boldsymbol X^m-3\boldsymbol X^{m-1}+\boldsymbol X^{m-2}}{\Delta t}$.
			\State \textbf{Corrector:} Solve the two linear subproblems for $\boldsymbol Y_1^{m+1}$ and $\boldsymbol Y_2^{m+1}$.
			\State Substitute $\boldsymbol Y^{m+1}=\boldsymbol Y_1^{m+1}+Q^{m+1}\boldsymbol Y_2^{m+1}$ into the auxiliary equation \eqref{sec6} and solve the resulting scalar nonlinear equation $F_m(Q^{m+1})=0$.
			\State Recover $\bigl(\boldsymbol X^{m+1},\varkappa^{m+1},V^{m+1},\kappa^{m+1}\bigr)$ from the  decomposition.
			\State Update the discrete curve $\Gamma^{m+1}$ and the associated geometric quantities.
			\EndFor
		\end{algorithmic}
	\end{algorithm}

	\subsection{Existence and uniqueness}
	The solvability of the proposed first- and second-order schemes
	is established in two steps. First, we prove the unique solvability
	of the linear subproblems. Second, we establish the local existence
	and uniqueness of solutions to the  nonlinear algebraic
	equations for the auxiliary variable $Q$. Let $\boldsymbol\omega^m\in[V^h]^2$ denote the vertex normal defined by \cite{garcke2025stable}.
	
	\[
	(\boldsymbol\omega^m,\boldsymbol\eta^h)^h_{\Gamma^m}
	=(\boldsymbol n^m,\boldsymbol\eta^h)^h_{\Gamma^m},
	\qquad
	\forall\,\boldsymbol\eta^h\in[V^h]^2.
	\]
	
	\subsubsection{Unique solvability of the linear subproblems}
	We first establish the unique solvability of the two linear subproblems.
	\begin{theorem}Assume that $\text{(A1)}~
		\lambda>0~\text{or}~
		\boldsymbol\omega^m(\rho_j)\neq\boldsymbol0, \text{for all }j=1,\ldots,N$, \text{(A2)}~dim $\operatorname{span}$$\{\boldsymbol\omega^m(\rho_j)\}_{j=1}^{N}=2.$
		Then the linear subproblems \eqref{ze1}--\eqref{ze4}
		and \eqref{zze1}--\eqref{zze4} for the first-order scheme admit unique solutions $	(V_\alpha^{m+1}, \varkappa_\alpha^{m+1},
		\boldsymbol{X}_\alpha^{m+1}, \kappa_\alpha^{m+1})
		\in V^h \times V^h \times [V^h]^2 \times V^h,~ \alpha=1,2,$	respectively.
	\end{theorem}
	
	\begin{proof}
		The two subproblems have the same square coefficient matrix. It therefore suffices to show that their common homogeneous problem has only the zero solution.
		
		Let $(V,\varkappa,\boldsymbol X,\kappa)$ solve the homogeneous problem. Testing its four equations with $\phi^h=V$, $\chi^h=\varkappa$, $\xi^h=-\lambda\kappa$, and $\boldsymbol\eta^h=\lambda\boldsymbol X/\Delta t$, respectively, then summing and multiplying by $\Delta t$, gives
		\begin{equation}\label{nn1_1}
			(\varkappa,\varkappa)^h_{\Gamma^m}+\Delta t(V,V)^h_{\Gamma^m}+\lambda(\partial_s\boldsymbol X,\partial_s\boldsymbol X)^h_{\Gamma^m}=0.
		\end{equation}
		Hence $V=\varkappa=0$.
		
		Testing the remaining two equations with $\xi^h=\kappa$ and $\boldsymbol\eta^h=\boldsymbol X$ yields
		\[
		(\partial_s\boldsymbol X,\partial_s\boldsymbol X)^h_{\Gamma^m}=-(\kappa\boldsymbol n^m,\boldsymbol X)^h_{\Gamma^m}=-\Delta t(\kappa,V)^h_{\Gamma^m}=0.
		\]
		Thus $\boldsymbol X=\boldsymbol X_c$ is constant. The normal-velocity equation and the definition of $\boldsymbol\omega^m$ imply
		$\boldsymbol X_c\cdot\boldsymbol\omega^m(\rho_j)=0$ for $j=1,\ldots,N$. The spanning assumption therefore gives $\boldsymbol X_c=\boldsymbol0$. 	If $\lambda>0$, the first homogeneous equation with $\phi^h=\kappa$ gives $\kappa=0$. If $\lambda=0$, the geometric equation implies
		$\kappa(\rho_j)\boldsymbol\omega^m(\rho_j)=\boldsymbol0$ for $j=1,\ldots,N$. The assumed nonvanishing of the vertex normals again gives $\kappa=0$. Thus the homogeneous kernel is trivial, which proves existence and uniqueness.
	\end{proof}
	
	\medskip
	We next establish the unique solvability of the linear subproblems  in the second-order scheme.
	
	\begin{theorem}
		Assume that $\lambda\geq0$ and that $\Gamma^m$ and
		$\widetilde\Gamma^{m+\frac12}$ are nondegenerate.
		Suppose that assumptions~(A1)--(A2) hold for the corresponding
		mass-lumped vertex normals on both curves.
		Then the linear subproblems associated with the predictor
		equations \eqref{ss1}--\eqref{ss4} and the corrector equations
		\eqref{sec1}--\eqref{sec4} admit unique solutions.
	\end{theorem}
	
	\begin{proof}
		The predictor linear subproblems are uniquely solvable
		by the first-order result.
		For the corrector, let
		$(V,\varkappa,\boldsymbol X,\kappa)$ solve the common
		homogeneous system associated with
		\eqref{sec1}--\eqref{sec4}.
		Testing its four equations with
		$\phi^h=V$, $\chi^h=\varkappa$,
		$\xi^h=-\lambda\kappa$, and
		$\boldsymbol\eta^h=2\lambda\boldsymbol X/\Delta t$,
		respectively, then summing and multiplying by $\Delta t$,
		gives
		\begin{equation}\label{nn1_2}
			\begin{aligned}
				0={}&(\varkappa,\varkappa)^h_
				{\widetilde\Gamma^{m+\frac12}}
				+\frac{\Delta t}{2}(V,V)^h_
				{\widetilde\Gamma^{m+\frac12}}+\lambda
				(\widetilde\partial_s\boldsymbol X,
				\widetilde\partial_s\boldsymbol X)^h_
				{\widetilde\Gamma^{m+\frac12}}.
			\end{aligned}
		\end{equation}
		Hence $V=\varkappa=0$.
		Testing the normal and geometric equations with
		$\xi^h=\kappa$ and $\boldsymbol\eta^h=\boldsymbol X$,
		respectively, gives $\widetilde\partial_s\boldsymbol X=0$.
		Thus $\boldsymbol X$ is constant, and assumption~(A2)
		on the predicted curve implies $\boldsymbol X=\boldsymbol0$.	If $\lambda>0$, the first homogeneous equation with
		$\phi^h=\kappa$ gives $\kappa=0$.
		If $\lambda=0$, the geometric equation and assumption~(A1)
		on the predicted curve give the same conclusion.
		Since the common coefficient matrix is square and has
		a trivial kernel, the corrector linear subproblems
		are uniquely solvable.
	\end{proof}
	
	\subsubsection{Local existence and uniqueness of the auxiliary variable}
	We establish sufficient conditions for the existence of
	a unique scalar solution near unity. We first prove a general
	contraction criterion and then apply it to the auxiliary
	equations of the first- and second-order schemes.
	
	\begin{lemma}[A local solvability criterion]\label{unique1}
		Let $F_m\in C^2(U_m)$, where $U_m$ is an open interval containing $1$, and define $\delta_m:=|F_m(1)|$, $\sigma_m:=|F_m'(1)|>0$. For $\gamma_m\in(0,1)$, set
		$\rho_m:=\delta_m/[\sigma_m(1-\gamma_m)]$, $I_m:=[1-\rho_m,1+\rho_m]$.
		Assume that $I_m\subset U_m$, and let $L_m:=\sup_{q\in I_m}|F_m''(q)|$. If
		\begin{equation}\label{eq:Q_general_condition}
			\frac{L_m\rho_m}{\sigma_m}\leq\gamma_m,
		\end{equation}
		then $F_m(q)=0$ admits a unique solution $Q^{m+1}\in I_m$, with $|Q^{m+1}-1|\leq\rho_m$.
	\end{lemma}
	
	\begin{proof}
		Define the simplified Newton map $\Phi_m(q):=q-F_m(q)/F_m'(1)$. For $q\in I_m$, the mean value theorem gives
		$|\Phi_m'(q)|=|F_m'(1)-F_m'(q)|/\sigma_m\leq L_m\rho_m/\sigma_m\leq\gamma_m<1$.
		Moreover,
		\[
		|\Phi_m(q)-1|\leq|\Phi_m(q)-\Phi_m(1)|+|\Phi_m(1)-1|\leq\gamma_m\rho_m+\delta_m/\sigma_m=\rho_m.
		\]
		Thus $\Phi_m$ maps $I_m$ into itself and is a contraction. The Banach fixed-point theorem gives a unique fixed point $Q^{m+1}\in I_m$, which is the unique zero of $F_m$ in this interval.
	\end{proof}
	
	\medskip
	\paragraph{\textbf{First-order scheme}}
	Substituting the decomposition
	\eqref{eq:first_order_decomposition} into \eqref{ful5}
	yields $F_m(Q^{m+1})=0$, where
	\begin{equation}\label{eq:FmQ}
		\begin{aligned}
			F_m(q):=\frac{q(q-Q^m)}{\Delta t}-\frac q2\left(
			(\partial_s\boldsymbol U^m\cdot
			\partial_s\boldsymbol X^m)\varkappa^m,
			\varkappa(q)
			\right)^h_{\Gamma^m}+\frac12\left(
			(\partial_s\boldsymbol W(q)\cdot
			\partial_s\boldsymbol X(q))\varkappa(q),
			\varkappa(q)
			\right)^h_{\Gamma^m},
		\end{aligned}
	\end{equation}
	with $\boldsymbol X(q):=\boldsymbol X_1^{m+1}+q\boldsymbol X_2^{m+1}, \varkappa(q):=\varkappa_1^{m+1}+q\varkappa_2^{m+1},
	\boldsymbol U^m
	:=\frac{\boldsymbol X^m-\boldsymbol X^{m-1}}{\Delta t},
	\boldsymbol W(q)
	:=\frac{\boldsymbol X(q)-\boldsymbol X^m}{\Delta t}.$
	The coefficients of $F_m$ are given in Appendix~\ref{app:quartic_coefficients}.
	For brevity, we define $\|u\|_{\Gamma^m,h}^2:=(u,u)^h_{\Gamma^m}$.
	
	\medskip
	\begin{theorem}[Local solvability of the first-order scheme]
		\label{thm:Q_first_order}
		Let $F_m(q)=\sum_{k=0}^{4}a_{k,m}q^k$, and define
		$\delta_m=|F_m(1)|$, $\sigma_m=|F_m'(1)|>0$, $\rho_m=2\delta_m/\sigma_m$.
		If
		\begin{equation}
			\frac{4\delta_m|F_m''(1)|}{\sigma_m^2}
			+\frac{8\delta_m^2|F_m'''(1)|}{\sigma_m^3}
			+\frac{192a_{4,m}\delta_m^3}{\sigma_m^4}\le1,
			\label{eq:Q_solvability_condition}
		\end{equation}
		then $F_m(Q^{m+1})=0$ admits a unique solution in $I_m=[1-\rho_m,1+\rho_m]$, satisfying $|Q^{m+1}-1|\le\rho_m$.
		
		Furthermore, let $\lambda>0$ and consider a fixed nondegenerate spatial mesh. Assume that
		\begin{equation}
			\|\varkappa^m\|_{\Gamma^m,h}+\|\partial_s\boldsymbol X^m\|_{\Gamma^m,h}+\|H_m\|_{\Gamma^m,h}\le C,
			\label{eq:Q_known_data_bound}
		\end{equation}
		where
		$H_m=\partial_s[(\boldsymbol X^m-\boldsymbol X^{m-1})/\Delta t]\cdot\partial_s\boldsymbol X^m\,\varkappa^m$,
		and $C$ is independent of $\Delta t$.
		
		If $|2-Q^m|\ge c_0>0$ and $\delta_m=o(\Delta t^{-1})$, then condition \eqref{eq:Q_solvability_condition} holds for sufficiently small $\Delta t$, and $|Q^{m+1}-1|=o(1)$. In particular, if $\delta_m=O(1)$, then $|Q^{m+1}-1|=O(\Delta t)$.
	\end{theorem}
	
	\begin{proof}
		Since $\boldsymbol X(q)$ and $\varkappa(q)$ depend affinely on $q$, and the discrete inner products are evaluated on the known curve $\Gamma^m$, $F_m$ is a polynomial of degree at most four. In particular,
		\[
		a_{4,m}=\frac{1}{2\Delta t}\left(|\partial_s\boldsymbol X_2^{m+1}|^2(\varkappa_2^{m+1})^2,1\right)^h_{\Gamma^m}\ge0.
		\]
		
		The exact expansion
		$F_m''(q)=F_m''(1)+F_m'''(1)(q-1)+12a_{4,m}(q-1)^2$
		implies
		\begin{equation}
			L_m\le|F_m''(1)|+|F_m'''(1)|\rho_m+12a_{4,m}\rho_m^2.
			\label{eq:Q_L_bound}
		\end{equation}
		Choosing $\gamma_m=\frac12$ in Lemma~\ref{unique1}, we obtain $\rho_m=2\delta_m/\sigma_m$. Hence
		\begin{equation}
			\frac{L_m\rho_m}{\sigma_m}
			\le\frac{2\delta_m|F_m''(1)|}{\sigma_m^2}
			+\frac{4\delta_m^2|F_m'''(1)|}{\sigma_m^3}
			+\frac{96a_{4,m}\delta_m^3}{\sigma_m^4}.
			\label{eq:Q_L_explicit}
		\end{equation}
		Thus, condition \eqref{eq:Q_solvability_condition} ensures that $L_m\rho_m/\sigma_m\le\frac12$. The first assertion follows from Lemma~\ref{unique1}.
		
		We next derive estimates for the polynomial coefficients using the two linear subproblems. Choosing
		$\phi^h=V_1^{m+1}$, $\chi^h=\varkappa_1^{m+1}$, $\xi^h=-\lambda\kappa_1^{m+1}$, and
		$\boldsymbol\eta^h=\lambda(\boldsymbol X_1^{m+1}-\boldsymbol X^m)/\Delta t$
		in \eqref{ze1}--\eqref{ze4}, respectively, and summing the resulting equations yields
		\begin{equation}
			\begin{aligned}
				&\|\varkappa_1^{m+1}\|_{\Gamma^m,h}^2
				+\lambda\|\partial_s\boldsymbol X_1^{m+1}\|_{\Gamma^m,h}^2
				+\|\varkappa_1^{m+1}-\varkappa^m\|_{\Gamma^m,h}^2\\
				&\quad+\lambda\|\partial_s(\boldsymbol X_1^{m+1}-\boldsymbol X^m)\|_{\Gamma^m,h}^2
				+2\Delta t\|V_1^{m+1}\|_{\Gamma^m,h}^2
				=\|\varkappa^m\|_{\Gamma^m,h}^2+\lambda\|\partial_s\boldsymbol X^m\|_{\Gamma^m,h}^2.
			\end{aligned}
			\label{eq:Q_linear_estimate_1}
		\end{equation}
		
		Similarly, testing \eqref{zze1}--\eqref{zze4} with $\phi^h=V_2^{m+1}$, $\chi^h=\varkappa_2^{m+1}$, $\xi^h=-\lambda\kappa_2^{m+1}$, and $\boldsymbol\eta^h=\lambda\boldsymbol X_2^{m+1}/\Delta t$, and applying Young's inequality gives
		\begin{equation}
			\frac12\|\varkappa_2^{m+1}\|_{\Gamma^m,h}^2
			+\Delta t\|V_2^{m+1}\|_{\Gamma^m,h}^2
			+\lambda\|\partial_s\boldsymbol X_2^{m+1}\|_{\Gamma^m,h}^2
			\le\frac{\Delta t^2}{8}\|H_m\|_{\Gamma^m,h}^2.
			\label{eq:Q_linear_estimate_2}
		\end{equation}
		
		For fixed $h$ and $\lambda>0$, inverse estimates and finite-dimensional norm equivalence imply
		\[
		\begin{aligned}
			&\|\varkappa_1^{m+1}\|_{L^\infty}+\|\partial_s\boldsymbol X_1^{m+1}\|_{L^\infty}
			+\|\partial_s(\boldsymbol X_1^{m+1}-\boldsymbol X^m)\|_{L^\infty}\le C,\\
			&\|\varkappa_2^{m+1}\|_{L^\infty}+\|\partial_s\boldsymbol X_2^{m+1}\|_{L^\infty}\le C\Delta t.
		\end{aligned}
		\]
		
		Substituting these estimates into the explicit expressions for $a_{k,m}$ yields
		\begin{equation}
			\begin{aligned}
				a_{0,m}&=O(\Delta t^{-1}),&
				a_{1,m}&=-\frac{Q^m}{\Delta t}+O(1),&
				a_{2,m}&=\frac1{\Delta t}+O(\Delta t),&
				a_{3,m}&=O(\Delta t^2),&
				a_{4,m}&=O(\Delta t^3).
			\end{aligned}
			\label{eq:Q_coefficient_estimates}
		\end{equation}
		Consequently,
		\begin{equation}
			F_m'(1)=\frac{2-Q^m}{\Delta t}+O(1),\quad
			F_m''(1)=\frac2{\Delta t}+O(\Delta t),\quad
			F_m'''(1)=O(\Delta t^2).
			\label{eq:Q_derivative_estimates}
		\end{equation}
		
		The bound $|2-Q^m|\geq c_0$ implies
		$\sigma_m\geq c_0/(2\Delta t)$ for sufficiently small $\Delta t$. Combining this estimate with \eqref{eq:Q_derivative_estimates}, we obtain
		\begin{equation}
			\frac{4\delta_m|F_m''(1)|}{\sigma_m^2}
			+\frac{8\delta_m^2|F_m'''(1)|}{\sigma_m^3}
			+\frac{192a_{4,m}\delta_m^3}{\sigma_m^4}
			\le C(\delta_m\Delta t+\delta_m^2\Delta t^5+\delta_m^3\Delta t^7).
			\label{eq:Q_condition_asymptotic}
		\end{equation}
		
		If $\delta_m=o(\Delta t^{-1})$, the right-hand side tends to zero as $\Delta t\to0$. Hence condition \eqref{eq:Q_solvability_condition} holds for sufficiently small $\Delta t$. By Lemma~\ref{unique1},
		$|Q^{m+1}-1|\le2\delta_m/\sigma_m\le C\delta_m\Delta t=o(1)$.
		If $\delta_m=O(1)$, this estimate improves to $|Q^{m+1}-1|=O(\Delta t)$.
	\end{proof}
	
	\begin{rem}
		The coefficients of $F_m$ are explicitly determined by the uniquely solvable linear subproblems. Hence, condition \eqref{eq:Q_solvability_condition} can be verified at each time step.	The energy estimates alone do not guarantee $\delta_m=o(\Delta t^{-1})$. Indeed,
		$F_m(1)=a_{0,m}+(1-Q^m)/\Delta t+O(1)$, so that
		\[
		\begin{aligned}
			\delta_m\Delta t
			&\le\frac12\left|\left(\partial_s(\boldsymbol X_1^{m+1}-\boldsymbol X^m)\cdot\partial_s\boldsymbol X_1^{m+1}(\varkappa_1^{m+1})^2,1\right)^h_{\Gamma^m}\right|+|Q^m-1|+C\Delta t.
		\end{aligned}
		\]
		Thus, an additional consistency estimate is needed to guarantee the residual condition. For example, if
		$\|\partial_s(\boldsymbol X_1^{m+1}-\boldsymbol X^m)\|_{L^\infty}\le C\Delta t$
		and $|Q^m-1|\le C\Delta t$, then $\delta_m=O(1)$ and $|Q^{m+1}-1|=O(\Delta t)$.
	\end{rem}
	
	\medskip
	\paragraph{\textbf{Second-order scheme}}We now establish local solvability of the global  auxiliary equation in the second-order corrector. Substituting the decomposition of the corrector unknowns into \eqref{sec6}, the auxiliary equation can be expressed as
	\begin{equation}
		F_m^{(2)}(q)=a_{0,m}^{(2)}+a_{1,m}^{(2)}q+a_{2,m}^{(2)}q^2+\mathcal R_m(q)=0.
		\label{eq:second_F_decomposition}
	\end{equation}
	The two linear subproblems determine the coefficients $a_{k,m}^{(2)}$,
	$k=0,1,2$. The remainder $\mathcal R_m(q)$ contains the nonlinear terms
	arising from the edge lengths of $\Gamma[\boldsymbol X(q)]$.
	Explicit expressions for the coefficients and the remainder are given
	in Appendix~\ref{app:second_auxiliary}.
	
	\begin{theorem}[Local solvability of the second-order scheme]
		\label{thm:Q_second_order}
		Let $F_m^{(2)}$ be defined by \eqref{eq:second_F_decomposition}, and set
		$\delta_m=|F_m^{(2)}(1)|$, $\sigma_m=|(F_m^{(2)})'(1)|>0$, $\rho_m=\frac{2\delta_m}{\sigma_m}$.
		Assume that all edges of the discrete curve $\Gamma[\boldsymbol X(q)]$ remain nondegenerate on an open neighborhood of $I_m=[1-\rho_m,1+\rho_m]$. Define $L_m=\sup_{q\in I_m}|(F_m^{(2)})''(q)|$. If
		\begin{equation}
			\frac{2\delta_m L_m}{\sigma_m^2}\le\frac12,
			\label{eq:Q_second_condition}
		\end{equation}
		then $F_m^{(2)}(Q^{m+1})=0$ has a unique solution in $I_m$, and $|Q^{m+1}-1| \leq \frac{2\delta_m}{\sigma_m}.$
	\end{theorem}
	
	\begin{proof}
		By the  decomposition and the edge nondegeneracy assumption, the geometric terms in $F_m^{(2)}$ are twice continuously differentiable on $I_m$. Hence $F_m^{(2)}\in C^2(I_m)$. 
		By \eqref{eq:Q_second_condition}, $	\frac{L_m\rho_m}{\sigma_m}
		= \frac{2\delta_m L_m}{\sigma_m^2}
		\leq \frac12.$ The result follows from Lemma~\ref{unique1} with $\gamma_m=\frac12$.
	\end{proof}
	
	\begin{rem}
		The first-order auxiliary equation is polynomial of degree at most four,
		whereas the second-order equation is generally nonpolynomial because
		the edge lengths depend on $q$. The formulas in
		Appendix~\ref{app:second_auxiliary} can be used to bound $L_m$ and
		check the local solvability condition~\eqref{eq:Q_second_condition}.
	\end{rem}

\section{Numerical results}\label{part5}
In this section, we present numerical experiments to investigate
the accuracy, energy dissipation, mesh quality and computational
efficiency of the proposed first- and second-order PFEM schemes.
The scalar nonlinear equations for the auxiliary variable
$Q^{m+1}$ are solved by the secant method.
Let $Q_k$ denote the $k$th iterate at a given time step.
The iteration is terminated when $|Q_{k+1}-Q_k|<10^{-10}
~\text{or}~
|F_m(Q_{k+1})|<10^{-10},$
with a maximum of 100 iterations per scalar solve.

\subsection{Accuracy test}
\paragraph{\textbf{Example 1: Convergence of the first- and second-order schemes}}
We first consider an expanding circle with $\lambda=0$. The exact solution is
\begin{equation}
	\begin{aligned}
		\boldsymbol X(\rho,t)=R(t)
		\begin{pmatrix}\cos g(\rho)\\ \sin g(\rho)\end{pmatrix},~
		R(t)=(1+2t)^{1/4},~\varkappa(t)=\kappa(t)=R(t)^{-1},~ V(t)=-\frac12R(t)^{-3}. \nonumber
	\end{aligned}
	\label{eq:circle_exact_numerics}
\end{equation}
We use $g(\rho)=2\pi\rho$ and
$g(\rho)=2\pi\rho+0.1\sin(2\pi\rho)$ to generate uniform and nonuniform
initial nodal distributions, respectively. We define the errors
\begin{align}
	e_{\boldsymbol X}&=\max_{0\le m\le M}\max_{1\le j\le N}
	\min_{\rho\in\mathbb I}|\boldsymbol X^m(\rho_j)-\boldsymbol X(\rho,t_m)|,\quad e_Q=\max_{0\le m\le M}|Q^m-1|,
	\label{eq:position_error}\nonumber\\
	e_{\varkappa}&=\max_{0\le m\le M}\max_{1\le j\le N}
	|\varkappa^m(\rho_j)-\varkappa(t_m)|,\quad
	e_V=\max_{0\le m\le M}\max_{1\le j\le N}|V^m(\rho_j)-V(t_m)|.\nonumber
\end{align}
The geometric curvature error $e_\kappa$ is defined analogously.

We further examine the convergence behavior for an initial ellipse defined by $\frac{x^2}{2}+y^2=1$.
We employ the manifold distance between two closed curves
$\Gamma_1$ and $\Gamma_2$, defined by
\begin{equation}
	M(\Gamma_1,\Gamma_2)
	:=
	2|\Omega_1\cup\Omega_2|
	-|\Omega_1|-|\Omega_2|,
	\label{eq:manifold_distance}
\end{equation}
where $|\Omega|$ denotes the area of  $\Omega$. We define the manifold Cauchy error by
\begin{equation}
	e_{h,\Delta t}^{\mathrm C}
	:=M\bigl(\Gamma_{h,\Delta t}(T),
	\Gamma_{h/2,\Delta t/r}(T)\bigr),
\end{equation}
where $\Gamma_{h,\Delta t}(T)$ is the numerical curve at time $T$,
and $r=4$ and $r=2$ for the first- and second-order schemes,
respectively.

\begin{figure}[htbp]
	\centering
	\includegraphics[width=0.32\linewidth]{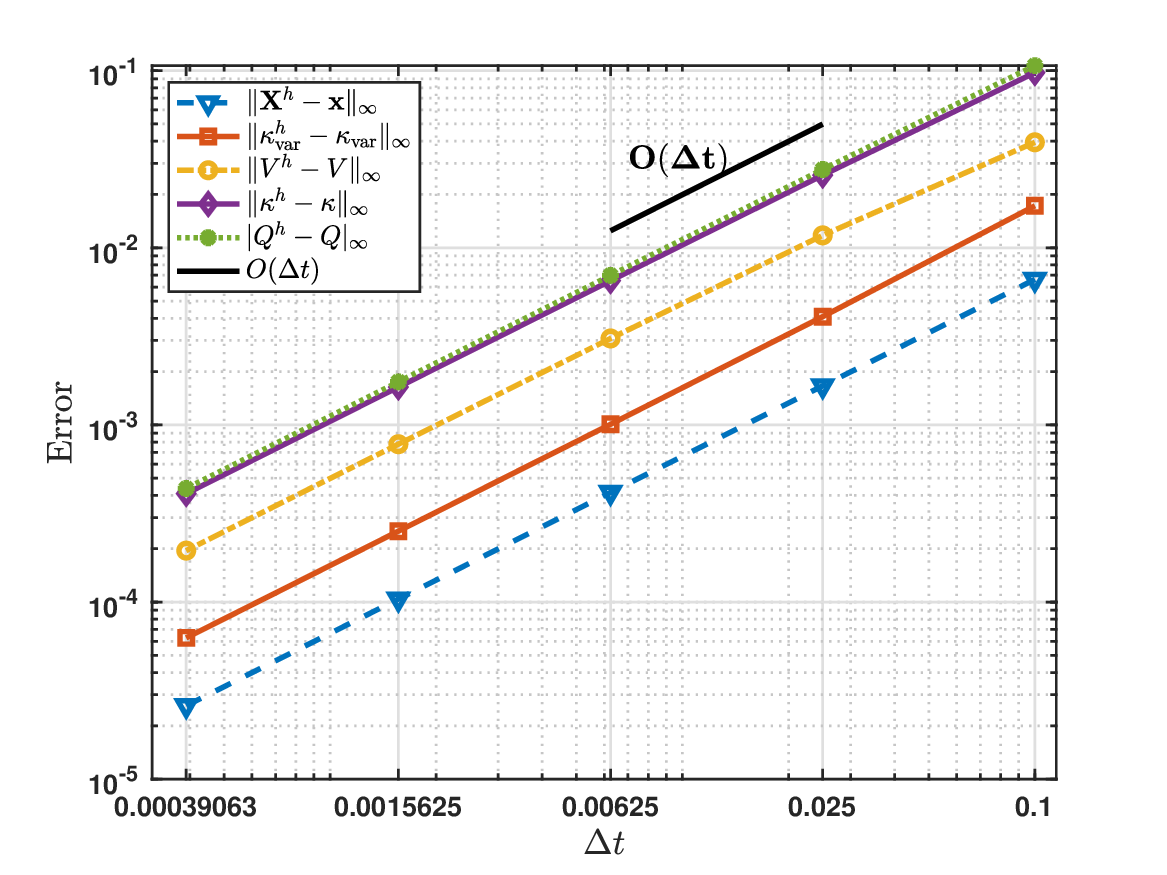}
	\includegraphics[width=0.32\linewidth]{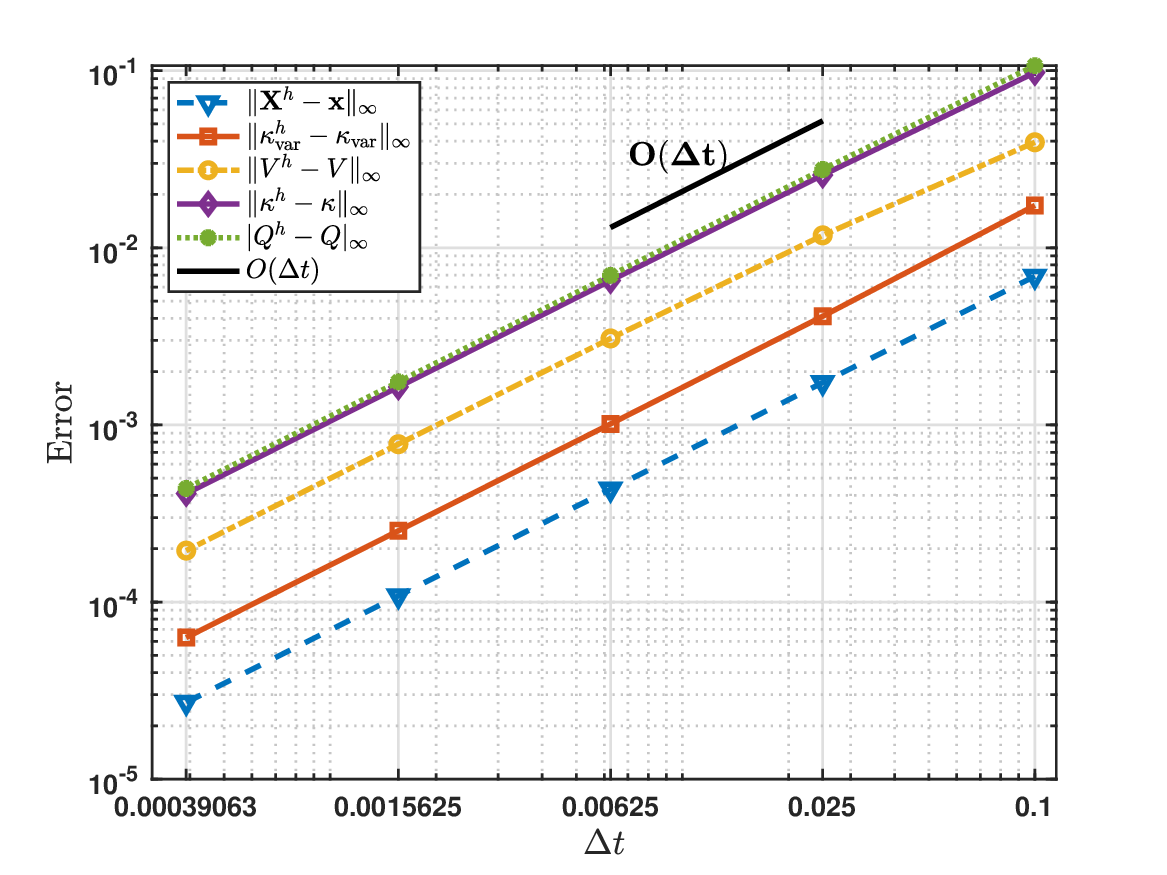}
	\includegraphics[width=0.32\linewidth]{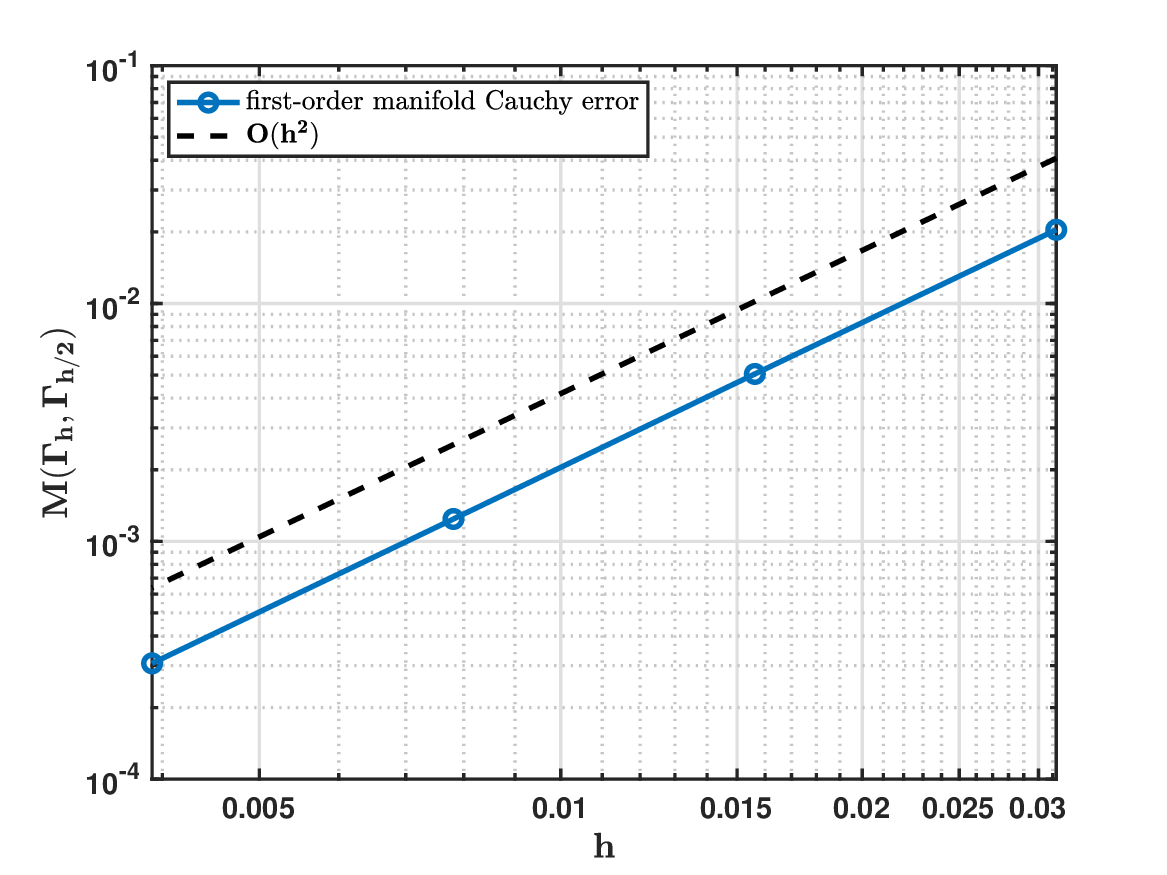}
	\caption{Example 1: convergence of the first-order scheme. Left and middle:
		errors for the circle with uniform and nonuniform initial nodal distributions. Right: manifold Cauchy error for the ellipse, with $\Delta t=O(h^2)$.}
	\label{f1}
\end{figure}

\begin{figure}[htbp]
	\centering
	\includegraphics[width=0.32\linewidth]{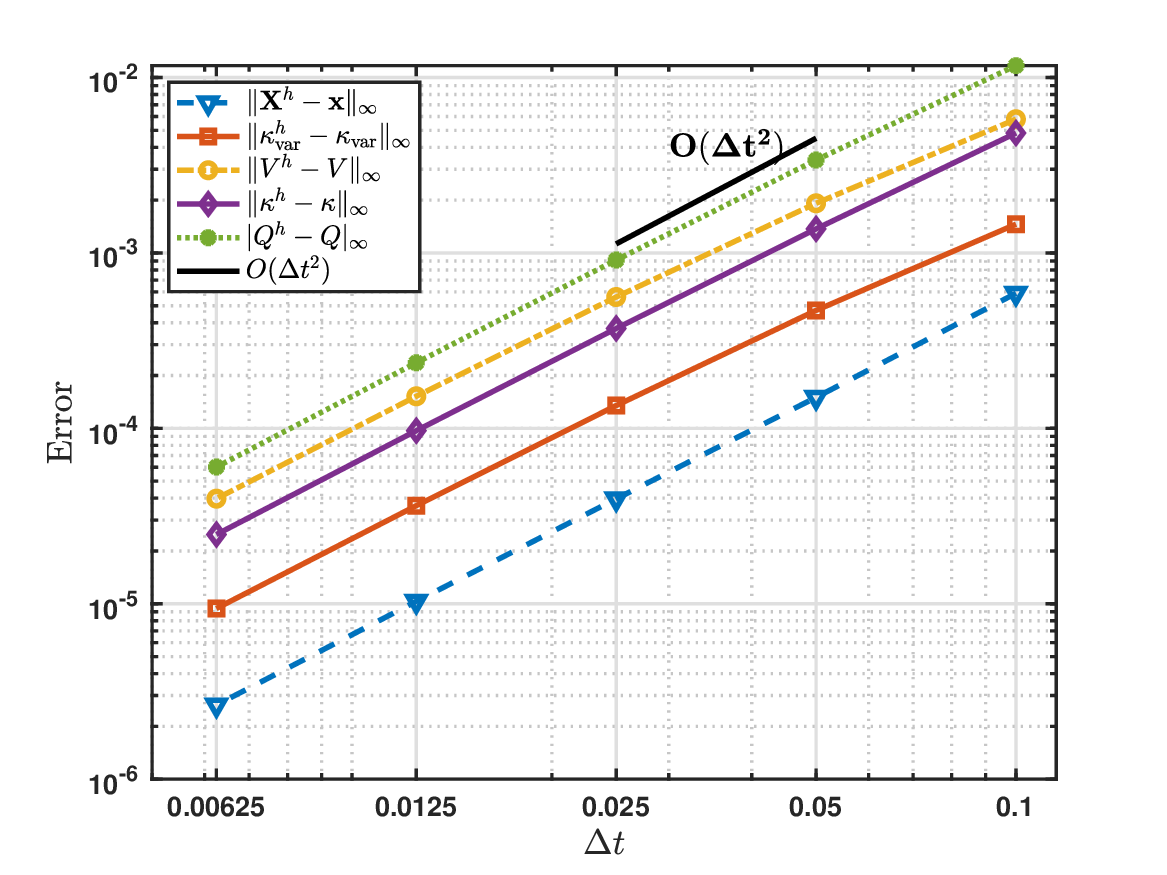}
	\includegraphics[width=0.32\linewidth]{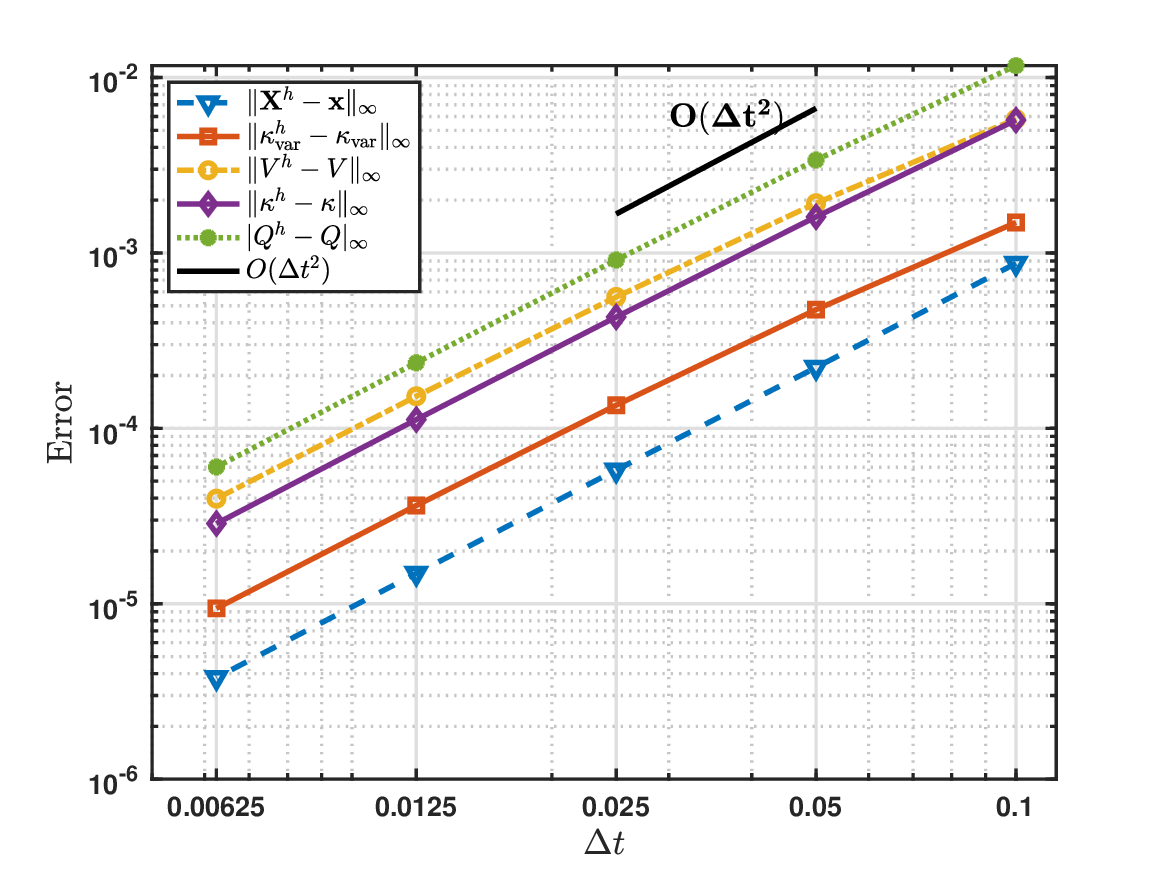}
	\includegraphics[width=0.32\linewidth]{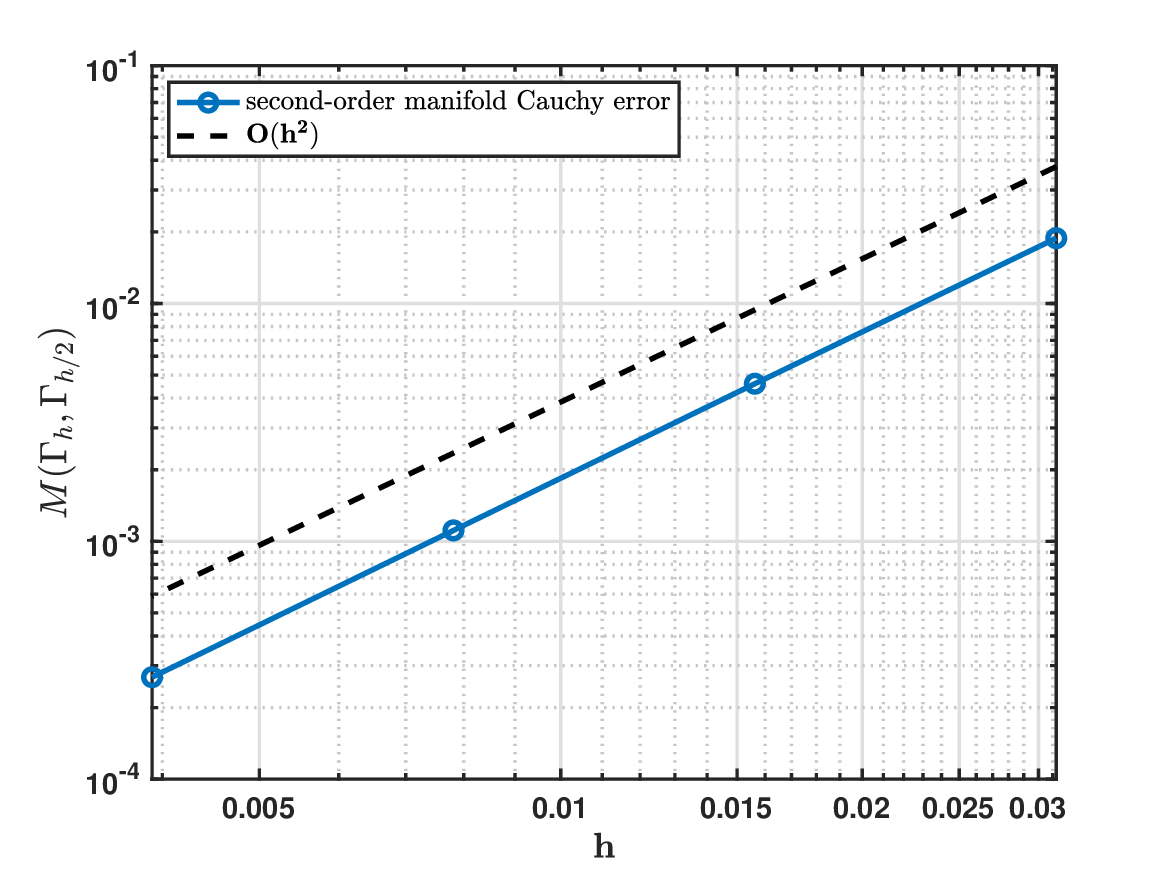}
	\caption{Example 1: convergence of the second-order scheme. Left and middle: errors for the circle with uniform and nonuniform initial nodal distributions. Right: manifold Cauchy error for
		the ellipse, with $\Delta t=O(h)$.}
	\label{f2}
\end{figure}

The numerical errors and convergence rates for the first- and second-order schemes are reported in Figs.~\ref{f1} and~\ref{f2}, respectively. For the expanding circle, the first-order scheme exhibits approximately first-order convergence in $\Delta t$, whereas the second-order scheme exhibits approximately second-order convergence. For the initially elliptical curve, the manifold Cauchy errors decay at approximately second order in $h$. These results are consistent with the expected first- and second-order temporal accuracy of the respective schemes, together with second-order spatial accuracy.


\subsection{Local solvability of the global auxiliary equation}
\paragraph{\textbf{Example 2: Local existence and uniqueness}}
We examine the sufficient conditions established in
Theorems~\ref{thm:Q_first_order} and~\ref{thm:Q_second_order}.
The initial curve is a unit circle with a uniform nodal
distribution, and we set $N=256$, $\lambda=0$, and $T=1$.
For the first-order scheme, the scalar function $F_m$ is a
polynomial of degree at most four, whereas the second-order function $F_m^{(2)}$ contains additional nonlinear geometric terms. For the first-order scheme, we evaluate $C_m^{(1)}=
\frac{4\delta_m|F_m''(1)|}{\sigma_m^2}
+\frac{8\delta_m^2|F_m^{(3)}(1)|}{\sigma_m^3}
+\frac{192a_{4,m}\delta_m^3}{\sigma_m^4}$ and check whether $C_m^{(1)}\le1$. For the second-order scheme, we evaluate $	C_m^{(2)}=\frac{2\delta_mL_m}{\sigma_m^2},~
L_m=\sup_{q\in I_m}|(F_m^{(2)})''(q)|$ and check whether $C_m^{(2)}\le1/2$, together with the edge nondegeneracy assumption in
Theorem~\ref{thm:Q_second_order}.

\begin{figure}[htbp]
	\centering
	\includegraphics[width=0.46\linewidth,height=4cm]{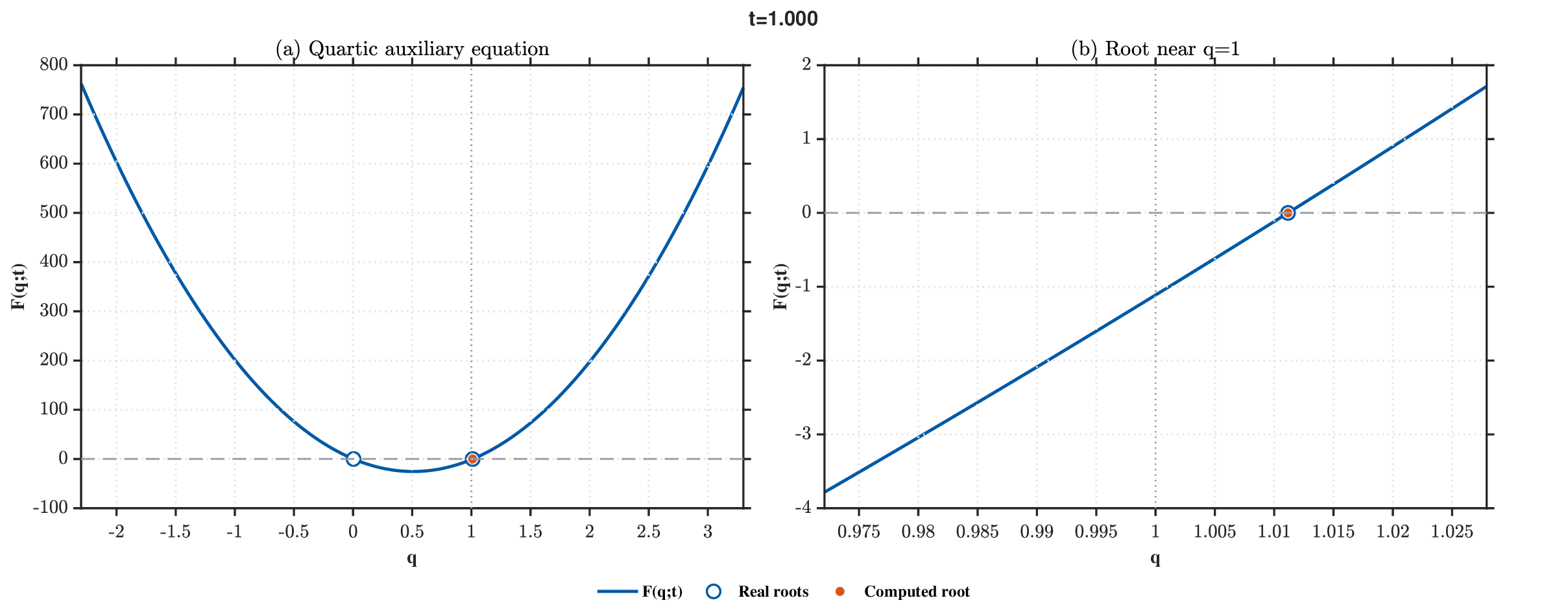}
	\includegraphics[width=0.46\linewidth,height=4cm]{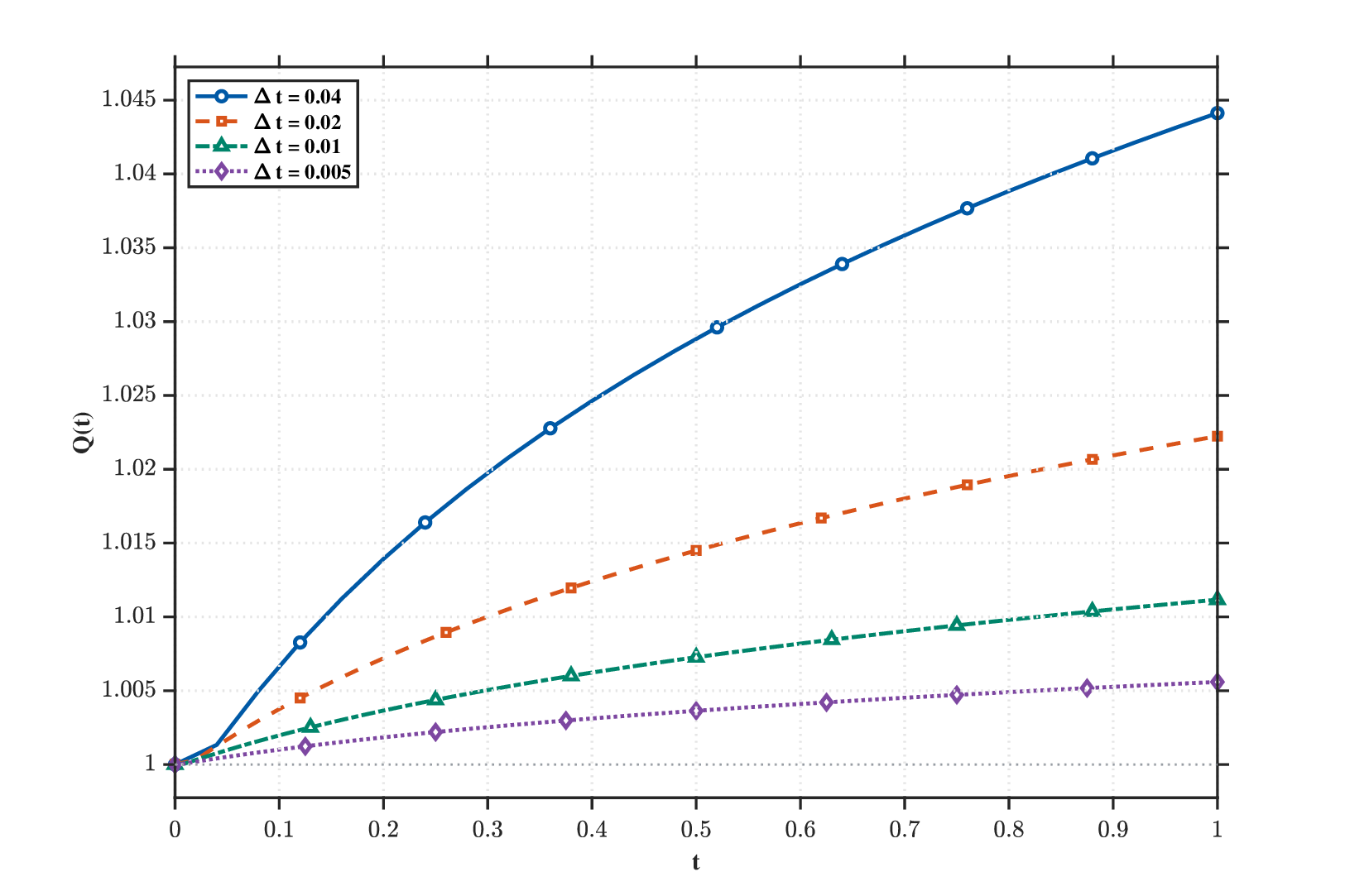}
	\caption{
		Example 2: global auxiliary equation for the first-order
		scheme.	Left: the nonlinear function and its computed root near
		$q=1$. Right: time histories of $Q$ for
		$\Delta t=0.04,0.02,0.01,0.005$.}	\label{eq1}
\end{figure}

\begin{figure}[htbp]
	\centering
	\includegraphics[width=0.8\linewidth,height=8cm]{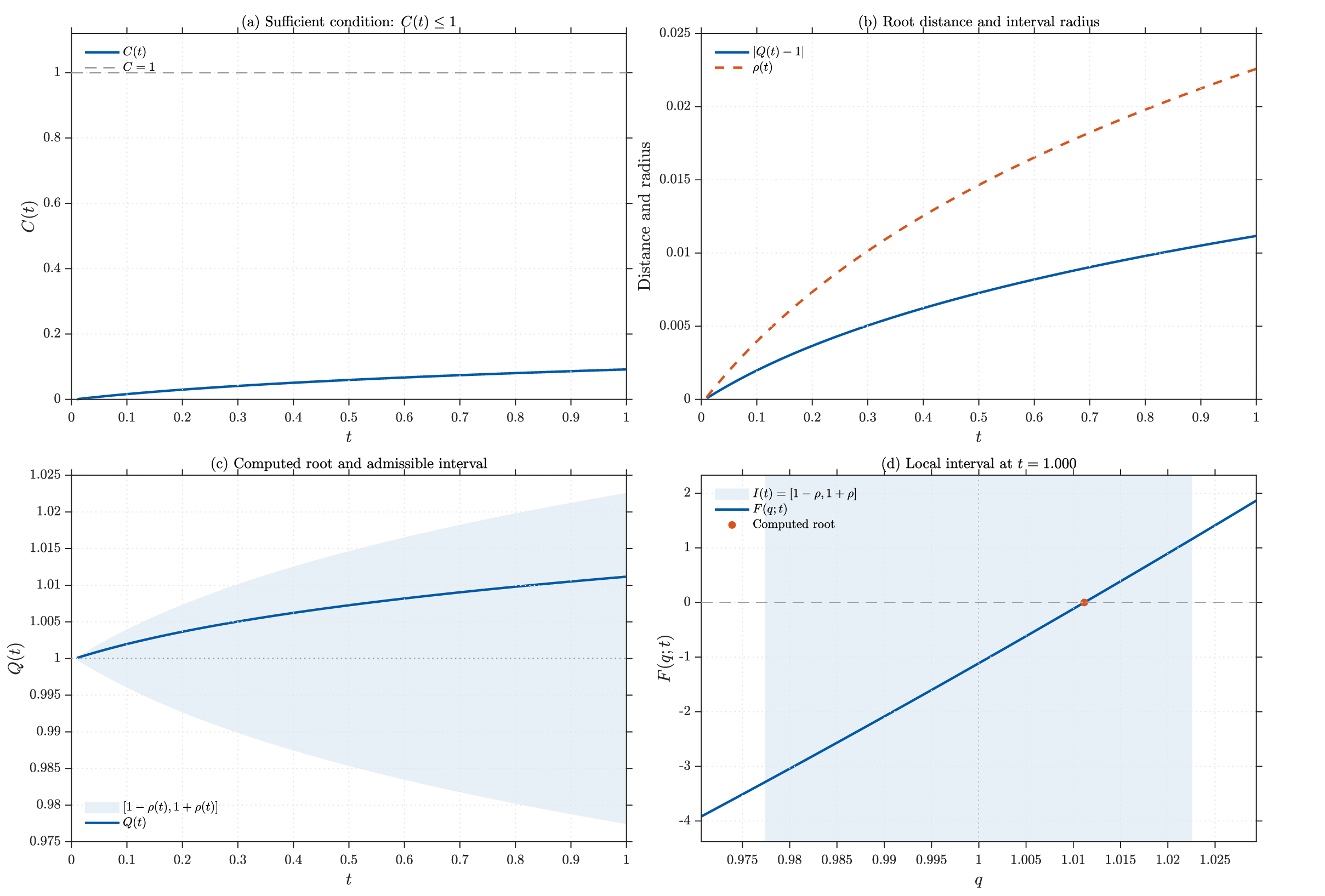}
	\caption{	Example 2: local-solvability diagnostics for the first-order	scheme.	From left to right and top to bottom:
		the sufficient-condition quantity $C_m^{(1)}$;
		the root deviation $|Q^{m+1}-1|$ and interval radius $\rho_m$;
		the computed roots and the intervals $I_m$;
		and a local view of the nonlinear function near the computed
		root.}
	\label{eq2}
\end{figure}

\begin{figure}[htbp]
	\centering
	\includegraphics[width=0.46\linewidth,height=4cm]{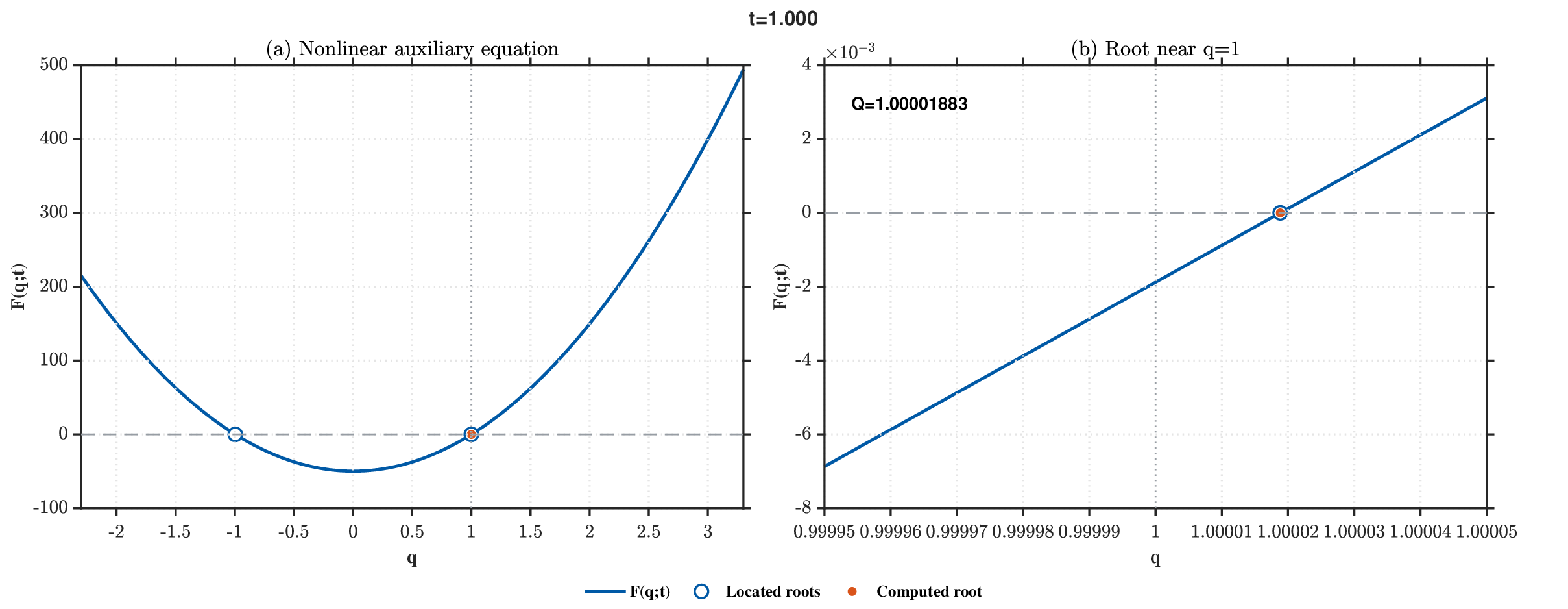}
	\includegraphics[width=0.46\linewidth,height=4cm]{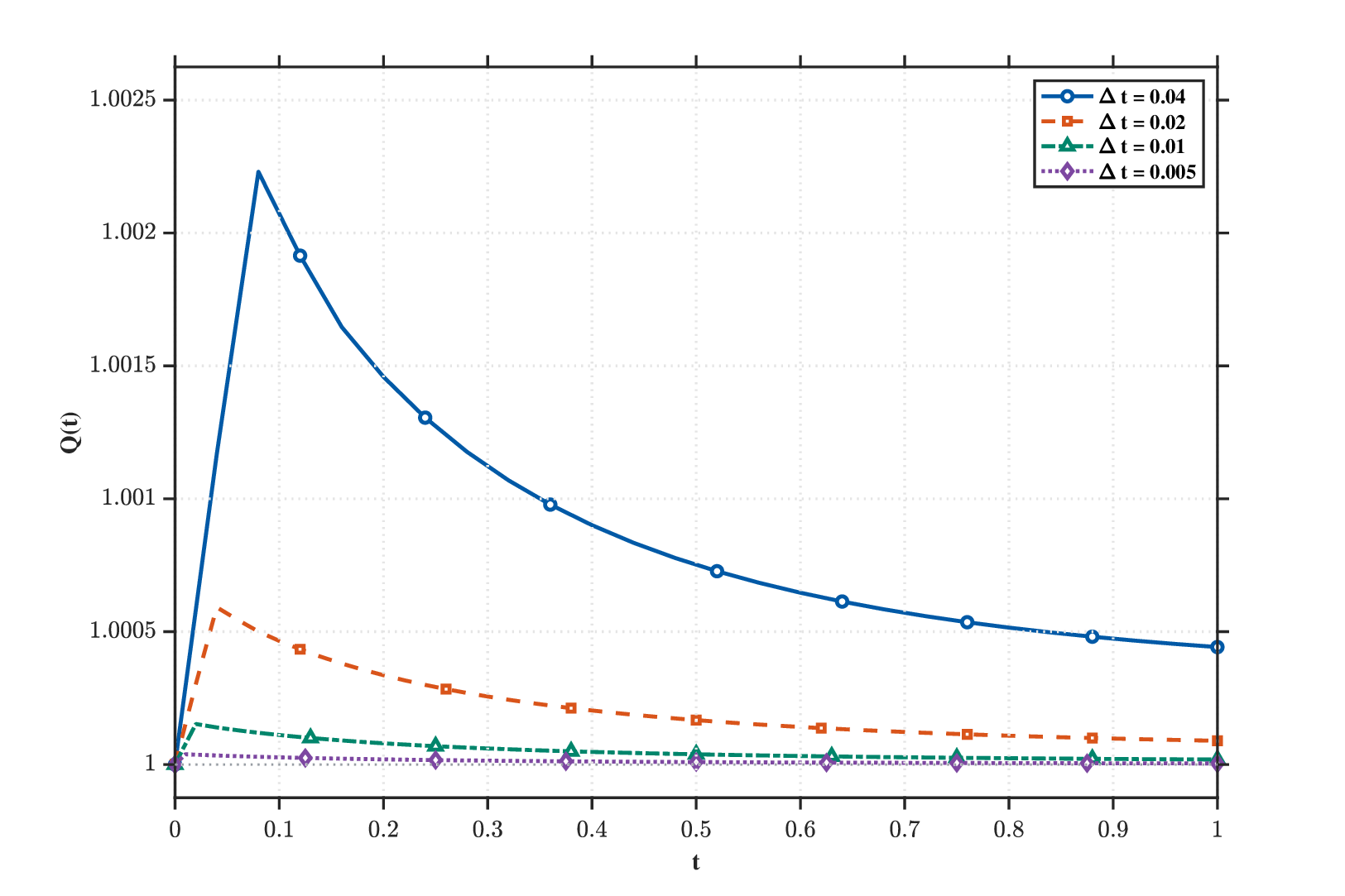}
	\caption{Example 2: global auxiliary equation for the second-order scheme. Left: the nonlinear function and its computed root near	$q=1$. Right: time histories of $Q$ for	$\Delta t=0.04,0.02,0.01,0.005$.}	\label{eq12}
\end{figure}

\begin{figure}[htbp]
	\centering
	\includegraphics[width=0.8\linewidth,height=8cm]{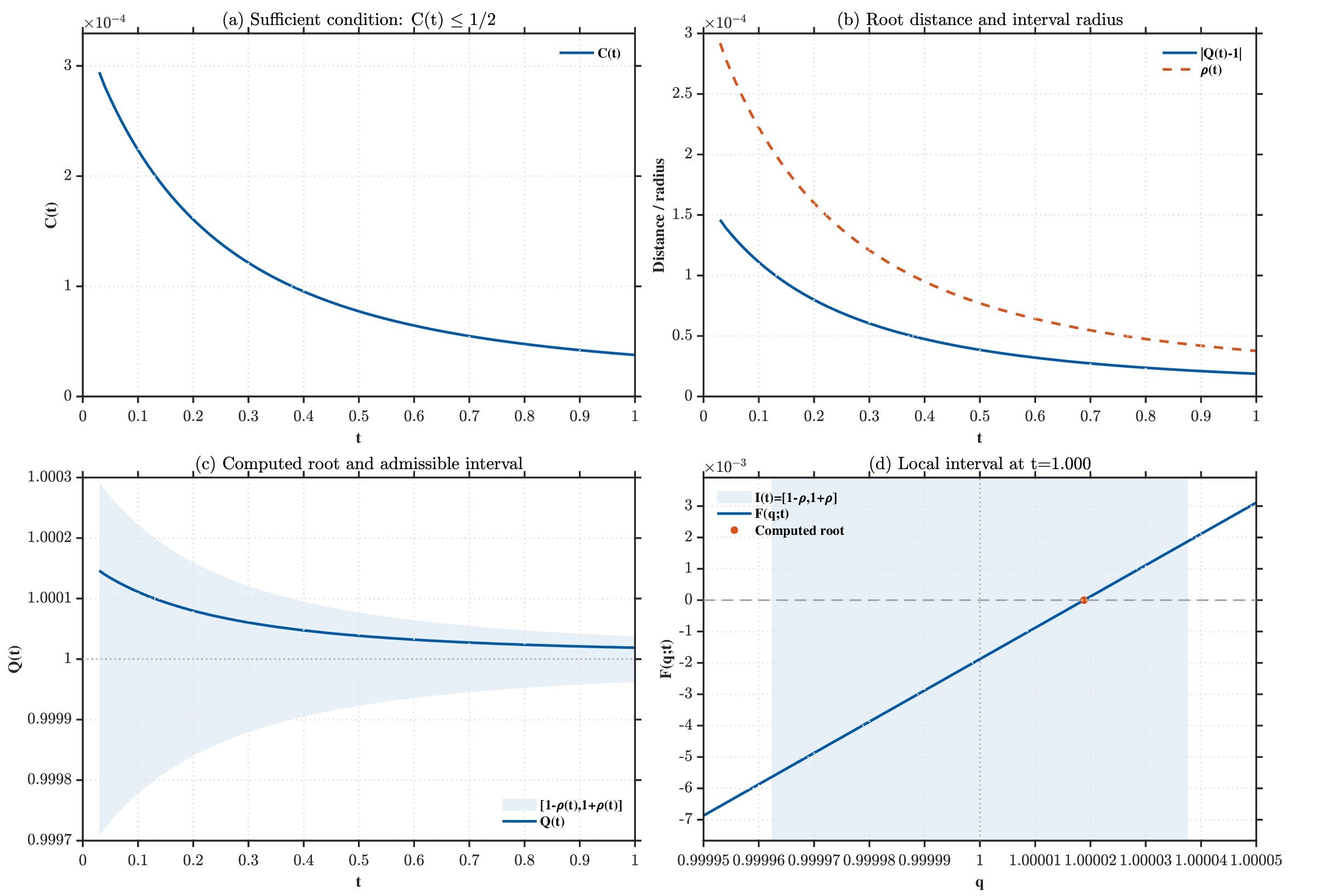}
	\caption{
		Example 2: local-solvability diagnostics for the second-order scheme.
		From left to right and top to bottom:
		the sufficient-condition quantity $C_m^{(2)}$;
		the root deviation $|Q^{m+1}-1|$ and interval radius $\rho_m$;
		the computed roots and the intervals $I_m$;
		and a local view of the nonlinear function near the computed
		root.}
	\label{eq22}
\end{figure}

The scalar nonlinear functions, computed roots, and quantities
corresponding to the local-solvability conditions are presented in
Figs.~\ref{eq1}--\ref{eq22}. For $\Delta t=10^{-2}$, we compare $C_m^{(1)}$
and $C_m^{(2)}$ with their respective thresholds and
$|Q^{m+1}-1|$ with the interval radius $\rho_m$.
Throughout the reported computations, both sufficient
conditions are satisfied, and the computed roots remain
within the corresponding intervals $I_m$.
These observations provide numerical support for the local
existence and uniqueness results.
We further plot the time histories of $Q$ for
$\Delta t=0.04,0.02,0.01,0.005$ to examine the dependence
of the computed auxiliary variable on the time-step size.


\subsection{Computational cost}
\paragraph{\textbf{Example 3: Comparison with a nonlinear Picard scheme}}
We assess the computational efficiency of the proposed schemes
by comparing them with the nonlinear energy-dissipative scheme
introduced in \cite{garcke2025stable}, solved by Picard iteration.
The benchmark is an expanding circle with a nonuniform initial
nodal distribution and $\lambda=0$. For comparison at the same resolution, all three schemes are
first tested on identical spatial meshes and with the same
time steps, using $\Delta t=O(h^2)$.
The numerical errors, experimental orders of convergence
with respect to $h$, maximum iteration counts, and CPU times
are reported in Table~\ref{tab:circle_comparison}(a)--(c).
We further test the second-order scheme with
$\Delta t=O(h)$ to assess its computational efficiency
when larger time steps are used.
The corresponding results are presented in
Table~\ref{tab:circle_comparison}(d).

As shown in Table~\ref{tab:circle_comparison}, all three schemes
exhibit approximately second-order convergence with respect
to $h$.
Compared with the nonlinear Picard scheme, the proposed
first-order scheme requires fewer iterations and achieves
lower computational costs. The second-order scheme also requires only a few
scalar nonlinear iterations per time step. Its higher temporal accuracy allows the use of larger time steps, substantially reducing the
computational cost while retaining comparable accuracy.
These results illustrate the computational efficiency of the
proposed first- and second-order schemes.

\begin{table}[htbp]
	\centering
	\caption{Numerical errors, experimental orders of convergence (EOCs), maximum iteration counts and CPU times  for an expanding circle with a nonuniform initial mesh and $\lambda=0$. The nonlinear Picard scheme and the proposed first- and
		second-order schemes are compared using $\Delta t=O(h^2)$
		in (a)--(c). The second-order scheme is additionally tested
		with $\Delta t=O(h)$ in (d).}
	\label{tab:circle_comparison}
	\resizebox{\textwidth}{!}{
		\begin{tabular}{cc|cc|cc|cc|cc|cc}
			\hline
			$N$ & $\Delta t$ & $e_{\boldsymbol X}$ & EOC & $e_{\varkappa}$ & EOC & $e_V$ & EOC & $e_\kappa$ & EOC & $N_{\mathrm{iter}}^{\max}$ & CPU (s)\\
			\hline\hline
			\multicolumn{12}{c}{\textbf{(a) Nonlinear Picard scheme \cite{garcke2025stable}}}\\
			\hline
			32  & $4.00\times10^{-2}$ & $4.40\times10^{-3}$ & --    & $4.90\times10^{-3}$ & --    & $2.53\times10^{-2}$ & --    & $4.41\times10^{-2}$ & --     & 7 & 0.0325\\
			64  & $1.00\times10^{-2}$ & $1.10\times10^{-3}$ & 2.00  & $1.20\times10^{-3}$ & 2.00  & $6.70\times10^{-3}$ & 1.92  & $1.12\times10^{-2}$ & 1.98  & 5 & 0.121\\
			128 & $2.50\times10^{-3}$ & $2.73\times10^{-4}$ & 2.00  & $3.07\times10^{-4}$ & 2.00  & $1.70\times10^{-3}$ & 1.98  & $2.80\times10^{-3}$ & 1.99   & 5 & 0.762\\
			256 & $6.25\times10^{-4}$ & $6.82\times10^{-5}$ & 2.00  & $7.67\times10^{-5}$ & 2.00  & $4.27\times10^{-4}$ & 1.99  & $7.01\times10^{-4}$ & 2.00  & 4 & 5.88\\
			512 & $1.56\times10^{-4}$ & $1.71\times10^{-5}$ & 2.00  & $1.92\times10^{-5}$ & 2.00  & $1.07\times10^{-4}$ & 2.00  & $1.75\times10^{-4}$ & 2.00   & 4 & 49.5\\
			\hline
			\multicolumn{12}{c}{\textbf{(b) First-order  scheme}}\\
			\hline
			32  & $4.00\times10^{-2}$& $2.10\times10^{-3}$ & --& $6.80\times10^{-3}$ & --& $1.81\times10^{-2}$ & --& $4.38\times10^{-2}$ & --	 & 4& 0.0862 \\
			64  & $1.00\times10^{-2}$	& $5.39\times10^{-4}$ & 1.99& $1.70\times10^{-3}$ & 2.03& $4.90\times10^{-3}$ & 1.89& $1.12\times10^{-2}$ & 1.97& 4& 0.0997 \\
			128 & $2.50\times10^{-3}$
			& $1.35\times10^{-4}$ & 2.00& $4.17\times10^{-4}$ & 2.01& $1.20\times10^{-3}$ & 1.97& $2.80\times10^{-3}$ & 1.99 &3& 0.4481 \\
			256 & $6.25\times10^{-4}$& $3.38\times10^{-5}$ & 2.00& $1.04\times10^{-4}$ & 2.00& $3.12\times10^{-4}$ & 1.99& $7.01\times10^{-4}$ & 2.00 &3 & 2.9354 \\
			512 & $1.56\times10^{-4}$& $8.44\times10^{-6}$ & 2.00& $2.60\times10^{-5}$ & 2.00& $7.81\times10^{-5}$ & 2.00& $1.75\times10^{-4}$ & 2.00 & 2 & 21.7246 \\
			\hline
			\multicolumn{12}{c}{\textbf{(c) Second-order scheme}}\\
			\hline
			32  & $4.00\times10^{-2}$ & $1.40\times10^{-3}$ & --    & $4.03\times10^{-4}$ & --    & $1.40\times10^{-3}$ & --    & $5.60\times10^{-3}$ & --     & 3 & 0.0195\\
			64  & $1.00\times10^{-2}$ & $3.80\times10^{-4}$ & 1.92  & $8.37\times10^{-5}$ & 2.27  & $1.08\times10^{-4}$ & 3.71  & $1.40\times10^{-3}$ & 1.97   & 3 & 0.1157\\
			128 & $2.50\times10^{-3}$ & $9.61\times10^{-5}$ & 1.98  & $2.21\times10^{-5}$ & 1.92  & $1.87\times10^{-5}$ & 2.53  & $3.64\times10^{-4}$ & 1.99   & 2 & 0.8837\\
			256 & $6.25\times10^{-4}$ & $2.41\times10^{-5}$ & 2.00  & $5.60\times10^{-6}$ & 1.98  & $4.42\times10^{-6}$ & 2.08  & $9.20\times10^{-5}$ & 1.98   & 2 & 6.8937\\
			512 & $1.56\times10^{-4}$ & $6.03\times10^{-6}$ & 2.00  & $1.40\times10^{-6}$ & 2.00  & $1.09\times10^{-6}$ & 2.02  & $2.31\times10^{-5}$ & 2.00  & 2 & 56.3238\\
			\hline
			\hline
			\multicolumn{12}{c}{\textbf{(d) Second-order scheme with large time step}}\\
			\hline
			32  & $4.00\times10^{-2}$ & $1.40\times10^{-3}$ & --    & $4.03\times10^{-4}$ & --    & $1.40\times10^{-3}$ & --    & $5.60\times10^{-3}$ & --    & 3 & 0.0228\\
			64  & $2.00\times10^{-2}$ & $3.67\times10^{-4}$ & 1.98  & $1.00\times10^{-4}$ & 2.01  & $3.85\times10^{-4}$ & 1.87  & $1.40\times10^{-3}$ & 1.98   & 3 & 0.0607\\
			128 & $1.00\times10^{-2}$ & $9.26\times10^{-5}$ & 1.99  & $2.50\times10^{-5}$ & 2.00  & $1.01\times10^{-4}$ & 1.93  & $3.61\times10^{-4}$ & 1.99   & 3 & 0.222\\
			256 & $5.00\times10^{-3}$ & $2.33\times10^{-5}$ & 1.99  & $6.25\times10^{-6}$ & 2.00  & $2.59\times10^{-5}$ & 1.96  & $9.07\times10^{-5}$ & 1.99   & 3 & 0.873\\
			512 & $2.50\times10^{-3}$ & $5.84\times10^{-6}$ & 2.00  & $1.56\times10^{-6}$ & 2.00  & $6.56\times10^{-6}$ & 1.98  & $2.27\times10^{-5}$ & 2.00 & 2 & 3.49\\
			\hline
	\end{tabular}	}
\end{table}

\subsection{Performance of the schemes}
\label{subsec:scheme_performance}
\begin{figure}[htbp]
	\centering
	\begin{minipage}[c]{0.32\linewidth}
		\centering
		\includegraphics[width=0.844\linewidth]{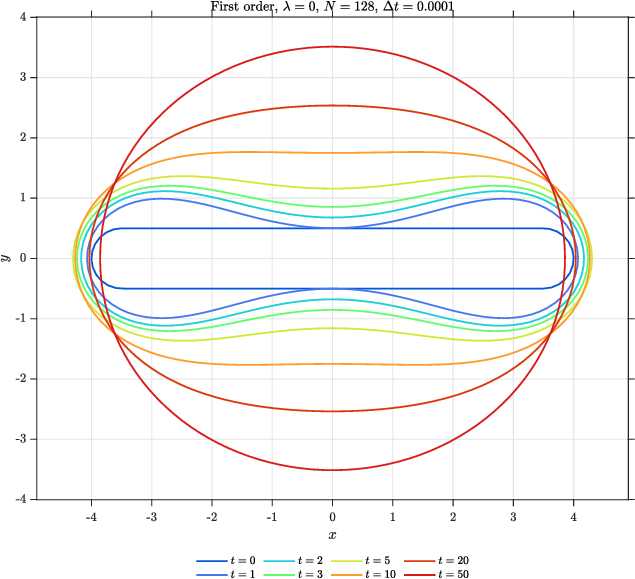}
	\end{minipage}\hfill
	\begin{minipage}[c]{0.32\linewidth}
		\centering
		\includegraphics[width=0.844\linewidth]{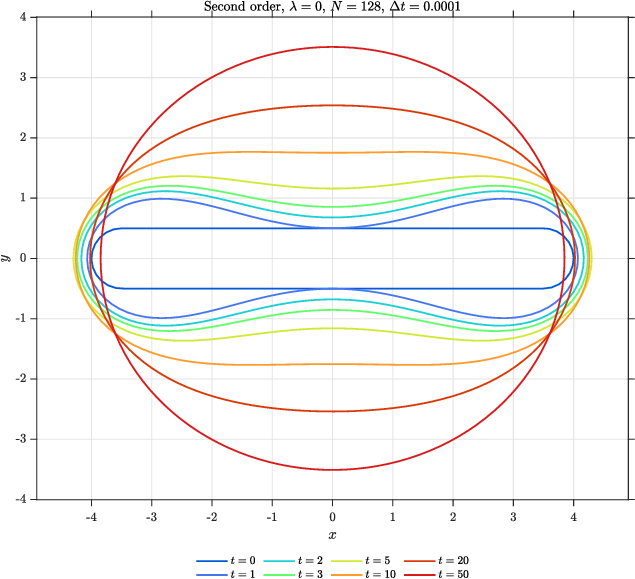}
	\end{minipage}\hfill
	\begin{minipage}[c]{0.32\linewidth}
		\centering
		\includegraphics[width=0.9375\linewidth]{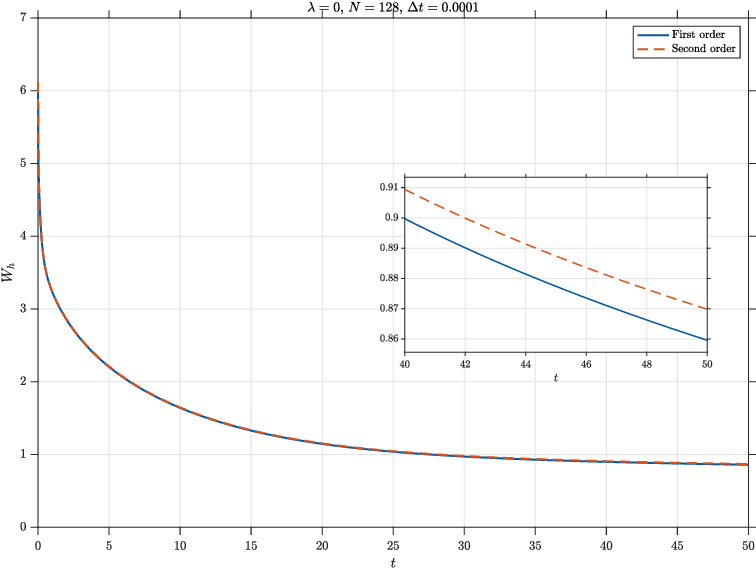}
	\end{minipage}
	\includegraphics[width=0.30\linewidth]{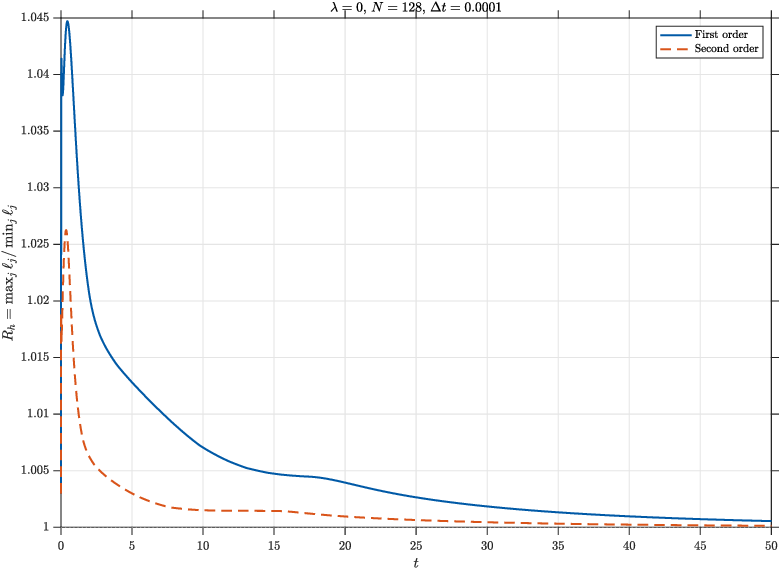}\hfill
	\includegraphics[width=0.30\linewidth]{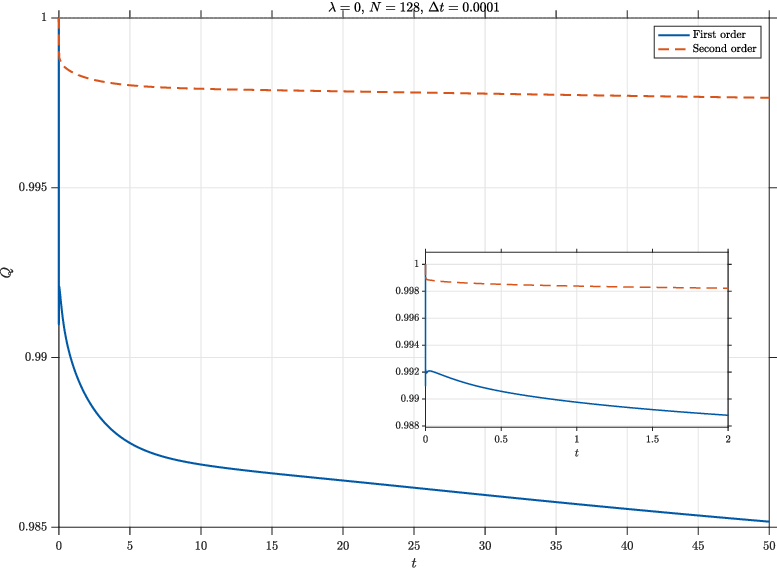}\hfill
	\includegraphics[width=0.30\linewidth]{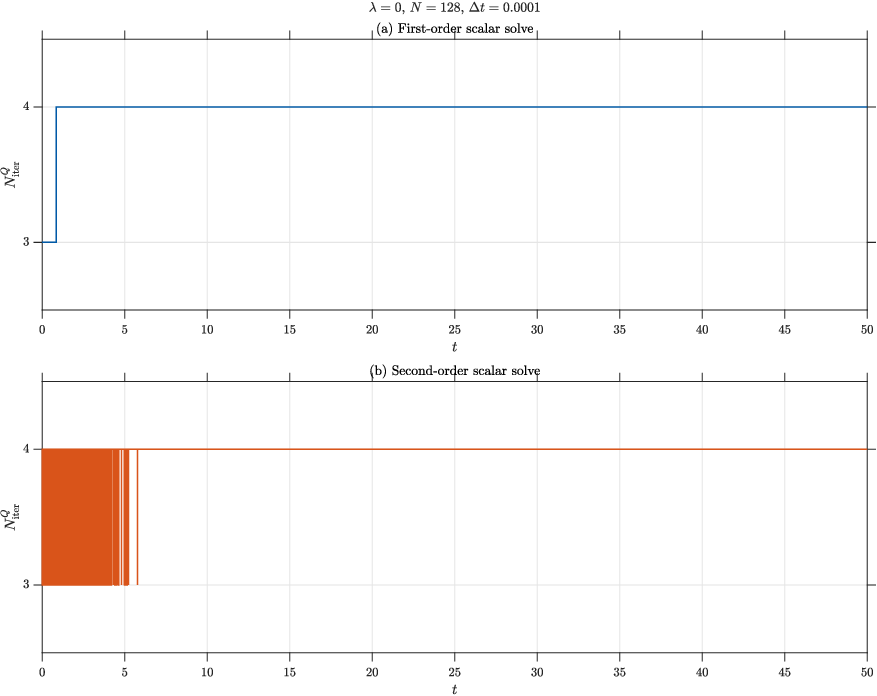}
	\caption{Evolution of an $8\times1$ tube with $\lambda=0$, $N=128$, and $\Delta t=10^{-4}$. 	Top row, from left to right: first- and second-order curve evolution at $t=0,1,2,3,5,10,20,50$, and modified discrete energy.
		Bottom row: mesh ratio, auxiliary variable $Q$, and scalar solver iteration counts.}
	\label{exp3}
\end{figure}

\begin{figure}[htbp]
	\centering
	\begin{minipage}[c]{0.32\linewidth}
		\centering
		\includegraphics[width=\linewidth]{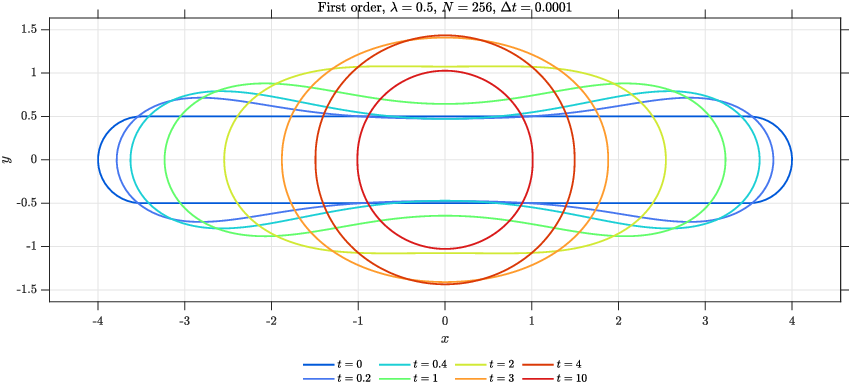}
	\end{minipage}\hfill
	\begin{minipage}[c]{0.32\linewidth}
		\centering
		\includegraphics[width=\linewidth]{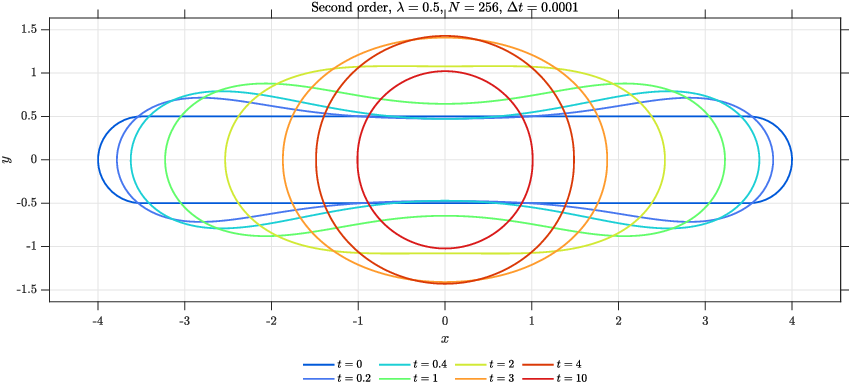}
	\end{minipage}\hfill
	\begin{minipage}[c]{0.30\linewidth}
		\centering
		\includegraphics[width=\linewidth]{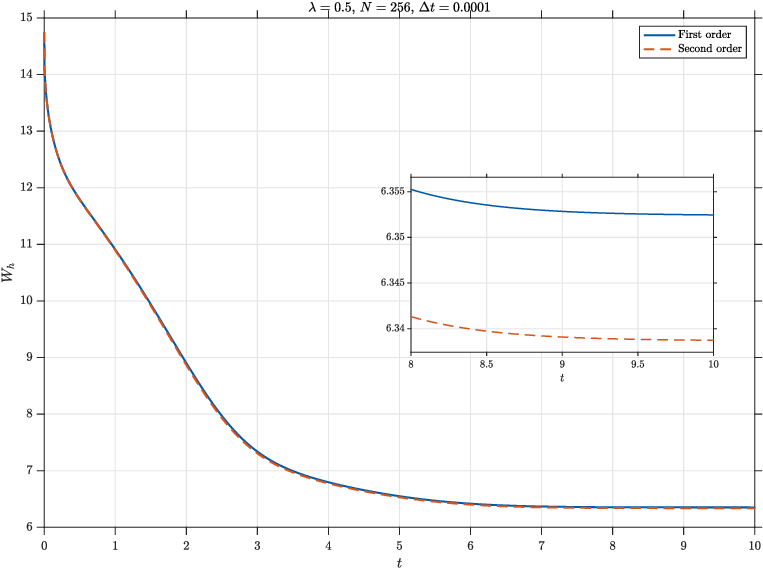}
	\end{minipage}
	
	\par\medskip
	\includegraphics[width=0.30\linewidth]{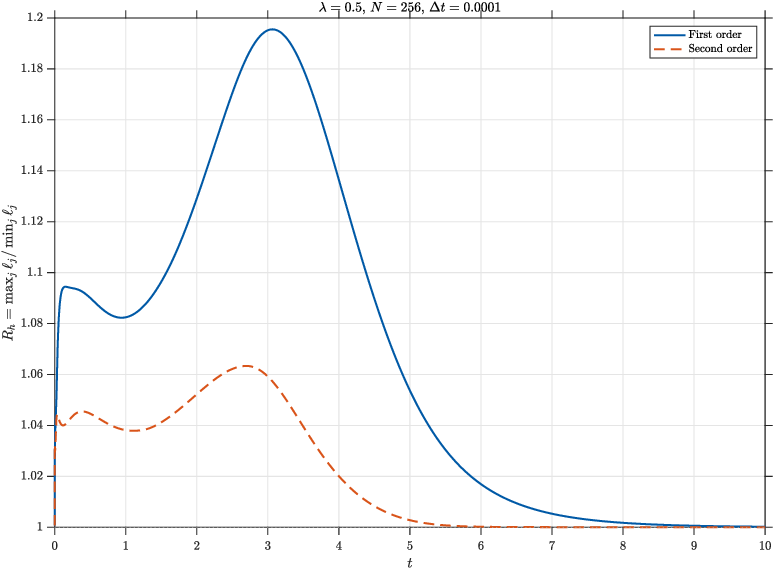}\hfill
	\includegraphics[width=0.30\linewidth]{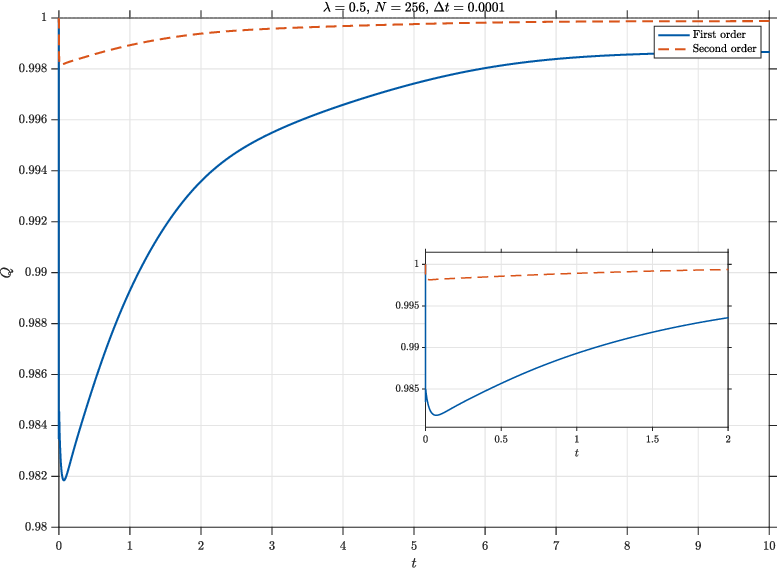}\hfill
	\includegraphics[width=0.30\linewidth]{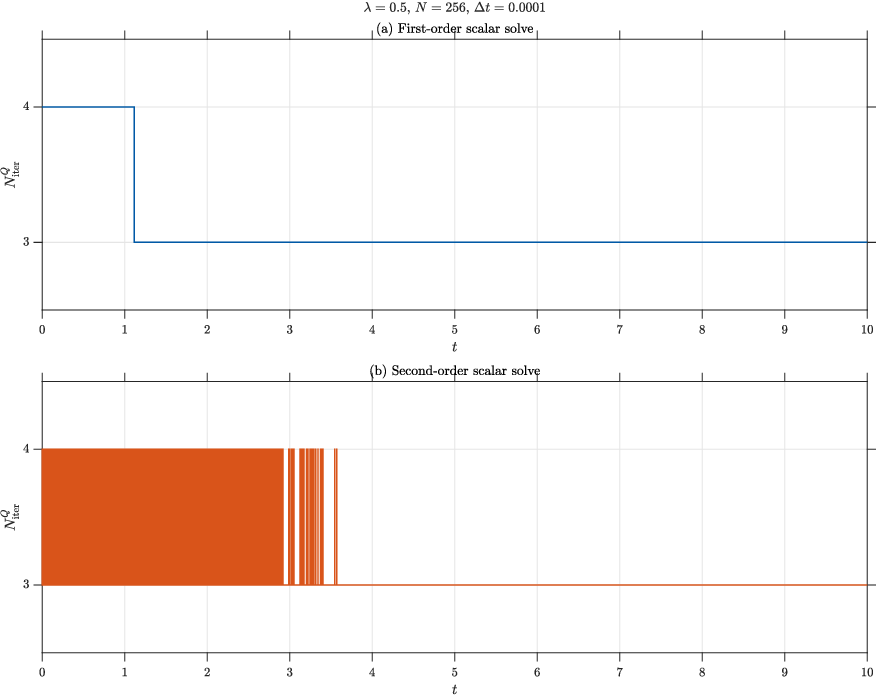}
	\caption{Evolution of an $8\times1$ tube with $\lambda=0.5$, $N=256$, and $\Delta t=10^{-4}$.
		Top row, from left to right: first- and second-order curve evolution at $t=	0,0.2,0.4,0.6,0.8,1,2,10$, and modified discrete energy.
		Bottom row: mesh ratio, auxiliary variable $Q$, and scalar solver iteration counts.}
	\label{exp4}
\end{figure}

\begin{figure}[htbp]
	\centering
	\begin{minipage}[c]{0.32\linewidth}
		\centering
		
		\includegraphics[width=\linewidth,	height=2cm, keepaspectratio]{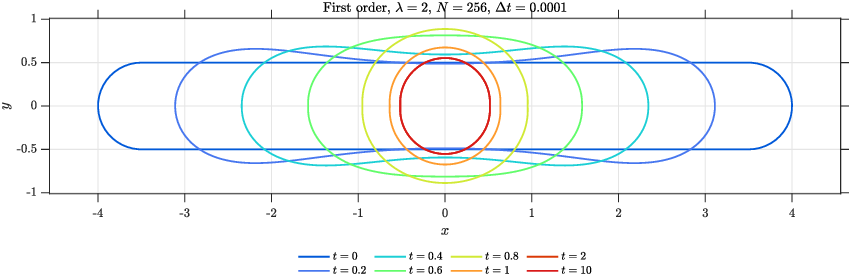}
	\end{minipage}\hfill
	\begin{minipage}[c]{0.32\linewidth}
		\centering
		\includegraphics[width=\linewidth,	height=2cm, keepaspectratio]{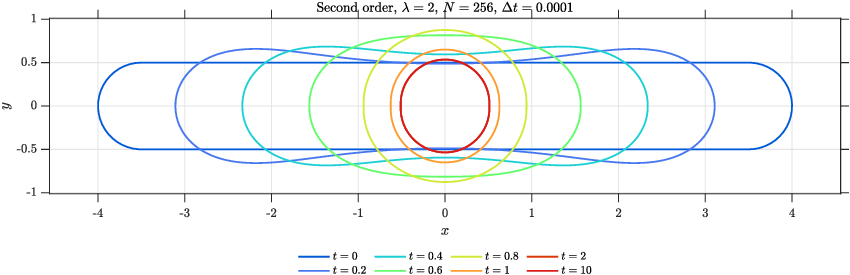}
	\end{minipage}\hfill
	\begin{minipage}[c]{0.30\linewidth}
		\centering
		\includegraphics[width=\linewidth]{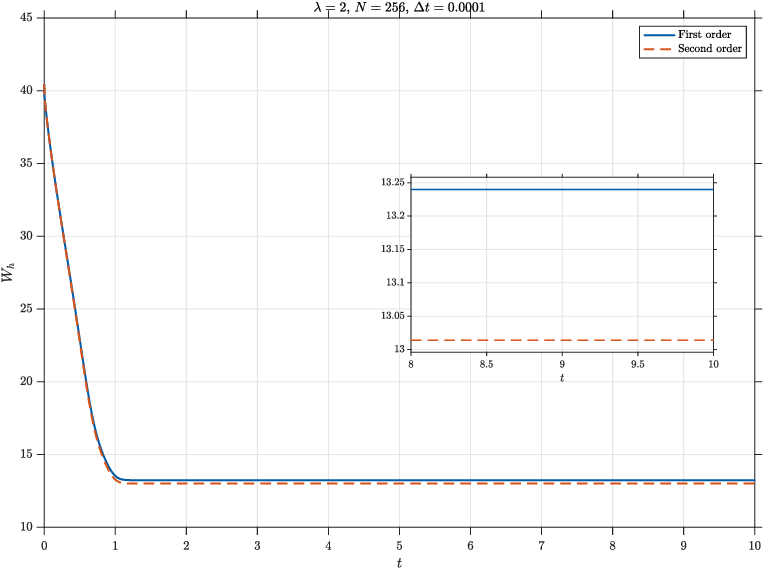}
	\end{minipage}
	
	\par\medskip
	\includegraphics[width=0.30\linewidth]{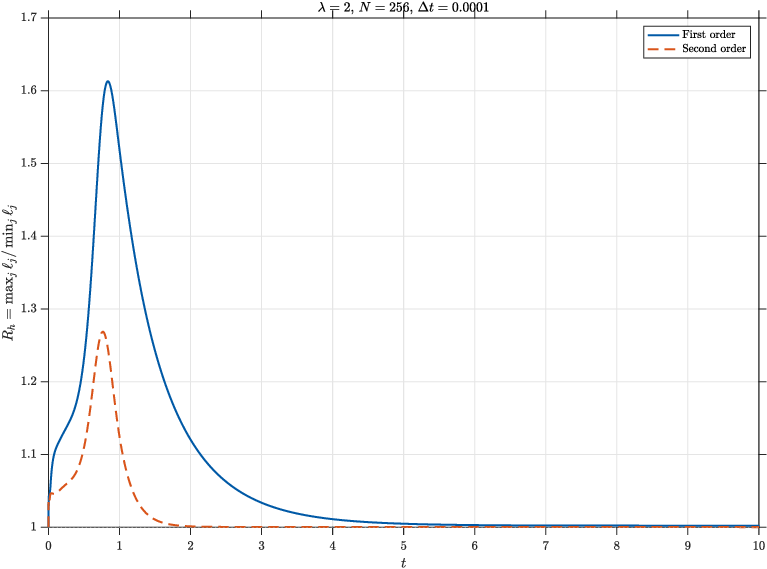}\hfill
	\includegraphics[width=0.30\linewidth]{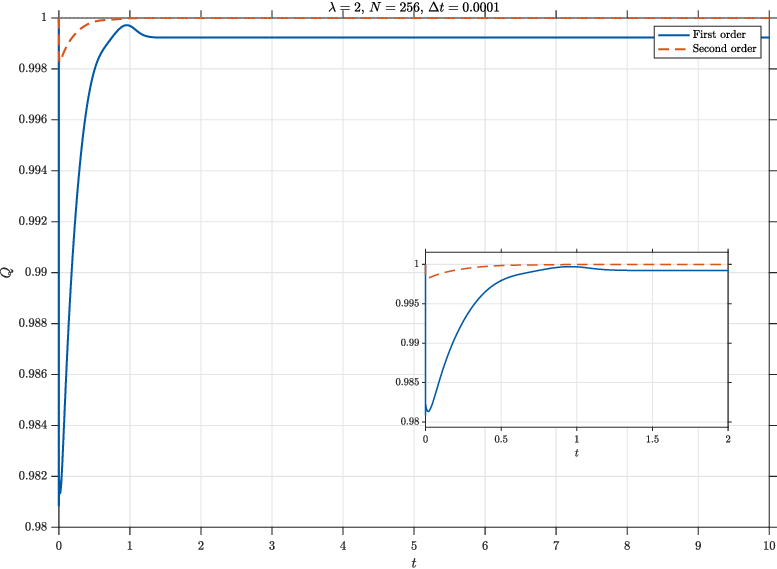}\hfill
	\includegraphics[width=0.30\linewidth]{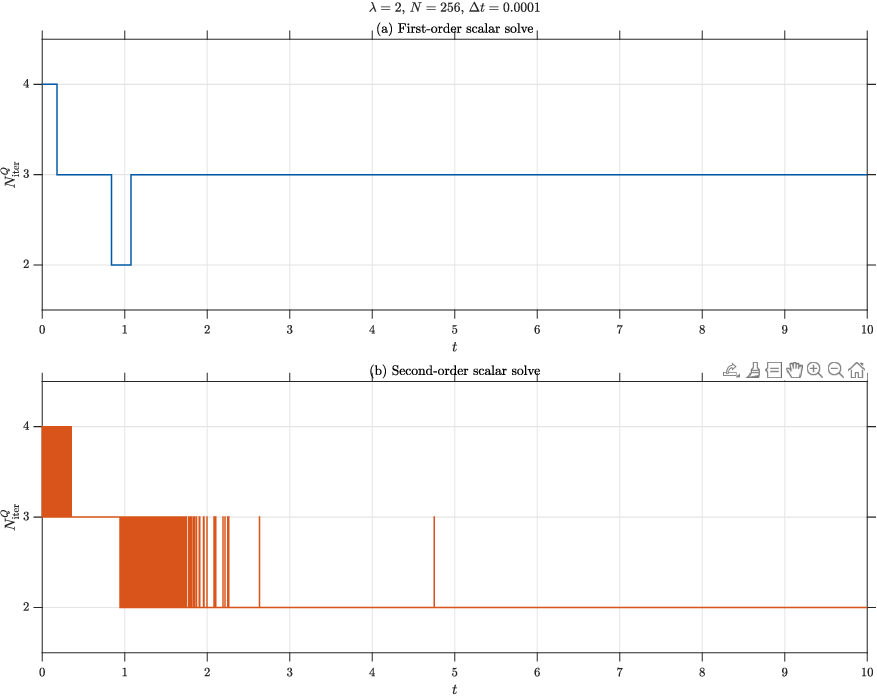}
	\caption{Evolution of an $8\times1$ tube with $\lambda=2$, $N=256$, and $\Delta t=10^{-4}$.
		Top row, from left to right: first- and second-order curve evolution at $t=	0,0.2,0.4,0.6,0.8,1,2,10$, and modified discrete energy.
		Bottom row: mesh ratio, auxiliary variable $Q$, and scalar solver iteration counts.}
	\label{exp5}
\end{figure}

\begin{figure}[htbp]
	\centering
	\begin{minipage}[c]{0.32\linewidth}
		\centering
		\includegraphics[width=\linewidth]{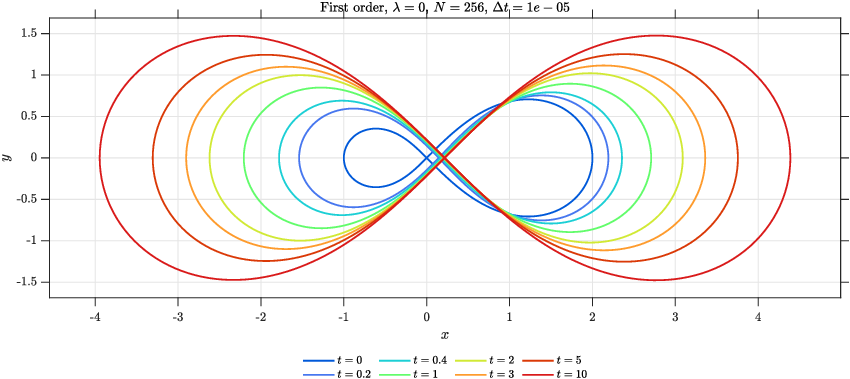}
	\end{minipage}\hfill
	\begin{minipage}[c]{0.32\linewidth}
		\centering
		\includegraphics[width=\linewidth]{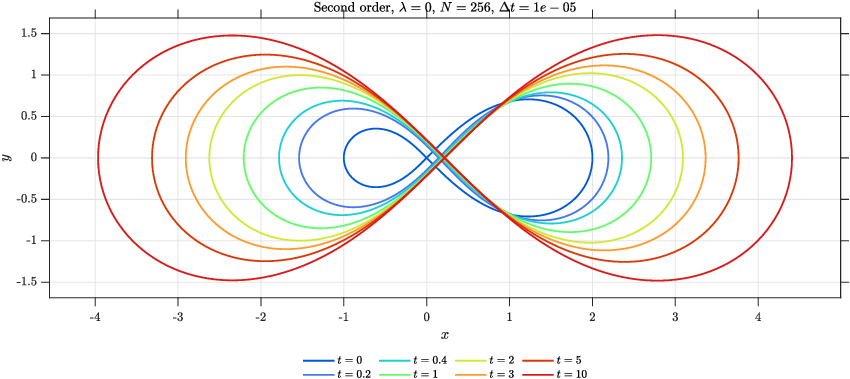}
	\end{minipage}\hfill
	\begin{minipage}[c]{0.30\linewidth}
		\centering
		\includegraphics[width=\linewidth]{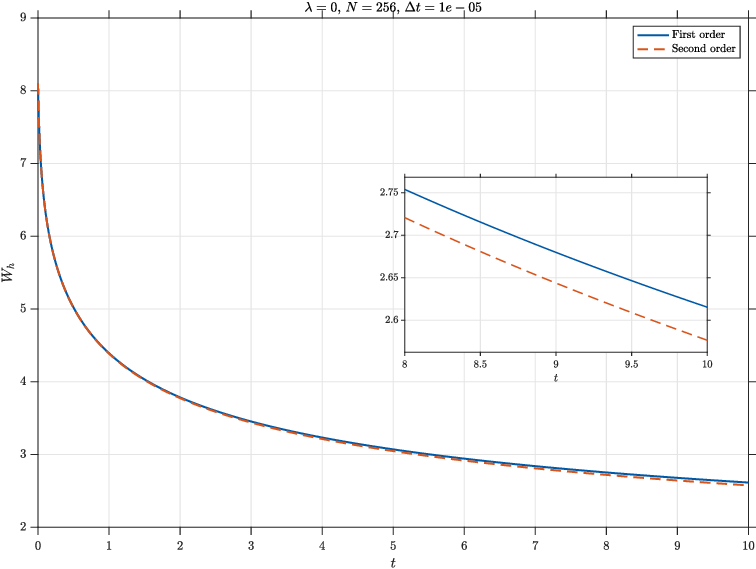}
	\end{minipage}
	
	\par\medskip

	\includegraphics[width=0.30\linewidth]{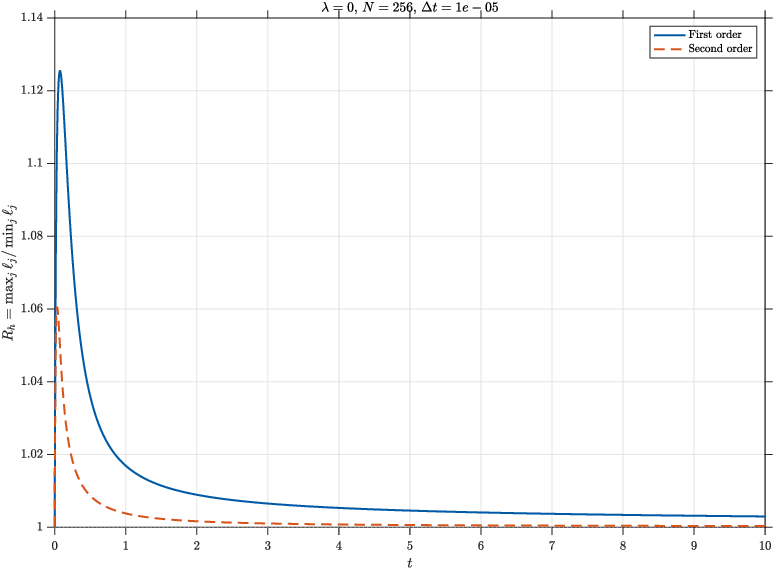}\hfill
	\includegraphics[width=0.30\linewidth]{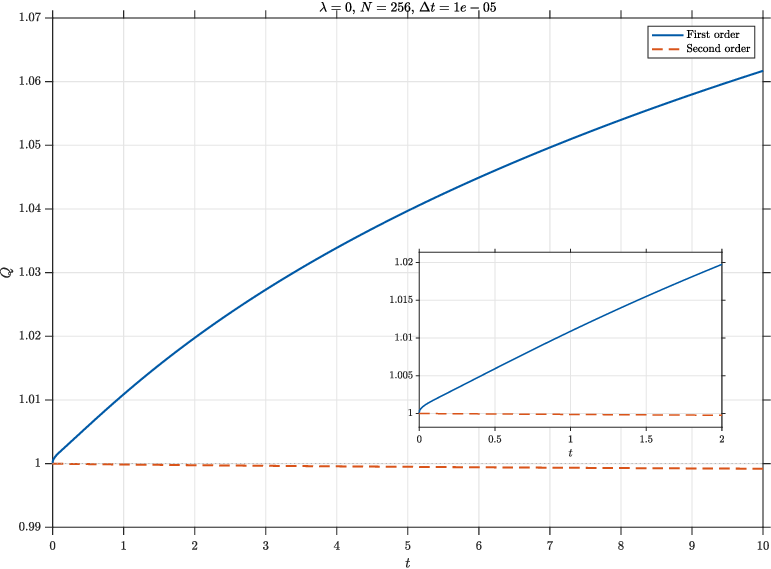}\hfill
	\includegraphics[width=0.30\linewidth]{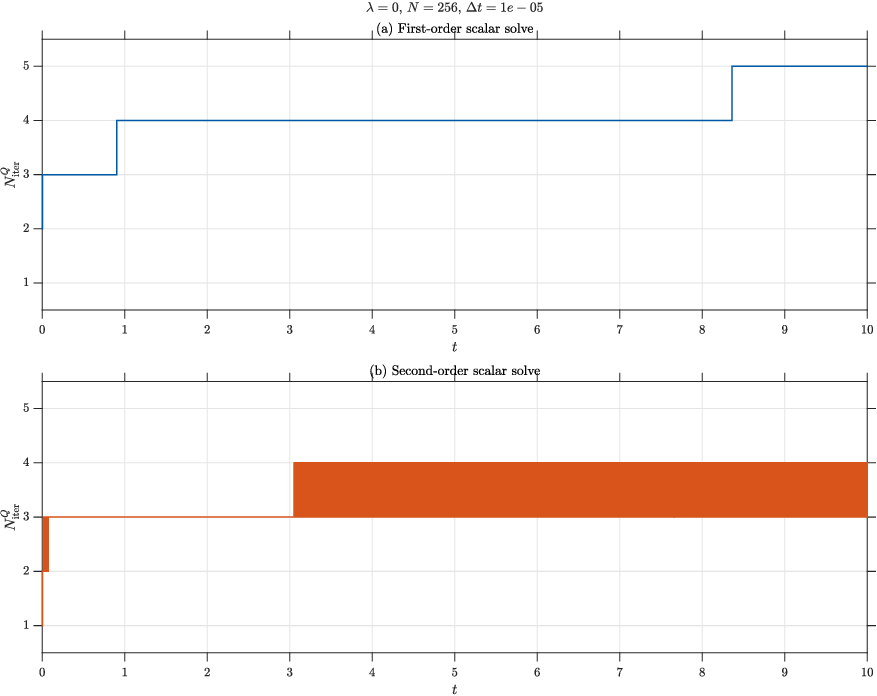}
	
	\caption{Evolution of a lemniscate with $\lambda=0$, $N=256$, and $\Delta t=10^{-5}$.
		Top row, from left to right: first- and second-order curve evolution, and modified discrete energy.
		Bottom row: mesh ratio, auxiliary variable $Q$, and scalar solver iteration counts.}
	\label{exp6}
\end{figure}

\begin{figure}[htbp]
	\centering
	\begin{minipage}[c]{0.32\linewidth}
		\centering
		\includegraphics[width=\linewidth]{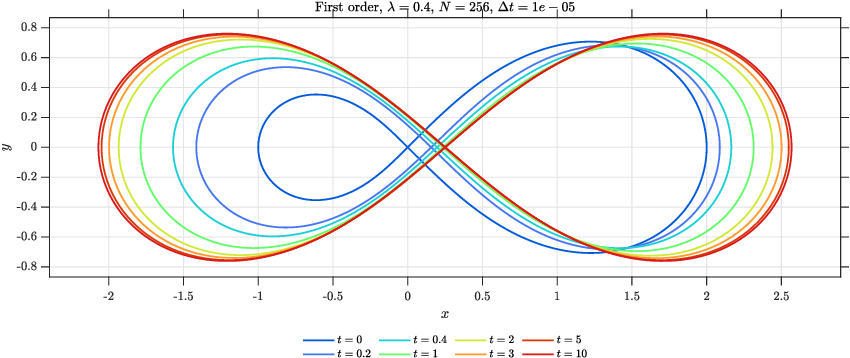}
	\end{minipage}\hfill
	\begin{minipage}[c]{0.32\linewidth}
		\centering
		\includegraphics[width=\linewidth]{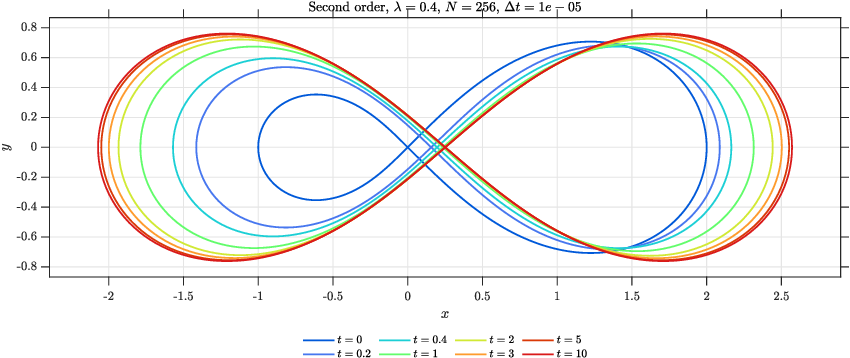}
	\end{minipage}\hfill
	\begin{minipage}[c]{0.30\linewidth}
		\centering
		\includegraphics[width=\linewidth]{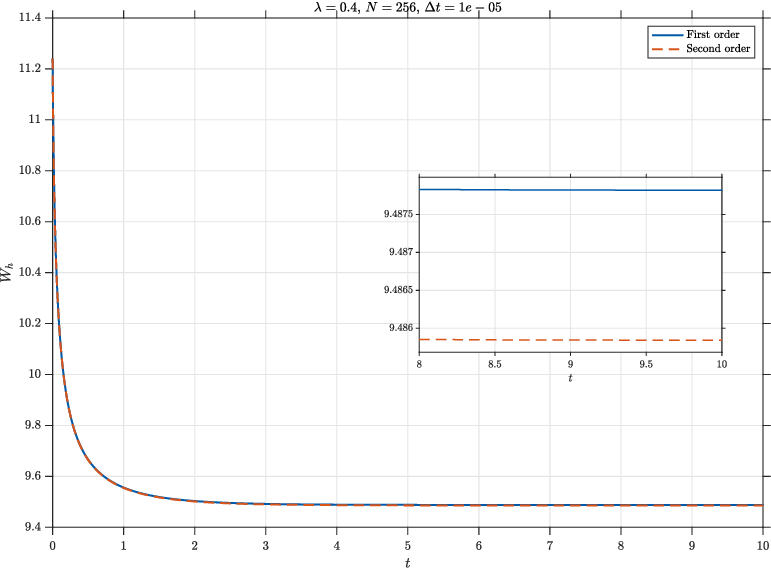}
	\end{minipage}
	
	\par\medskip
	
	\includegraphics[width=0.30\linewidth]{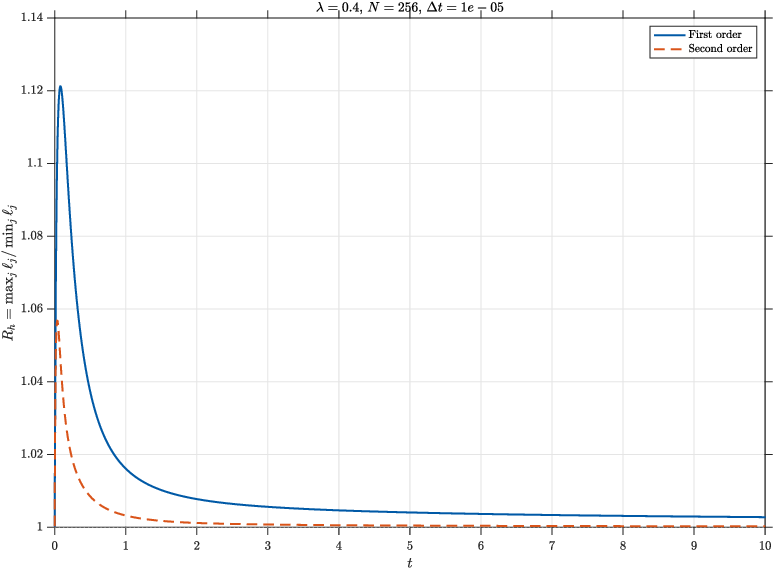}\hfill
	\includegraphics[width=0.30\linewidth]{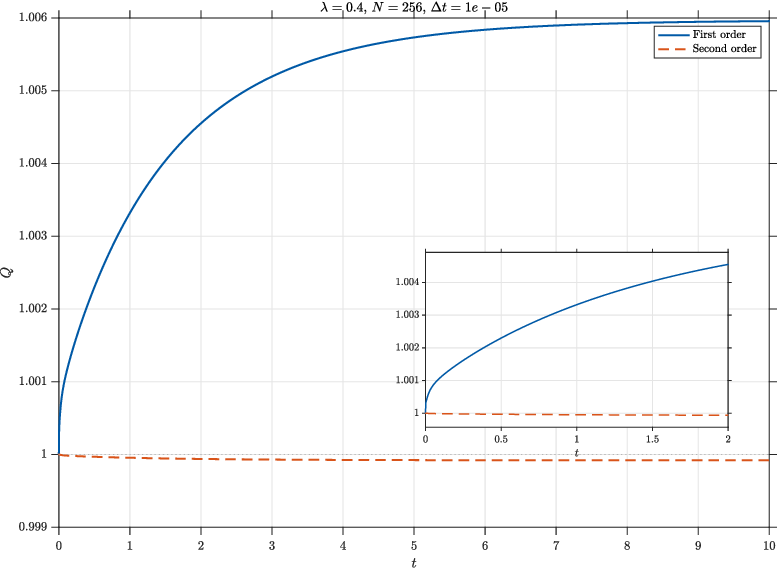}\hfill
	\includegraphics[width=0.30\linewidth]{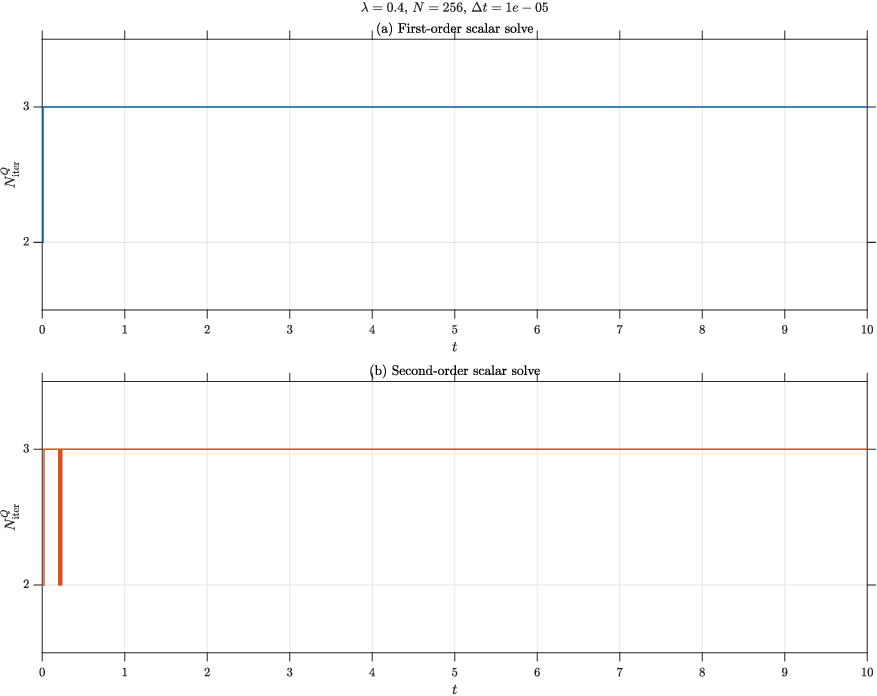}
	
	\caption{Evolution of a lemniscate with $\lambda=0.4$, $N=256$, and $\Delta t=10^{-5}$.
		Top row, from left to right: first- and second-order curve evolution, and modified discrete energy.
		Bottom row: mesh ratio, auxiliary variable $Q$, and scalar solver iteration counts.}
	\label{exp7}
\end{figure}

We finally investigate the performance of the proposed schemes for several representative initial curves, including an elongated tube, a 2:1 lemniscate, and a 3-fold  curve.  To test the mesh quality, we define the mesh ratio
\begin{equation}
	R_h^m:=\frac{\max_{1\le j\le N}|\boldsymbol h_j^m|}{\min_{1\le j\le N}|\boldsymbol h_j^m|}.
	\label{eq:mesh_ratio}
\end{equation}
A value of $R_h^m$ close to one indicates a nearly uniform distribution of mesh points. We also use the modified discrete energy 
to examine the discrete energy-dissipation law. In Figs.~\ref{exp3}--\ref{exp7}, the top row shows, from left to right,
the first-order evolution, the second-order evolution, and the modified
energy. The bottom row compares the mesh ratio, the auxiliary variable,
and the scalar iteration counts. Initial and final curves are highlighted
in blue and red, respectively.

\paragraph{\textbf{Example 4: Evolution of an elongated tube}}
The initial curve is an elongated tube of dimensions $8\times1$.
We use $\Delta t=10^{-4}$, with $N=128$ for $\lambda=0$ and $N=256$
for $\lambda=0.5$ and $\lambda=2$.
For $\lambda=0$, the curve develops a nonconvex intermediate shape before
approaching an expanding circular configuration in Fig.~\ref{exp3}.
The modified discrete energy decreases monotonically throughout the evolution, consistent with the theoretical results. Moreover, the mesh ratio gradually approaches one, indicating an asymptotically uniform distribution of
mesh points.  For $\lambda=0.5$ and $\lambda=2$, the evolving curves approach
the corresponding circular equilibria, as shown in Figs.~\ref{exp4} and~\ref{exp5}, respectively.


\paragraph{\textbf{Example 5: Evolution of a 2:1 lemniscate}}The initial immersed curve is parameterized by
\begin{equation}
	\boldsymbol X(\rho,0)=
	\left(\frac{a(\rho)\cos(2\pi\rho)}{1+\sin^2(2\pi\rho)},
	\frac{a(\rho)\cos(2\pi\rho)\sin(2\pi\rho)}{1+\sin^2(2\pi\rho)}\right)^T,
	\quad
	a(\rho)=\begin{cases}1,&\rho\in[1/4,3/4],\\2,&\text{otherwise}.
	\end{cases}
	\label{eq:initial_lemniscate}
\end{equation}
We set $N=256$ and $\Delta t=10^{-5}$ and consider
$\lambda=0$ and $\lambda=0.4$.
As shown in Figs.~\ref{exp6} and~\ref{exp7}, the initially
asymmetric lemniscate gradually approaches a symmetric
configuration. For $\lambda=0$, the curve continues to expand,
whereas for $\lambda=0.4$, the length penalty limits the
expansion and the curve approaches an approximately
stationary configuration.

\paragraph{\textbf{Example 6: Time-step sensitivity for a three-fold curve}}
We consider the initial curve
\begin{equation}
	\boldsymbol X(\rho,0)=
	\left(1+\frac1{15}\cos(6\pi\rho)\right)
	\begin{pmatrix}
		\cos(2\pi\rho)\\
		\sin(2\pi\rho)
	\end{pmatrix},
	\qquad \rho\in\mathbb I,
	\label{eq:initial_threefold}
\end{equation}
with $N=256$, $\lambda=0$, and $T=1$.
We apply the first-order scheme with
$\Delta t=10^{-4},10^{-3},5\times10^{-3}$
to investigate the influence of the time-step size.

\begin{figure}[htbp]
	\centering
	\includegraphics[width=0.79\linewidth,height=4cm]{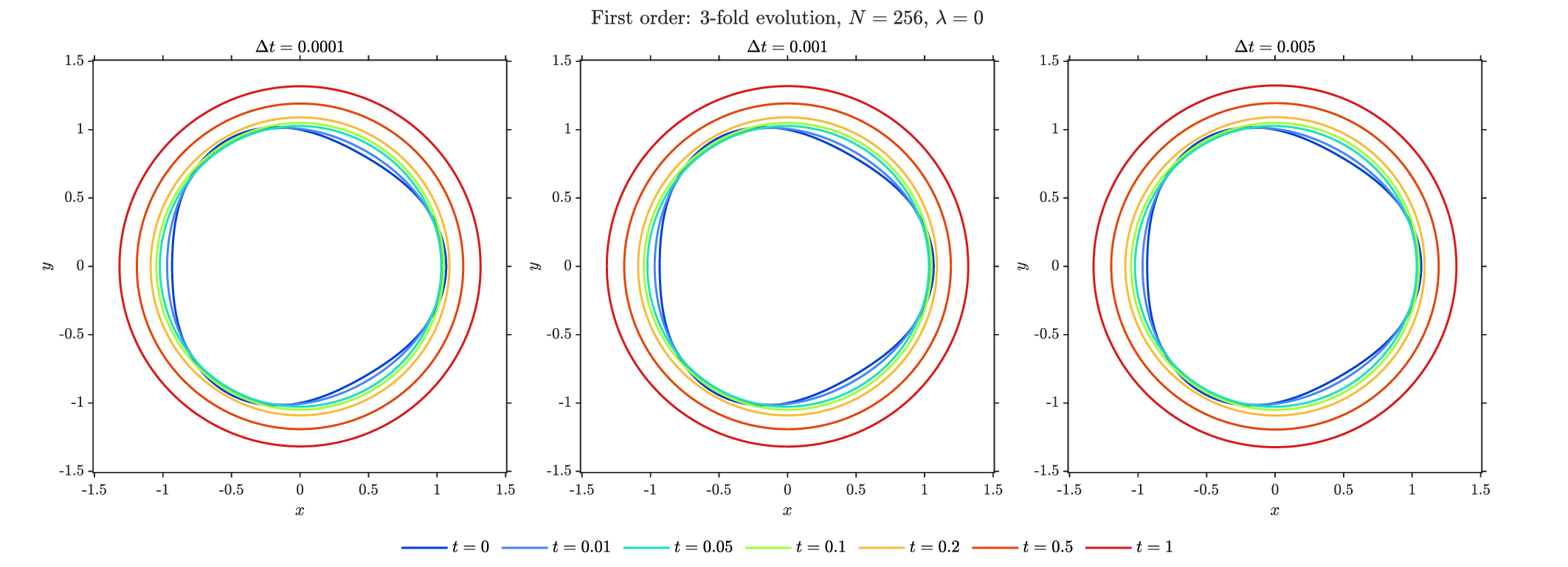}\\
	\includegraphics[width=0.30\linewidth]{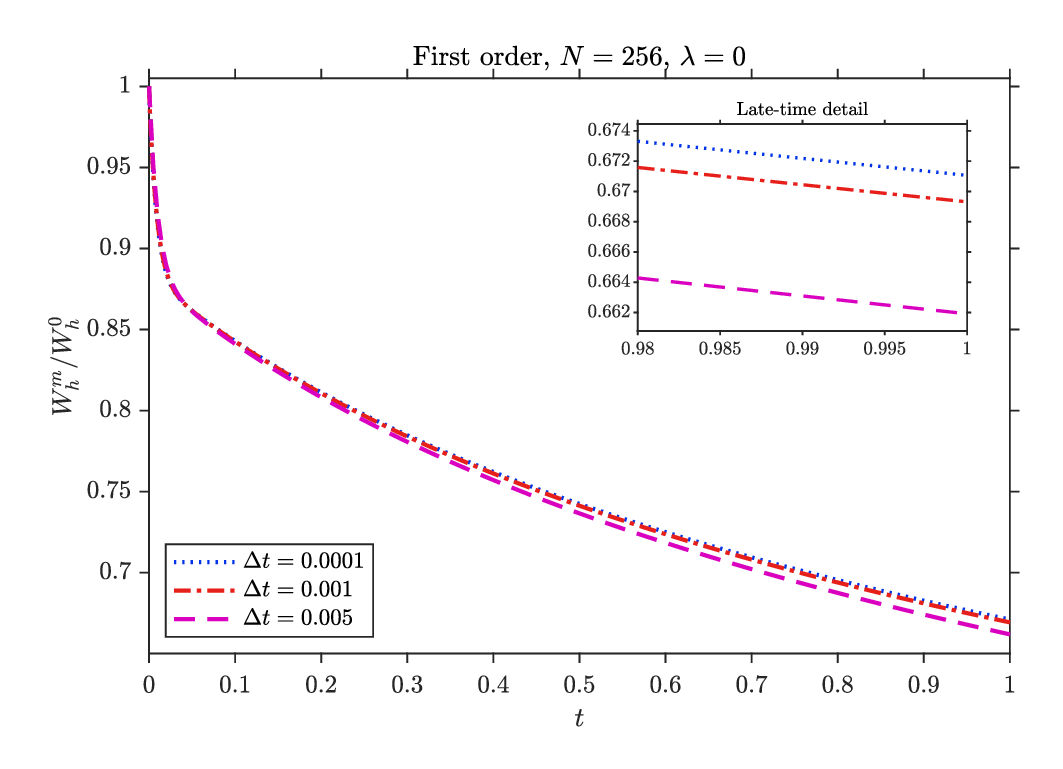}
	\includegraphics[width=0.30\linewidth]{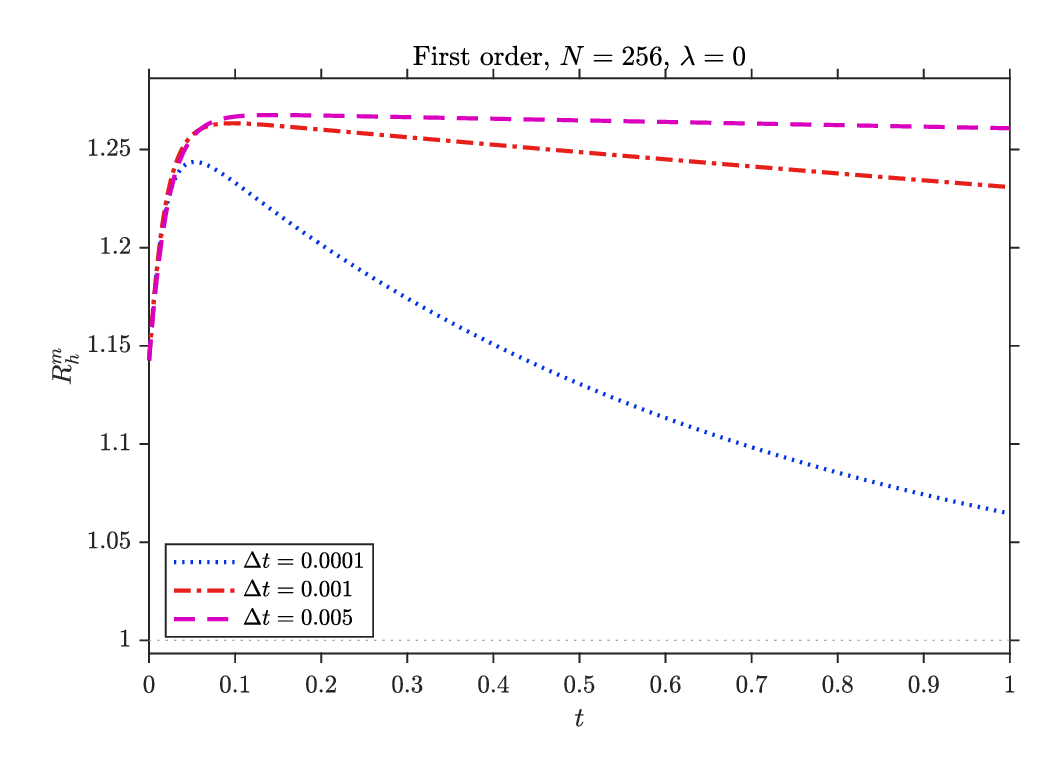}
	\caption{Example 6: time-step sensitivity of the first-order scheme for a three-fold curve with $\lambda=0$ and $N=256$.
		Top, from left to right: curve evolution for $\Delta t=10^{-4},10^{-3},5\times10^{-3}$, with snapshots at $t=0,0.01,0.05,0.1,0.2,0.5,1$.
		Bottom: normalized modified energy (left) and mesh ratio $R_h^m$ (right).}
	\label{f8}
\end{figure}

The curve evolution, normalized modified energy, and mesh
ratio are presented in Fig.~\ref{f8}.
The initial three-fold perturbation gradually diminishes,
and the curve approaches an expanding circular configuration.
The modified energy decreases monotonically for all three
time steps, while the inset reveals a greater cumulative
energy decrease for larger time steps.
The mesh-ratio histories show a more pronounced dependence
on $\Delta t$. After an initial increase, the mesh ratio
decreases more rapidly for smaller time steps, indicating
a faster approach toward a uniform nodal distribution.

\section{Conclusions}\label{part6}
In this paper, we have developed first- and second-order energy-stable
parametric finite element methods for planar Willmore flow.
By introducing a global auxiliary variable, we obtained an
equivalent reformulation that preserves the continuous
energy-dissipation structure. Based on this reformulation,
we constructed a first-order scheme and a second-order
predictor--corrector Crank--Nicolson scheme, both satisfying
unconditional stability estimates for a modified discrete
energy.  Under suitable  assumptions, we established the unique
solvability of the linear subproblems and derived computable
sufficient conditions for local existence and uniqueness
of the global auxiliary solution near unity.
Each stage requires only two linear solves and one scalar nonlinear solve.
Numerical experiments exhibit convergence behavior consistent
with the expected temporal orders and second-order spatial
accuracy, together with energy decay and good mesh quality. The computational comparisons further
illustrate the efficiency gained by using larger time steps
with the second-order scheme. Future work will focus on rigorous convergence and error	analysis for both schemes, including the control of the
auxiliary-variable error over time. Extensions of the proposed
framework to more general geometric flows will also be investigated.


\bibliographystyle{siamplain}
\bibliography{Willmore_flow}

@article{droske2004level,
	author  = {Droske, Marc and Rumpf, Martin},
	title   = {A level set formulation for {Willmore} flow},
	journal = {Interfaces and Free Boundaries},
	volume  = {6},
	pages   = {361--378},
	year    = {2004},
	doi = {10.4171/IFB/105}
}

@article{du2004phase,
	author  = {Du, Qiang and Liu, Chun and Wang, Xiaoqiang},
	title   = {A phase field approach in the numerical study of the
	elastic bending energy for vesicle membranes},
	journal = {Journal of Computational Physics},
	volume  = {198},
	number  = {2},
	pages   = {450--468},
	year    = {2004},
	doi = {10.1016/j.jcp.2004.01.029}
}

@article{franken2013phase,
	author  = {Franken, Martina and Rumpf, Martin and Wirth, Benedikt},
	title   = {A phase field based {PDE} constrained optimization approach
	to time discrete {Willmore} flow},
	journal = {International Journal of Numerical Analysis and Modeling},
	volume  = {10},
	number  = {1},
	pages   = {116--138},
	year    = {2013}
}

@article{bretin2015phasefield,
	author  = {Bretin, {\'E}lie and Masnou, Simon and Oudet, {\'E}douard},
	title   = {Phase-field approximations of the {Willmore} functional and flow},
	journal = {Numerische Mathematik},
	volume  = {131},
	number  = {1},
	pages   = {115--171},
	year    = {2015},
	doi = {10.1007/s00211-014-0683-4}
}

@article{rumpf2025hybrid,
	author  = {Rumpf, Martin and Sassen, Josua and Smoch, Christoph},
	title   = {A hybrid minimizing movement and neural network approach to {Willmore} flow},
	journal = {Interfaces and Free Boundaries},
	year    = {2026},
	note    = {Published online first},
	doi     = {10.4171/IFB/563}
}

@article{bgn2016spontaneous,
	author  = {Barrett, John W. and Garcke, Harald and N{\"u}rnberg, Robert},
	title   = {Computational parametric {Willmore} flow with spontaneous
	curvature and area difference elasticity effects},
	journal = {SIAM Journal on Numerical Analysis},
	volume  = {54},
	number  = {3},
	pages   = {1732--1762},
	year    = {2016},
	doi = {10.1137/16M1065379}
}

@article{bgn2017boundary,
	author  = {Barrett, John W. and Garcke, Harald and N{\"u}rnberg, Robert},
	title   = {Stable variational approximations of boundary value problems
	for {Willmore} flow with {Gaussian} curvature},
	journal = {IMA Journal of Numerical Analysis},
	volume  = {37},
	number  = {4},
	pages   = {1657--1709},
	year    = {2017},
	doi = {10.1093/imanum/drx006}
}

@article{dziuk2008computational,
	author  = {Dziuk, Gerhard},
	title   = {Computational parametric {Willmore} flow},
	journal = {Numerische Mathematik},
	volume  = {111},
	number  = {1},
	pages   = {55--80},
	year    = {2008},
	doi = {10.1007/s00211-008-0179-1}
}

@article{gao2026mdr,
	author  = {Gao, Guangwei and Garcke, Harald and Li, Buyang and Tang, Rong},
	title   = {An energy-stable minimal deformation rate scheme
	for mean curvature flow and surface diffusion},
	journal = {SIAM Journal on Scientific Computing},
	volume  = {48},
	number  = {1},
	pages   = {A103--A131},
	year    = {2026}
}

@article{jiang2026structure,
	author  = {Jiang, Wei and Su, Chunmei and Tang, Kaike},
	title   = {Structure-preserving parametric finite element methods for elastic flows of planar curves with an explicit tangential velocity},
	journal = {Journal of Computational Physics},
	volume  = {563},
	pages   = {115087},
	year    = {2026},
	doi     = {10.1016/j.jcp.2026.115087}
}

@article{jiang2024bdf,
	author  = {Jiang, Wei and Su, Chunmei and Zhang, Ganghui},
	title   = {Stable backward differentiation formula time discretization
	of {BGN}-based parametric finite element methods
	for geometric flows},
	journal = {SIAM Journal on Scientific Computing},
	volume  = {46},
	number  = {5},
	pages   = {A2874--A2898},
	year    = {2024},
	doi = {10.1137/23M1625597}
}

@misc{jiang2026imex,
	author        = {Jiang, Wei and Su, Chunmei and Tang, Kaike},
	title         = {Implict-Explicit {Runge--Kutta} schemes for geometric flows},
	year          = {2026},
	eprint        = {2610.08019},
	archivePrefix = {arXiv},
	primaryClass  = {math.NA}
}

@misc{zhang2026arbitrary,
	author        = {Zhang, Ganghui and Andrews, Boris D. and Farrell, Patrick E.},
	title         = {Arbitrary-order structure-preserving discretizations
	for geometric curvature flows},
	year          = {2026},
	eprint        = {2605.20371},
	archivePrefix = {arXiv},
	primaryClass  = {math.NA}
}

@misc{bao2026pullback,
	author        = {Bao, Weizhu and Li, Yifei and Wang, Dongmin},
	title         = {Arbitrarily high-order structure-preserving parametric
	finite element methods for geometric {PDEs}
	via local pullback flow maps},
	year          = {2026},
	eprint        = {2610.04389},
	archivePrefix = {arXiv},
	primaryClass  = {math.NA}
}

@article{kemmochi2025tangential,
	author  = {Kemmochi, Tomoya and Miyatake, Yuto and Sakakibara, Koya},
	title   = {Structure-preserving numerical methods for constrained gradient
	flows of planar closed curves with explicit tangential velocities},
	journal = {Japan Journal of Industrial and Applied Mathematics},
	volume  = {42},
	number  = {2},
	pages   = {575--603},
	year    = {2025},
	doi = {10.1007/s13160-024-00684-0}
}

@article{helfrich1973elastic,
	author  = {Helfrich, Wolfgang},
	title   = {Elastic properties of lipid bilayers: Theory and possible experiments},
	journal = {Zeitschrift f{\"u}r Naturforschung C},
	volume  = {28},
	number  = {11--12},
	pages   = {693--703},
	year    = {1973}
}

@article{clarenz2004restoration,
	author  = {Clarenz, U. and Diewald, U. and Dziuk, G. and Rumpf, M. and Rusu, R.},
	title   = {A finite element method for surface restoration with smooth boundary conditions},
	journal = {Computer Aided Geometric Design},
	volume  = {21},
	number  = {5},
	pages   = {427--445},
	year    = {2004},
	doi = {10.1016/j.cagd.2004.02.004}
}

@article{torabi2009phase,
	author  = {Torabi, Solmaz and Lowengrub, John and Voigt, Axel and Wise, Steven},
	title   = {A new phase-field model for strongly anisotropic systems},
	journal = {Proceedings of the Royal Society A},
	volume  = {465},
	pages   = {1337--1359},
	year    = {2009},
	doi = {10.1098/rspa.2008.0385}
}

@book{willmore1993riemannian,
	title={Riemannian geometry},
	author={Willmore, Thomas James},
	year={1993},
	publisher={Oxford University Press}
}

@techreport{levien2008elastica,
	author      = {Levien, Raph},
	title       = {The elastica: a mathematical history},
	institution = {EECS Department, University of California, Berkeley},
	number      = {UCB/EECS-2008-103},
	year        = {2008},
	url         = {https://www2.eecs.berkeley.edu/Pubs/TechRpts/2008/EECS-2008-103.html}
}

@article{pan2026lagrange,
	author  = {Pan, Zhiqing and Jia, Jiwei},
	title   = {An Energy-Stable Parametric Finite Element Framework with {Lagrange} Multiplier for the Planar {Willmore} Flow},
	journal = {Journal of Scientific Computing},
	volume  = {108},
	pages   = {77},
	year    = {2026},
	doi = {10.1007/s10915-026-03399-1}
}

@article{JIANG2024113220,
	author  = {Jiang, Wei and Su, Chunmei and Zhang, Ganghui},
	title   = {A second-order in time, {BGN}-based parametric finite element method for geometric flows of curves},
	journal = {Journal of Computational Physics},
	volume  = {514},
	pages   = {113220},
	year    = {2024}
}

@article{jiang2025predictor,
	author  = {Jiang, Wei and Su, Chunmei and Zhang, Ganghui and Zhang, Lian},
	title   = {Predictor-corrector, {BGN}-based parametric finite element methods for surface diffusion},
	journal = {Journal of Computational Physics},
	volume  = {530},
	pages   = {113901},
	year    = {2025}
}

@article{li2025sinum,
	author  = {Bao, Weizhu and Li, Yifei},
	title   = {An Energy-Stable Parametric Finite Element Method for the Planar {Willmore} Flow},
	journal = {SIAM Journal on Numerical Analysis},
	volume  = {63},
	number  = {1},
	pages   = {103--121},
	year    = {2025},
	doi = {10.1137/24M1633893}
}

@article{garcke2025stable,
	title={Stable fully discrete finite element methods with {BGN} tangential motion for {Willmore} flow of planar curves},
	author={Garcke, Harald and N{\"u}rnberg, Robert and Zhao, Quan},
	journal={Journal of Scientific Computing},
	volume={105},
	number={2},
	pages={45},
	year={2025},
	publisher={Springer},
	doi = {10.1007/s10915-025-03068-9}
}

@article{BARRETT2007441,
	title = {A parametric finite element method for fourth order geometric evolution equations},
	journal = {Journal of Computational Physics},
	volume = {222},
	number = {1},
	pages = {441--467},
	year = {2007},
	issn = {0021-9991},
	author = {John W. Barrett and Harald Garcke and Robert N{\"u}rnberg},
	doi = {10.1016/j.jcp.2006.07.026}
}

@article{bgn2008willmore,
	author = {Barrett, John W. and Garcke, Harald and N{\"u}rnberg, Robert},
	title = {Parametric Approximation of {Willmore} Flow and Related Geometric Evolution Equations},
	journal = {SIAM Journal on Scientific Computing},
	volume = {31},
	number = {1},
	pages = {225--253},
	year = {2008},
	doi = {10.1137/070700231}
}

@article{bgn2012elastic,
	author = {Barrett, John W. and Garcke, Harald and N{\"u}rnberg, Robert},
	title = {Parametric Approximation of Isotropic and Anisotropic Elastic Flow for Closed and Open Curves},
	journal = {Numerische Mathematik},
	volume = {120},
	pages = {489--542},
	year = {2012},
	doi = {10.1007/s00211-011-0416-x}
}

@article{deckelnick2009elastic,
	author = {Deckelnick, Klaus and Dziuk, Gerhard},
	title = {Error Analysis for the Elastic Flow of Parametrized Curves},
	journal = {Mathematics of Computation},
	volume = {78},
	number = {266},
	pages = {645--671},
	year = {2009},
	doi = {10.1090/S0025-5718-08-02176-5}
}

@article{gnz2026splitting,
	author = {Garcke, Harald and N{\"u}rnberg, Robert and Zhao, Quan},
	title = {An Energy-Stable Parametric Finite Element Method for {Willmore} Flow with Normal-Tangential Velocity Splitting},
	journal = {SIAM Journal on Scientific Computing},
	volume = {48},
	number = {3},
	pages = {A1235--A1259},
	year = {2026},
	doi = {10.1137/25M1773878}
}

@article{yang2021zec,
	author = {Yang, Xiaofeng},
	title = {On a Novel Fully Decoupled, Second-Order Accurate Energy Stable Numerical Scheme for a Binary Fluid-Surfactant Phase-Field Model},
	journal = {SIAM Journal on Scientific Computing},
	volume = {43},
	number = {2},
	pages = {B479--B507},
	year = {2021}
}

@article{bonito2010biomembranes,
	author = {Bonito, Andrea and Nochetto, Ricardo H. and Pauletti, M. Sebastian},
	title = {Parametric {FEM} for Geometric Biomembranes},
	journal = {Journal of Computational Physics},
	volume = {229},
	number = {9},
	pages = {3171--3188},
	year = {2010}
}

@article{xu2006surface,
	author = {Xu, Guoliang and Pan, Qing},
	title = {{$G^1$} Surface Modelling Using Fourth Order Geometric Flows},
	journal = {Computer-Aided Design},
	volume = {38},
	number = {4},
	pages = {392--403},
	year = {2006}
}

@article{chen2018willmore,
	author = {Chen, Ying and Lowengrub, John and Shen, Jie and Wang, Cheng and Wise, Steven},
	title = {Efficient Energy Stable Schemes for Isotropic and Strongly Anisotropic {Cahn--Hilliard} Systems with the {Willmore} Regularization},
	journal = {Journal of Computational Physics},
	volume = {365},
	pages = {56--73},
	year = {2018},
	doi = {10.1016/j.jcp.2018.03.024}
}

\newpage

\appendix
\section{First-order auxiliary equation}
\label{app:quartic_coefficients}

Substituting the affine decomposition \eqref{eq:first_order_decomposition} into \eqref{ful5} gives
\begin{equation}
	F_m(q)=\sum_{k=0}^{4}a_{k,m}q^k=0,\qquad q=Q^{m+1}.
	\label{eq:quartic_polynomial}
\end{equation}
The coefficients $a_{k,m}$ depend only on the subproblem solutions $\boldsymbol X_\alpha^{m+1}$, $\varkappa_\alpha^{m+1}$, $\alpha=1,2$, and data from previous time levels. Collecting powers of $q$ gives
\begin{align}
	a_{0,m}
	&=
	\frac{1}{2\Delta t}\left(\partial_s(\boldsymbol X_1^{m+1}-\boldsymbol X^m)\cdot\partial_s\boldsymbol X_1^{m+1}(\varkappa_1^{m+1})^2,1\right)^h_{\Gamma^m},
	\label{eq:a0}
	\\[-2pt]
	a_{1,m}
	&=
	-\frac{Q^m}{\Delta t}
	-\frac12\left(\partial_s\frac{\boldsymbol X^m-\boldsymbol X^{m-1}}{\Delta t}\cdot\partial_s\boldsymbol X^m\,\varkappa^m,\varkappa_1^{m+1}\right)^h_{\Gamma^m}
	\nonumber\\[-2pt]
	&\quad
	+\frac{1}{2\Delta t}\left(\partial_s(\boldsymbol X_1^{m+1}-\boldsymbol X^m)\cdot\partial_s\boldsymbol X_2^{m+1}(\varkappa_1^{m+1})^2,1\right)^h_{\Gamma^m}
	\nonumber\\[-2pt]
	&\quad
	+\frac{1}{2\Delta t}\left(\partial_s\boldsymbol X_2^{m+1}\cdot\partial_s\boldsymbol X_1^{m+1}(\varkappa_1^{m+1})^2,1\right)^h_{\Gamma^m}
	+\frac{1}{\Delta t}\left(\partial_s(\boldsymbol X_1^{m+1}-\boldsymbol X^m)\cdot\partial_s\boldsymbol X_1^{m+1}\varkappa_1^{m+1}\varkappa_2^{m+1},1\right)^h_{\Gamma^m},
	\label{eq:a1}
	\\[-2pt]
	a_{2,m}
	&=
	\frac{1}{\Delta t}
	-\frac12\left(\partial_s\frac{\boldsymbol X^m-\boldsymbol X^{m-1}}{\Delta t}\cdot\partial_s\boldsymbol X^m\,\varkappa^m,\varkappa_2^{m+1}\right)^h_{\Gamma^m}
	+\frac{1}{2\Delta t}\left(|\partial_s\boldsymbol X_2^{m+1}|^2(\varkappa_1^{m+1})^2,1\right)^h_{\Gamma^m}
	\nonumber\\[-2pt]
	&\quad
	+\frac{1}{\Delta t}\left(\partial_s(\boldsymbol X_1^{m+1}-\boldsymbol X^m)\cdot\partial_s\boldsymbol X_2^{m+1}\varkappa_1^{m+1}\varkappa_2^{m+1},1\right)^h_{\Gamma^m}
	\nonumber\\[-2pt]
	&\quad
	+\frac{1}{\Delta t}\left(\partial_s\boldsymbol X_2^{m+1}\cdot\partial_s\boldsymbol X_1^{m+1}\varkappa_1^{m+1}\varkappa_2^{m+1},1\right)^h_{\Gamma^m}
	+\frac{1}{2\Delta t}\left(\partial_s(\boldsymbol X_1^{m+1}-\boldsymbol X^m)\cdot\partial_s\boldsymbol X_1^{m+1}(\varkappa_2^{m+1})^2,1\right)^h_{\Gamma^m},
	\label{eq:a2}
	\\[-2pt]
	a_{3,m}
	&=
	\frac{1}{\Delta t}\left(|\partial_s\boldsymbol X_2^{m+1}|^2\varkappa_1^{m+1}\varkappa_2^{m+1},1\right)^h_{\Gamma^m}
	+\frac{1}{2\Delta t}\left(\partial_s(\boldsymbol X_1^{m+1}-\boldsymbol X^m)\cdot\partial_s\boldsymbol X_2^{m+1}(\varkappa_2^{m+1})^2,1\right)^h_{\Gamma^m}
	\nonumber\\[-2pt]
	&\quad
	+\frac{1}{2\Delta t}\left(\partial_s\boldsymbol X_2^{m+1}\cdot\partial_s\boldsymbol X_1^{m+1}(\varkappa_2^{m+1})^2,1\right)^h_{\Gamma^m},
	\label{eq:a3}
	\\[-2pt]
	a_{4,m}
	&=
	\frac{1}{2\Delta t}\left(|\partial_s\boldsymbol X_2^{m+1}|^2(\varkappa_2^{m+1})^2,1\right)^h_{\Gamma^m}
	\geq0.
	\label{eq:a4}
\end{align}

\section{Second-order auxiliary equation}
\label{app:second_auxiliary}

With $q=Q^{m+1}$, write
$\boldsymbol X(q)=\boldsymbol X_1^{m+1}+q\boldsymbol X_2^{m+1}$,
$\varkappa(q)=\varkappa_1^{m+1}+q\varkappa_2^{m+1}$.
Set
$Q^{m+\frac12}(q)=\frac{q+Q^m}{2}$,
$\varkappa^{m+\frac12}(q)=\frac{\varkappa(q)+\varkappa^m}{2}$,
and
$\boldsymbol X^{m+\frac12}(q)=\frac{\boldsymbol X(q)+\boldsymbol X^m}{2}$.

Write
$(u,v)^h_{\widetilde\Gamma}:=(u,v)^h_{\widetilde\Gamma^{m+\frac12}}$,
and define
$H_m^{(2)}=\widetilde\partial_s\widetilde{\boldsymbol U}^{m+\frac12}\cdot\widetilde\partial_s\widetilde{\boldsymbol X}^{m+\frac12}\,\widetilde\varkappa^{m+\frac12}$.

The discrete bending-energy and length identities give
\begin{align}
	F_m^{(2)}(q)
	={}&
	\frac{q^2-(Q^m)^2}{2\Delta t}
	-\frac{q+Q^m}{4}\left(H_m^{(2)},\varkappa^{m+\frac12}(q)\right)^h_{\widetilde\Gamma}
	-\left(\frac{\varkappa(q)-\varkappa^m}{\Delta t},\varkappa^{m+\frac12}(q)\right)^h_{\widetilde\Gamma}
	\nonumber\\
	&+
	\frac{1}{2\Delta t}\left[\|\varkappa(q)\|_{\Gamma[\boldsymbol X(q)],h}^2-\|\varkappa^m\|_{\Gamma^m,h}^2\right]
	+\frac{\lambda}{\Delta t}\left[|\Gamma[\boldsymbol X(q)]|-|\Gamma^m|\right]
	\nonumber\\
	&-
	\lambda\left(\widetilde\partial_s\boldsymbol X^{m+\frac12}(q),\widetilde\partial_s\frac{\boldsymbol X(q)-\boldsymbol X^m}{\Delta t}\right)^h_{\widetilde\Gamma}.
	\label{eq:app_second_F}
\end{align}

Expanding the terms on the predicted curve gives
$F_m^{(2)}(q)=a_{0,m}^{(2)}+a_{1,m}^{(2)}q+a_{2,m}^{(2)}q^2+\mathcal R_m(q)$.
In \eqref{eq:second_F_decomposition}, the coefficients are
\begin{align}
	a_{0,m}^{(2)}
	={}&
	-\frac{(Q^m)^2}{2\Delta t}
	-\frac{Q^m}{8}\left(H_m^{(2)},\varkappa_1^{m+1}+\varkappa^m\right)^h_{\widetilde\Gamma}
	-\frac{1}{2\Delta t}\left[\|\varkappa_1^{m+1}\|_{\widetilde\Gamma,h}^2-\|\varkappa^m\|_{\widetilde\Gamma,h}^2\right]
	\nonumber\\
	&-
	\frac{1}{2\Delta t}\|\varkappa^m\|_{\Gamma^m,h}^2
	-\frac{\lambda}{\Delta t}|\Gamma^m|
	-\frac{\lambda}{2\Delta t}\left[\|\widetilde\partial_s\boldsymbol X_1^{m+1}\|_{\widetilde\Gamma}^2-\|\widetilde\partial_s\boldsymbol X^m\|_{\widetilde\Gamma}^2\right],
	\label{eq:app_second_a0}
	\\
	a_{1,m}^{(2)}
	={}&
	-\frac18\left(H_m^{(2)},\varkappa_1^{m+1}+\varkappa^m+Q^m\varkappa_2^{m+1}\right)^h_{\widetilde\Gamma}
	-\frac1{\Delta t}(\varkappa_1^{m+1},\varkappa_2^{m+1})^h_{\widetilde\Gamma}
	-\frac{\lambda}{\Delta t}\left(\widetilde\partial_s\boldsymbol X_1^{m+1},\widetilde\partial_s\boldsymbol X_2^{m+1}\right)^h_{\widetilde\Gamma},
	\label{eq:app_second_a1}
	\\
	a_{2,m}^{(2)}
	={}&
	\frac1{2\Delta t}
	-\frac18(H_m^{(2)},\varkappa_2^{m+1})^h_{\widetilde\Gamma}
	-\frac1{2\Delta t}\|\varkappa_2^{m+1}\|_{\widetilde\Gamma,h}^2
	-\frac{\lambda}{2\Delta t}\|\widetilde\partial_s\boldsymbol X_2^{m+1}\|_{\widetilde\Gamma}^2.
	\label{eq:app_second_a2}
\end{align}

The remainder is
\begin{equation}
	\mathcal R_m(q)=\frac{1}{2\Delta t}\|\varkappa(q)\|_{\Gamma[\boldsymbol X(q)],h}^2+\frac{\lambda}{\Delta t}|\Gamma[\boldsymbol X(q)]|.
	\label{eq:app_second_R}
\end{equation}

For edge $j$, set
$\boldsymbol h_{\alpha,j}^{m+1}=\boldsymbol X_\alpha^{m+1}(\rho_j)-\boldsymbol X_\alpha^{m+1}(\rho_{j-1})$, $\alpha=1,2$,
and
$\boldsymbol d_j(q)=\boldsymbol h_{1,j}^{m+1}+q\boldsymbol h_{2,j}^{m+1}$, $\ell_j(q)=|\boldsymbol d_j(q)|$.
Define
\[
A_j=(\varkappa_{1,j-1}^{m+1})^2+(\varkappa_{1,j}^{m+1})^2,\quad
B_j=2\left(\varkappa_{1,j-1}^{m+1}\varkappa_{2,j-1}^{m+1}+\varkappa_{1,j}^{m+1}\varkappa_{2,j}^{m+1}\right),\quad
C_j=(\varkappa_{2,j-1}^{m+1})^2+(\varkappa_{2,j}^{m+1})^2.
\]
The remainder then takes the edgewise form
\begin{equation}
	\mathcal R_m(q)=\frac1{4\Delta t}\sum_{j=1}^{N}\ell_j(q)\left(A_j+B_jq+C_jq^2+4\lambda\right).
	\label{eq:app_second_R_explicit}
\end{equation}

For nondegenerate edges, differentiating $\ell_j(q)=|\boldsymbol h_{1,j}^{m+1}+q\boldsymbol h_{2,j}^{m+1}|$ gives
\begin{align}
	\ell_j'(q)=\frac{\boldsymbol d_j(q)\cdot\boldsymbol h_{2,j}^{m+1}}{\ell_j(q)},~
	\ell_j''(q)=\frac{|\boldsymbol h_{1,j}^{m+1}|^2|\boldsymbol h_{2,j}^{m+1}|^2-(\boldsymbol h_{1,j}^{m+1}\cdot\boldsymbol h_{2,j}^{m+1})^2}{\ell_j(q)^3}.
	\label{eq:app_edge_second}
\end{align}

Differentiating \eqref{eq:app_second_R_explicit} yields
\begin{align}
	\mathcal R_m'(q)
	={}&
	\frac1{4\Delta t}\sum_{j=1}^{N}\left[\ell_j'(q)(A_j+B_jq+C_jq^2+4\lambda)+\ell_j(q)(B_j+2C_jq)\right],
	\label{eq:app_R_first}
	\\
	\mathcal R_m''(q)
	={}&
	\frac1{4\Delta t}\sum_{j=1}^{N}\left[\ell_j''(q)(A_j+B_jq+C_jq^2+4\lambda)+2\ell_j'(q)(B_j+2C_jq)+2C_j\ell_j(q)\right].
	\label{eq:app_R_second}
\end{align}

Therefore,
\begin{equation}
	(F_m^{(2)})'(1)=a_{1,m}^{(2)}+2a_{2,m}^{(2)}+\mathcal R_m'(1),\qquad
	(F_m^{(2)})''(q)=2a_{2,m}^{(2)}+\mathcal R_m''(q).
	\label{eq:app_F_derivatives}
\end{equation}

These derivatives can be used to check the local solvability condition \eqref{eq:Q_second_condition}.

\end{document}